\documentclass[a4paper,10pt]{amsart}
\usepackage{amsfonts}
\usepackage{amsmath}
\usepackage{amstext}
\usepackage{amsxtra}
\usepackage{amssymb}
\usepackage{epsfig}
\usepackage{latexsym}

\usepackage{a4wide}

\usepackage{palatino}

\numberwithin{equation}{section}

\theoremstyle{plain}
\newtheorem{teo}{Theorem}[section]
\newtheorem{lemma}[teo]{Lemma}
\newtheorem{prop}[teo]{Proposition}
\newtheorem{cor}[teo]{Corollary}
\newtheorem{ackn}{Acknowledgement\!}

\theoremstyle{definition}
\newtheorem{dfnz}[teo]{Definition}
\newtheorem{conge}[teo]{Conjecture}

\theoremstyle{remark}
\newtheorem{rem}[teo]{Remark}

\def\R{\mathbb R}

\def\HH{\mathcal H}
\def\NN{\mathbb N}
\def\composed{\circ}
\def\eps{\varepsilon}

\def\dert{\partial_t}
\def\ders{\partial_s}
\def\pol{\mathfrak p}
\def\qol{\mathfrak q}
\def\rol{\mathfrak r}
\def\loc{_{\operatorname{loc}}}
\def\div{\operatornamewithlimits{div}\nolimits}
\def\limnf{\operatornamewithlimits{\underline{\lim}}}
\def\limup{\operatornamewithlimits{\overline{\lim}}}

\def\composed{\circ}

\def\eps{\varepsilon}

\def\TTT{\mathbb T}
\def\VVV{\mathbb V}

\def\tt{\mathfrak t}
\def\zz{\tt}
\def\dist{\mathrm {dist}}
\def\res  {
  \begin{picture}(9,8)
    \put (1,0){\line(0,1){8}}
    \put (1,0){\line(1,0){5}}
  \end{picture}}

\begin{document}

\title[Motion by Curvature of Planar Networks]
{Motion by Curvature of Planar Networks}

\author[Carlo Mantegazza]{Carlo Mantegazza}
\address[Carlo Mantegazza]{Scuola
  Normale Superiore, Pisa, 56126, Italy}
\email[C.~Mantegazza]{mantegazza@sns.it}

\author[Matteo Novaga]{Matteo Novaga}
\address[Matteo Novaga]{Math.
  Dept. Univ. Pisa, Pisa, 56127, Italy}
\email[M.~Novaga]{novaga@dm.unipi.it}

\author[Vincenzo Maria Tortorelli]{Vincenzo
  Maria Tortorelli}
\address[Vincenzo Maria Tortorelli]{Math.
  Dept. Univ. Pisa, Pisa, 56127, Italy}
\email[V.~M.~Tortorelli]{tortorel@dm.unipi.it}

\date{\today}

\begin{abstract} We consider the motion by curvature of a network of
  smooth curves with multiple junctions in the plane, that is, the
  geometric gradient flow associated to the length
  functional.\\
  Such a flow represents the evolution of a two--dimensional
  multiphase system where the energy is simply the sum of the
  lengths of the interfaces, in particular it is a possible model for
  the  growth of grain boundaries.\\
  Moreover, the motion of these networks of curves is the simplest
  example of curvature flow for sets which are ``essentially'' non
  regular.
  
  As a first step, in this paper we study in detail the case of three curves in the
  plane concurring at a single triple junction and with the other ends
  fixed. We show some results about the existence, uniqueness and,
  in particular, the global regularity of the flow, following
  the line of analysis carried on in the last years for the
  evolution by mean curvature of smooth curves and hypersurfaces.
\end{abstract}

\maketitle

\tableofcontents

\section{Introduction and Basic Definitions}

In this work we address the problem of the motion by curvature
of a network of curves in the plane, where by network of
curves we mean a connected planar graph without
self--intersections.\\
The evolution by curvature of such a network is the geometric
gradient flow with respect to the energy given by
the {\em Length functional} which is simply the sum of the lengths of all the
curves of the network (see~\cite{brakke}).

We point out two motivations to study this evolution. The first is 
the analysis of models of two--dimensional multiphase
systems, where the problem of the structure and regularity of the
interfaces between different phases arises naturally.
As an example, the model where the energy of a configuration is
simply the total length has proven useful in the analysis of the
growth of grain boundaries,
see~\cite{brakke,bronsard,gurtin2,kinderliu}, the papers by Herring and Mullins
in~\cite{hermul} and {\em http://mimp.mems.cmu.edu}.\\
The second motivation is more theoretical: the evolution of such a network of
curves in the plane is the simplest example
of motion by mean curvature of a set which is {\em essentially}
singular. In literature there are 
various generalized definitions of flow by mean curvature
for a singular set (see~\cite{altawa,brakke,degio4,es,ilman1,soner1},
for instance). All of them are fairly general, but usually lack
uniqueness and a satisfactory regularity theory, even for simple
situations.

Inspired by Grayson's Theorem in~\cite{gray1}, stating 
that any smooth closed curve embedded in $\R^2$ evolves by curvature without 
singularities before vanishing, and by the new approach to such result by Huisken
in~\cite{huisk2}, one can reasonably expect that an ``embedded'' network
of smooth curves does not develop singularities during the
flow if its ``topological structure'' does not change (we will be more precise
about this point in the sequel) and asymptotically converges to a critical
configuration for the Length functional.\\
Moreover, in~\cite{huisk2} it is also shown that the motion by curvature of a
single embedded curve in a strip of $\R^2$,
with its end points fixed to be some $P^1$ and $P^2$ on the
boundary of the strip, evolves smoothly and approaches 
the segment connecting the two points $P^1, P^2$. We can see this case as a very
special positive example of motion of a network and we will try to
follow the line of analysis traced in that paper.

We consider a connected network ${\mathbb S}=\cup_{i=1}^n \sigma^i$,
composed of a finite family of smooth curves
$\sigma^i(x):[0,1]\to\overline{\Omega}$, where $\Omega$ is a
smooth open convex subset of $\R^2$, that can intersect each other or
self--intersect only at their end points. We call ``multi--points''
the vertices $O^1, O^2,\dots, O^n\in\Omega$ of such smooth graph
${\mathbb S}$, where the order is greater than one. Moreover, 
we assume that  all the other ends of the curves (if present) 
have to coincide with some
points $P^l$ on the boundary of $\Omega$.\\
The problem is then to analyse the existence, uniqueness,
regularity and asymptotic behavior of the
evolution by curvature of such a network, under the constrain that the
end points $P^l\in\partial\Omega$ stay fixed.\\
Clearly, one can also set an analogous ``Neumann'' problem,
  requiring that, instead of being fixed, at the ``free'' ends the
  curves intersect orthogonally the boundary of $\Omega$.

\begin{rem}
  If $\Omega$ is a generic smooth open subset of
  $\R^2$, possibly non convex, it could happen that during the
  evolution one of the curves ``hits'' the boundary of $\Omega$ with a
  point different from its end point. 
  Then, for sake of simplicity, we assume 
  the convexity of the domain since
  this condition excludes a priori such an event
  if the lengths of the three curves are bounded away from zero 
  (Proposition~\ref{omegaok}).

  The hypothesis that the points $P^l$ stay on the boundary and not
  inside $\Omega$, is also in this spirit. If the fixed end points are
  inside the domain the interior of a curve could possibly ``touch''
  one of them, forming a loop and losing the ``embeddedness'' of the
  network.\\
  The convexity assumption and the fact that $P^l\in\partial\Omega$
  avoid this possibility, as we will see in Section~\ref{disuellesec}.
\end{rem}

In  such a generality, although simplified in this way, 
this problem shares various complications related
to the multi--points.\\
As previously underlined, the existing weak definitions of curvature
motion do not give uniqueness of the flow, or allow
``fattening'' phenomena (see~\cite{es},
for instance) which we would like to avoid, as they
seems quite extraneous to our setting. Among the existing notions,
the most suitable to our point of view is Brakke's
one (see Definition~\ref{brk} and the subsequent discussion), which
also lacks uniqueness but maintains the (Hausdorff) dimension of the sets,
excluding at least the event of fattening. In what follows 
this definition is the only that we consider
in relation to the evolution of networks, in particular, our flows are
{\em Brakke flows}.

In Section~\ref{kestimates} we show a 
satisfactory {\em small time} existence result
(Theorems~\ref{smoothexist},~\ref{brakkeevolution} and
Remark~\ref{manytrip}) of a smooth motion for a
special class of networks, that is,
the ones having only multi--points with three concurring curves
forming angles of $120$ degrees (this last property is called {\em
  Herring condition}).\\
We have to say that the uniqueness problem is less clear at the moment.\\
In the case of the presence in the initial network 
of a ``bad'' 3--point, not satisfying the Herring condition, 
we are not able at the moment to show the existence of a
flow, smooth for every positive time, satisfying a ``robust''
definition (at least as  Definition~\ref{brk}).\\
Actually, one would expect that
the desired good definition should give uniqueness of the motion and force,
by an instantaneous regularization, the three angles to become
immediately of $120$ degrees and to remain so.
This is sustained by the fact that,  by an energy argument
(\cite{brakke}), any smooth Brakke flow has to share
such a property (which is also suggested by numerical and physical
experiments, see at~{\em http://mimp.mems.cmu.edu}
and also the discussions
in~\cite{hermul,brakke,bronsard,gurtin2,kinderliu}).\\
Notice that, by the variational nature of the problem 
it is appealing to guess that some sort of 
parabolic regularization could play a role here.
We remark that if a multi--point has only
two curves concurring, it can be shown, by the regularizing effect
of the evolution by curvature (see~\cite{altgra,angen2,angen3,gray1}),
that the two curves together become instantaneously a single smooth
curve moving by curvature. Hence, the 2--point has vanished but this
particular event is so ``soft'' (and topologically null) 
that we can avoid to consider it as a real structural change.

We discuss now some other difficulties of geometric character 
in having a good definition of flow for a generic network.

\begin{enumerate}
\item {\em The presence of multi--points $O^j$ of order greater than
    three:}\\
  In the case of a 4--point (and clearly also of a higher order
  multi--point), for instance, considering the network described 
  by two curves crossing each other, 
  there are  really several possible candidates for the flow, even
  excluding a priori ``fattening'' phenomena. One cannot easily decide
  how the angles must behave, like in the 3--point case above,
  moreover, one can allow the four concurring curves to separate in
  two pairs of curves moving independently each other and it could even be
  taken into account the ``creation'' of new multi--points
  from such a single one (these events are actually possible in Brakke's
  definition).\\
  In these latter cases, the topology of the network changes
  dramatically, forcing us to change the structure of the equations
  governing the evolution or the family of curves composing the
  network.

\item {\em The presence of several multi--points $O^j$:}\\
  during the flow some of them can ``collapse'' together,
  again modifying the topological structure of the network,
  when the length of at least one curve
  of the network goes to zero (which can actually happen). In this
  case, like at the previous point, one possibly has to ``restart'' the
  evolution with a different set of curves.\\
  Notice that even if one starts with a network such that all the
  multi--points are 3--points, in the event of a collapse, one could
  have to face a situation with multi--points of order higher than
  three (consider for instance two 3--points collapsing along a single
  curve connecting them) or to deal with ``bad'' 3--points
  (think of three 3--points collapsing
  together along three curves connecting them).
\end{enumerate}

\begin{rem}\label{unstable}
  Actually, it seems reasonable that the configurations with
  multi--points of order greater than three or 3--points with angles
  different from $120$ degrees should be {\rm unstable} (actually, they
  are unstable for the {\em Length} functional), with the meaning that they
  can appear at some discrete set of times (and probably in some cases are
  unavoidable), but they must vanish immediately after.
\end{rem}

Because of all these complications, as a first step, in this paper 
we are going to
analyze the simplest possible network which, by construction,
rules out all the troubles related to the cases above.\\
Anyway, we will point out when the results can be extended to more
general networks.

\begin{dfnz}\label{trddfnz}
We call {\em triod} $\TTT$ in $\Omega$ a special network composed only
of three  regular, embedded $C^2$ curves 
$\sigma^i:[0,1]\to\overline{\Omega}$ (here {\em regular} means
$\sigma^i_x(x)\not=0$ for every $x\in[0,1]$ and $i\in\{1, 2, 3\}$),
where $\Omega$ is a  smooth open convex subset of $\R^2$, moreover,
these curves (sometimes we will call them also {\em edges of the triod}) 
intersect each other only at a single 3--point
$O=\sigma^1(0)=\sigma^2(0)=\sigma^3(0)$ (the {\em 3--point of the
  triod $\TTT$}) and have the other three end
points coinciding with three distinct points $P^i=\sigma^i(1)$ (the
{\em end points of $\TTT$}) belonging to the boundary of $\Omega$.\\
Finally, we assume that the tangents of the three curves form angles of
$120$ degrees at the 3--point $O$.\\
We say that the triod is (of class) $C^k$ or $C^\infty$ if the three
curves are respectively $C^k$ or $C^\infty$. We remark that with this
definition we assume that, unless explicitly otherwise stated, 
every triod is at least of class $C^2$, in
order to speak of its curvature.
\end{dfnz}
Notice that the angular condition can be expressed as
$$
\sum_{i=1}^3\frac{\sigma_x^i(0)}{\vert{\sigma_x^i(0)}\vert}=0.
$$

It is obvious that during the motion of such a triod
the ``bad'' configurations we discussed before
are a priori excluded.

\begin{rem}\label{noendp}
Since in all the paper we will consider only triods with angles of
$120$ degrees between the curves, for sake of simplicity, we chose to
call these sets simply triods with the meaning ``triods with angles of $120$
degrees''.\\
Moreover, sometimes we will speak of a {\em triod without end points}
composed of three curves $\sigma^i:[0,+\infty)\to\R^2$, for instance,
when we will need to allow them go to infinity.\\
Actually, we could also consider the possibility that the three
curves intersect each other, so when ambiguity is possible, we will
underline the property of non self--intersection saying that the triod is
{\em embedded}.
\end{rem}

The paper is devoted to study the existence, uniqueness, regularity
and asymptotic behavior of these embedded triods
in $\Omega$, moving by curvature keeping fixed the end points
$P^1, P^2, P^3$ .\\
This problem has been considered by Bronsard and
Reitich in~\cite{bronsard}, where they prove
an existence result which is the core of
Theorem~\ref{smoothexist} and by Kinderlehrer and Liu
in~\cite{kinderliu} showing the global existence of a 
flow for an initial triod sufficiently close to the minimal
configuration connecting the three points $P^i$ (Steiner
configuration).\\
We want to extend all this to {\em any} embedded initial triod,
concentrating in particular on the global existence and regularity of
the flow. Even if this is the simplest case, its
understanding is clearly crucial in analyzing more general networks,
taking also into account Remark~\ref{unstable}.

\smallskip

{\em Our conjecture is that any embedded
initial triod evolves in time without
singularities and asymptotically converges to the minimal connection
between the three points $P^i$ if the lengths
of the three curves stay away from zero, being the ``vanishing'' of a 
curve the only possible ``catastrophic'' event during the flow.}

\begin{rem}
It should be noticed that if the triangle formed by the points
$P^1, P^2, P^3$ has an angle of more than $120$ degrees,
then a triod composed of three segments forming angles of $120$ and
connecting its vertices does not exist.
Then, if the conjecture is true, necessarily one of the
lengths is not uniformly bounded away from zero.
\end{rem}

After discussing the existence (and uniqueness) of a smooth flow on
some maximal time interval in the first part of the paper, in order to
prove such conjecture we try to generalize
the analysis of the motion by mean curvature of smooth closed
curves and hypersurfaces in the Euclidean space, employing a
mix of PDE's and differential geometry techniques.\\
Essentially, what is needed is the understanding of the structure of
the possible {\em blow up} around the singularities, in order to actually
exclude these latter by means of geometric arguments. 
Some key references for this line of research 
are~\cite{altsch,hamilton4,huisk3,huisk2}.\\
The most relevant difference between our case and the smooth one, is the
difficulty in using the {\em maximum principle}, which is the
main tool to get estimates on the geometric quantities during the
flow. Indeed, the 3--point (by the $120$ degrees condition) is 
nice from the {\em distributional} point of view: in a sense, it
is an interior point, but it is troublesome for any argument based on
the maximum principle, since it behaves like a boundary point.\\
For this reason, some important pointwise estimates which are 
almost trivial applications of the maximum principle in the smooth case,
are here much more complicated to prove (sometimes we do
not even know if they actually hold) and we will have
to resort to integral estimates. These latter are similar to the ones
in~\cite{altsch,angen2,angen3,huisk1} for instance, but require some
extra work in order to deal with this strange ``boundary'' point.

If the length of the curves do not reduce to zero, 
by means of these latter estimates, we can see that at the time $T$ of
singularity the curvature has to explode, then like in the smooth
case, we separate the analysis according to its rate of blow up.\\
We say that a singularity is of {\em Type I} if for some
constant $C$ we have $\max_{\TTT_t}k^2\leq C/(T-t)$ as $t\to T$ and it
is of {\em Type II} otherwise.\\
Rescaling properly the flow around a hypothetical {\em Type I} singularity
one gets an evolution of embedded triods (unbounded and without end
points) shrinking homothetically during the motion by 
curvature. Classifying all such particular evolutions, 
we will show that none of them can arise as a blow up of the flow
$\TTT_t$, this clearly implies that {\em Type I} singularities cannot
develop.\\
With the same idea, rescaling the flow around a {\em Type II} singularity,
one gets an {\em eternal} motion by curvature, that is, an evolution
of triods defined for every time $t\in\R$.\\
What is missing at the moment is that this eternal flow is
actually simply given by a {\em translating} triod (unbounded and
without end points), like it happens in the case of a single smooth
curve. 
Here also, the main difficulty resides in replacing some maximum
principle arguments.\\ 
If the blow up would be translating, after classification, we could
exclude also this case by means of an argument based on the
monotonicity of a geometric quantity 
(see Sections~\ref{blw} and~\ref{typeIIsecn} for details), 
hence, no singularity at all could 
appear during the flow if the lengths of the
three curves of the triod stay away from zero.\\
The conjecture then would follow.

\smallskip

\begin{ackn} We are grateful to Alessandra Lunardi for helping us in
  the proof of the small time existence of a smooth flow in
  Section~\ref{kestimates}.
\end{ackn}

\section{Definitions and Preliminaries}
\label{due}

We start setting in a precise analytical way the curvature evolution problem
for an embedded initial triod $\TTT_0=\cup_{i=1}^3\sigma^i$ in
$\Omega$.

\begin{dfnz} We say that the triods
  $\TTT_t=\cup_{i=1}^3\gamma^i(\cdot, t)$ evolve
by curvature (remaining embedded) in the time interval $[0,T)$
if the three functions $\gamma^i:[0,1]\times[0,T)\to\overline{\Omega}$
are of class $C^2$ in space and $C^1$ in time, at least, and satisfy
the following quasilinear parabolic system
\begin{equation}\label{problema}
\begin{cases}
\begin{array}{lll}
\gamma_x^i(x,t)\not=0\qquad &&\text{ regularity}\\
\gamma^i(x,t)\not=\gamma^i(y,t)\qquad &\text{{ if} $x\not=y$}&\text{ simplicity}\\
\gamma^i(x,t)=\gamma^j(y,t)\,\Leftrightarrow\,x,y=0\qquad&\text{{    if} $i\not=j$}
&\text{ intersection only  at the 3--point}\\
\sum_{i=1}^3\frac{\gamma_x^i(0,t)}{\vert{\gamma_x^i(0,t)}\vert}=0
\qquad&&\text{ angles of $120$ degrees at the 3--point}\\
\gamma^i(1,t)=P^i
\qquad &&\text{ fixed end points condition}\\
\gamma^i(x,0)=\sigma^i(x)\qquad &&\text{ initial data}\\
\gamma^i_t(x,t)=\frac{\gamma_{xx}^i(x,t)}{{\vert{\gamma_x^i(x,t)}\vert}^2}\qquad
&&\text{ motion by curvature}
\end{array}
\end{cases}
\end{equation}
for every $x\in[0,1]$, $t\in[0,T)$ and $i, j\in\{1, 2, 3\}$.
\end{dfnz}

\begin{rem}
We spend some words on the evolution equation
\begin{equation}
\gamma^i_t=\frac{\gamma_{xx}^i}{{\vert{\gamma_x^i}\vert}^2}\,,\label{evoleq}
\end{equation}
which is not the usual way to describe the motion by curvature, that
is,
$$
\gamma^i_t=\frac{\langle\gamma_{xx}^i\,\vert\,\nu^i\rangle}
{{\vert{\gamma_x^i}\vert}^2}\,\nu^i=k^i\nu^i
$$
where we denoted with $\nu^i$ the unit normal to the curve
$\gamma^i$ and $k^i$ its curvature.\\
The two velocities differ only by a tangential component which
actually affects the motions of the single points (Lagrangian point of
view), but it does not affect the local motion of a curve as a whole
subset of $\R^2$ (Eulerian point of view).\\
We remark that this property of global invariance is not peculiar to
this particular tangential term, it holds  for any tangential
modification of the velocity. This fact is well
known for the curvature evolution of a smooth curve, 
hence also for a triod, any tangential
contribution to the velocity does not modify the flow outside the
3--point.\\
In our situation such extra term becomes
necessary in order to allow the motion of the 3--point
$O(t)=\gamma^i(0,t)$. 
Indeed, since we look for a $C^2$
solution of Problem~\eqref{problema}, if the velocity would be in normal
direction at every point of the three curves, the 3--point should
move in a direction which is normal to all of them, then the only
possibility would be that it does not move at all (see also the discussions
and examples in~\cite{brakke,bronsard,kinderliu}).
\end{rem}

\begin{rem}
It should be noticed that this definition of flow of a $C^2$ triod is
very strong, indeed, as the maps $\gamma^i$ have to be $C^1$ in time
and $C^2$ is space till the parabolic boundary, the compatibility
conditions of order $2$ have to be satisfied, that is,
$$
\text{$\frac{\sigma_{xx}^i(1)}{{\vert{\sigma_x^i(1)}\vert}^2}=0$
{ for every} $i\in\{1,2,3\}$ { and}
$\frac{\sigma_{xx}^i(0)}{{\vert{\sigma_x^i(0)}\vert}^2}=
\frac{\sigma_{xx}^j(0)}{{\vert{\sigma_x^j(0)}\vert}^2}$
for every $i,j\in\{1,2,3\}$}\,.
$$
(the compatibility conditions of order $0$ and $1$ are automatically
satisfied, since they are equivalent to say  that the three curves
$\sigma^i$ form a 3--point with angles of $120$ degrees, which is
assumed by construction).\\
This means, for instance, that for the initial triod the curvature 
at the end points $P^i$ and the sum of the three curvatures at the
3--point have to be zero (see later).\\
Notice that these two conditions on $\TTT_0$ are really geometric,
independent of the parametrization of the curves $\gamma^i$, but
intrinsic to the {\em set} $\TTT_0$ (see Definition~\ref{geocomp} and
Remark~\ref{geocomprem}) and are not satisfied by a generic $C^2$ triod.\\
Since we are interested in an  existence theorem for {\em any} 
triod composed of $C^2$ curves (Theorem~\ref{brakkeevolution}), we
consider also Brakke's definition of curvature flow.
\end{rem}

\begin{dfnz}\label{smoothflow}
  We say that a triod $\TTT_0$, composed of 
  three $C^\infty$ curves $\sigma^i$ is {\em smooth}
  if it satisfies the {\em compatibility conditions of every order} of
  parabolic system~\eqref{problema}
  (this is clearly much more striking that being simply a $C^\infty$
  triod according to Definition~\ref{trddfnz}).\\
  Precisely, this means that at the end points and at the
  3--point, there hold all the relations on the space derivatives of
  the functions $\sigma^i$ obtained differentiating in time, at $t=0$,
  the  boundary conditions.\\
  We say that a solution $\TTT_t$ of Problem~\eqref{problema},
  given by the curves  $\gamma^i(x,t)\in C^\infty([0,1]\times[0,T))$
  is a {\em smooth flow} if the compatibility conditions of every
  order are satisfied 
  at every time $t\in[0,T)$, that is, all the triods $\TTT_t$ are 
  {\em  smooth}.
\end{dfnz}

\begin{dfnz}\label{geocomp} 
We say that a triod $\TTT_0$ is {\em geometrically} smooth
  if there exist a regular parametrization of its three curves such
  that the resulting triod is smooth.
\end{dfnz}

To denote a flow we will often
write simply $\TTT_t$ instead of letting explicit the curves $\gamma^i$
which compose the triods.\\
Moreover, it will be also useful to describe a triod as a map
$F:\TTT\to\overline{\Omega}$ from a fixed {\em standard} triod $\TTT$
in $\R^2$,  composed of three unit segments from the origin in the
plane, forming angles of $120$ degrees. In this case we will still
denote with $O$ the 3--point of $\TTT$ and with $P^i$ the three
end points of such standard triod.\\
The evolution then will be given by a map
$F:\TTT\times[0,T)\to\overline{\Omega}$,
constructed naturally from  the curves
$\gamma^i$, so $\TTT_t=F(\TTT,t)$.

We adopt the following notation:
$$
\begin{array}{ll}
\tau^i=\tau^i(x,t)= \frac{\gamma_x^i}{\vert\gamma_x^i\vert}
&\qquad\text{{ unit tangent vector to} $\gamma^i$}\,,\\
\nu^i=\nu^i(x,t)= {\mathrm R} \tau^i(x,t)={\mathrm R} \frac{\gamma_x^i}{\vert\gamma_x^i\vert}
&\qquad\text{{ unit normal vector to} $\gamma^i$}\,,\\
O=O(t)=\gamma^i(0,t)
&\qquad\text{{ 3--point of the triod} $\TTT_t$}\,,\\
\underline{v}^i=\underline{v}^i(x,t)= \frac{\gamma_{xx}^i}{{\vert{\gamma_x^i}\vert}^2}
&\qquad\text{{ velocity of the point} $\gamma^i(x,t)$}\,,\\
\lambda^i=\lambda^i(x,t)= \frac{\langle\gamma_{xx}^i\,\vert\,\tau^i\rangle}
{{\vert{\gamma_x^i}\vert}^2}=\frac{\langle\gamma_{xx}^i\,\vert\,\gamma^i_x\rangle}
{{\vert{\gamma_x^i}\vert}^3}
&\qquad\text{{ tangential velocity of the point} $\gamma^i(x,t)$}\,,\\
k^i=k^i(x,t)= \frac{\langle\gamma_{xx}^i\,\vert\,\nu^i\rangle}
{{\vert{\gamma_x^i}\vert}^2}=
\langle\partial_s\tau^i\,\vert\,\nu^i\rangle=
-\langle\partial_s \nu^i\,\vert\,\tau^i\rangle
&\qquad\text{{ curvature at the point} $\gamma^i(x,t)$}\,,
\end{array}
$$
where $s$ is the arclength parameter on  the relative curve, defined by
$s(x,t)=\int_0^x\vert\gamma^i_x(\xi,t)\vert\,d\xi$, and
${\mathrm R}:\R^2\to\R^2$ is the counterclockwise rotation of
$\pi/2$. Notice that
$\partial_s=\vert\gamma_x^i\vert^{-1}\partial_x$.\\
Moreover, we set $\underline{\lambda}^i=\lambda^i\tau^i$ and
$\underline{k}^i=k^i\nu^i$, then, it clearly follows that
$\underline{v}^i=\underline{\lambda}^i+\underline{k}^i$ and $\vert
\underline{v}^i\vert^2=(\lambda^i)^2+(k^i)^2$.

Here and in the sequel, we denote with $\partial_sf$ and $\partial_tf$
the derivatives of a function $f$ along a curve with respect to the
relative arclength parameter and the time, $\partial^n_sf$, $\partial^n_tf$ are
the higher order partial derivatives which often we will also write as
$f_{s}, f_{ss},\dots$ and  $f_t, f_{tt},\dots$.\\
We adopt the following convention for integrals,
$$
\int_{{\TTT_t}} f(t,\gamma,\tau,\nu,k,k_s,\dots,\lambda,\lambda_s\dots)\,ds =
\sum_{i=1}^3 \int_0^1
f(t,\gamma^i,\tau^i,\nu^i,k^i,k^i_s,\dots,\lambda^i,\lambda^i_s\dots)\,\vert
\gamma^i_x\vert\,dx
$$
as the arclength measure
is given by $ds=\vert\gamma^i_x\vert\,dx$ on the curve $\gamma^i$.\\
In general, if there is no need to make explicit the three curves
composing a triod, we simply write $\tau$, $\nu$, $\underline{v}$,
$\lambda$ and $k$ for the previous quantities, omitting the indices.

We suppose now to have a smooth flow $\TTT_t$ on some positive time interval
$[0,T)$ and we write the evolution equations for the geometric
quantities.

\begin{lemma} If $\gamma$ is a curve moving by
$$
\gamma_t=\frac{\gamma_{xx}}{{\vert\gamma_x\vert}^2}=\lambda\tau+k\nu\,.
$$
then the following commutation rule holds,
\begin{equation}\label{commut}
\dert\ders=\ders\dert +
(k^2 -\lambda_s)\ders
\end{equation}
\end{lemma}
\begin{proof}
Let $f:[0,1]\times[0,T)\to\R$ be a smooth function, then
\begin{align*}
\dert\ders f - \ders\dert f =&\, \frac{f_{tx}}{\vert\gamma_x\vert} -
\frac{\langle \gamma_x\,\vert\,\gamma_{xt}\rangle f_x}
{\vert\gamma_x\vert^3}
- \frac{f_{tx}}{\vert\gamma_x\vert} = - {\langle
  \tau\,\vert\,\partial_s\gamma_t\rangle}\partial_sf\\
=&\, - {\langle\tau\,\vert\,\partial_s(\lambda\tau+k\nu)\rangle}\partial_sf=
 (k^2 - \lambda_s)\ders f
\end{align*}
and the formula is proved.
\end{proof}
Then we can compute,
\begin{align}
\dert\tau=&\,
\dert\ders\gamma=\ders\dert\gamma+(k^2-\lambda_s)\ders\gamma =
\ders(\lambda\tau+k\nu)+(k^2-\lambda_s)\tau =
(k_s+k\lambda)\nu\label{derttau}\\
\dert\nu=&\, \dert({\mathrm R}\tau)={\mathrm
R}\,\dert\tau=-(k_s+k\lambda)\tau\label{dertdinu}\\
\dert k=&\, \dert\langle \ders\tau\,|\, \nu\rangle=
\langle\dert\ders\tau\,|\, \nu\rangle\label{dertdik}
= \langle\ders\dert\tau\,|\, \nu\rangle +
(k^2-\lambda_s)\langle\ders\tau\,|\, \nu\rangle\\
=&\, \ders\langle\dert\tau\,|\, \nu\rangle + k^3-k\lambda_s =
\ders(k_s+k\lambda) + k^3-k\lambda_s\nonumber\\
=&\, k_{ss}+k_s\lambda + k^3\nonumber\\
\dert\lambda =&\, -\dert\partial_x\frac{1}{\vert\gamma_x\vert}=
\partial_x \frac{\langle\gamma_x\,\vert\,\gamma_{tx}\rangle}
{\vert\gamma_x\vert^3}=
\partial_x \frac{\langle\tau\,\vert\,\ders (\lambda\tau+k\nu)\rangle}
{\vert\gamma_x\vert}=\partial_x \frac{(\lambda_s - k^2)}
{\vert\gamma_x\vert}\label{dertdilamb}\\
=&\, \ders(\lambda_s - k^2) -\lambda(\lambda_s -
k^2)=\lambda_{ss} -\lambda\lambda_s - 2kk_s +\lambda k^2\,.\nonumber
\end{align}
As the triods $\TTT_t$ are smooth, 
differentiating in time the concurrency condition
$\gamma^i(0,t)=\gamma^{j}(0,t)$ we obtain
$\lambda^i\tau^i+k^i\nu^i=\lambda^j\tau^j+k^j\nu^j$ at
the 3--point for every pair of indices $i, j$. Multiplying these vector equalities
for $\tau^l$ and $\nu^l$ and varying $i, j, l$ we get the relations
\begin{gather*}
\lambda^i=-\lambda^{i+1}/2-\sqrt{3}k^{i+1}/2\\
\lambda^i=-\lambda^{i-1}/2+\sqrt{3}k^{i-1}/2\\
k^i=-k^{i+1}/2+\sqrt{3}\lambda^{i+1}/2\\
k^i=-k^{i-1}/2-\sqrt{3}\lambda^{i-1}/2
\end{gather*}
with the convention that the superscripts are
considered modulus three. Solving this system we get
\begin{gather*}
\lambda^i=\frac{k^{i-1}-k^{i+1}}{\sqrt{3}}\\
k^i=\frac{\lambda^{i+1}-\lambda^{i-1}}{\sqrt{3}}
\end{gather*}
which implies
\begin{equation*}
\sum_{i=1}^3 k^i=\sum_{i=1}^3\lambda^i=0
\end{equation*}
at the 3--point of the triods $\TTT_t$.\\
Moreover, considering ${\mathrm K}=(k^1,k^2,k^3)$ and
${\mathrm\Lambda}=(\lambda^1,\lambda^2,\lambda^3)$ as vectors in
$\R^3$, we have seen that ${\mathrm K}$ and ${\mathrm\Lambda}$ belong
to the plane orthogonal to the vector $(1,1,1)$ and 
${\mathrm K}={\mathrm \Lambda} \wedge (1,1,1)/\sqrt{3}$, 
${\mathrm  \Lambda}=-{\mathrm K} \wedge (1,1,1)/\sqrt{3}$ that
is, ${\mathrm K}={\mathrm S}{\mathrm \Lambda}$ and ${\mathrm
  \Lambda}=-{\mathrm S}{\mathrm K}$ where ${\mathrm
  S}$ is the rotation in $\R^3$ of an angle of $\pi/2$ around the axis
${\mathrm I}=\langle(1,1,1)\rangle$. Hence, it also follows that 
\begin{equation*}
\sum_{i=1}^3 (k^i)^2=\sum_{i=1}^3 (\lambda^i)^2\qquad\text{ { and}
  }\qquad \sum_{i=1}^3 k^i\lambda^i=0\,.
\end{equation*}
at the 3--point of the triods $\TTT_t$.\\
Now we differentiate in time the angular condition
$\sum_{i=1}^3\tau^i(0,t)=0$ at the 3--point, by equation~\eqref{derttau} we get
\begin{equation}\label{topolino4}
k_s^i+\lambda^i k^i=k_s^j+\lambda^j k^j
\end{equation}
for every pair $i, j$. In terms of vectors $\R^3$ as before, we can write
$$
{\mathrm K}_s+{\mathrm {\Lambda K}}=(k^1_s+\lambda^1
  k^1,k^2_s+\lambda^2 k^2,k^3_s+\lambda^3 k^3)\in {\mathrm I}\,.
$$
Differentiating repeatedly in time all these vector relations we have
\begin{gather*}
\partial_t^l{\mathrm K}\,,\,\partial_t^l{\mathrm \Lambda}\perp{\mathrm
  I}\, \text{ and }\, \partial_t^l({\mathrm {K\dot\Lambda}})=0\,,\\
\dert^l{\mathrm \Lambda}=-\dert^l{\mathrm S}{\mathrm K}=-{\mathrm
  S}\dert^l{\mathrm K}\,,\\
\partial_t^m({\mathrm K}_s+{\mathrm {\Lambda K}})\in {\mathrm I}\,,
\end{gather*}
which, making explicit the indices, give the following equalities at
the 3-point,
\begin{gather}
\dert^l\sum_{i=1}^3 k^i=\sum_{i=1}^3\dert^l k^i=\dert^l\sum_{i=1}^3 
\lambda^i=\sum_{i=1}^3\dert^l\lambda^i=\dert\sum_{i=1}^3k^i\lambda^i=0\,,\nonumber\\
\sum_{i=1}^3 (\dert^l k^i)^2=\sum_{i=1}^3
(\dert^l\lambda^i)^2\label{lambdakappa}\,\,\text{ { for every} $l\in\NN$,}\\
\dert^m(k^i_s+\lambda^i k^i)=\dert^m(k^j_s+\lambda^j k^j)
\,\,\text {  for every pair $i, j$ and $m\in\NN$.}\nonumber
\end{gather}
By the orthogonality relations with respect to the axis ${\mathrm I}$, we
get also
$$
\partial^l_tK\partial_t^m({\mathrm K}_s+{\mathrm {\Lambda
    K}})=\partial^l_t\Lambda\partial_t^m({\mathrm K}_s+{\mathrm
  {\Lambda K}})=0\,,$$
that is,
\begin{equation}\label{orto}
\sum_{i=1}^3 \dert^lk^i\,\dert^m(k^i_s+\lambda^i k^i)
=\sum_{i=1}^3 \dert^l\lambda^i\,\dert^m(k^i_s+\lambda^i k^i)=0
\,\,\text { for  every $l, m\in\NN$.}
\end{equation}

Looking then at the three end points, 
by Lemma~\ref{evenly} in the next section, we
have that all the {\em even} space derivatives of $k$ and $\lambda$
are zero.

\begin{rem}\label{geocomprem}
As the ``shape'' of the curves is fixed under reparametrization, if a triod is
geometrically smooth, the three curves have to be $C^\infty$ once
parametrized in arclength. Moreover, necessarily $\sum_{i=1}^3\nu^i=0$, 
and $\sum_{i=1}^3k^i=0$ at the 3--point must hold, like 
all the relations one gets from these two, iteratively, differentiating
in time by means of formulas~\eqref{dertdinu},~\eqref{dertdik} 
and substituting every occurrence of $\lambda^i$ according to the formula
$(\lambda^1,\lambda^2,\lambda^3)=-(k^1,k^2,k^3)\wedge(1,1,1)/\sqrt{3}$
(notice that in this way $\lambda_s$ never appears). 
Working analogously at the end points, the compatibility conditions
reduce to require that every even space derivative of the curvature is
zero (by Lemma~\ref{evenly} in the next section).\\
These necessary conditions are actually also sufficient, indeed 
the {\em geometrical smoothness} is a matter of curvature, not
involving $\lambda$. If a triod satisfies such conditions, we can
parametrize every curve $\gamma^i$ in a way that $\lambda^i$ at
the 3--point has the right value given by the relation
$\lambda^i=-({\mathrm K}\wedge{\mathrm
  I})^i/\sqrt{3}$ and it is zero at the end point of the
curve, for instance, setting $\lambda^i_s$ constantly equal to
$-\lambda^i(0)/L^i$, where $L^i$ is the length of the curve.\\
It is straightforward to check that this can be
done and the resulting parametrization gives a smooth triod.
\end{rem}

\begin{dfnz}\label{brk}
We will speak of {\em Brakke flow with equality} of an initial triod
$\TTT_0$ in $[0,T)$,
for a family of $C^2$ triods $\TTT_t$ in $\Omega$ all with the same end
points as $\TTT_0$ and satisfying the equation
\begin{equation}\label{brakkeqqq}
\frac{d\,}{dt}\int_{\TTT_t}\varphi(\gamma,t)\,ds=
-\int_{\TTT_t}\varphi(\gamma,t) k^2\,ds
+\int_{\TTT_t}\langle\nabla\varphi(\gamma,t)\,\vert\,\underline{k}\rangle\,ds
+ \int_{\TTT_t}\varphi_t(\gamma,t)\,ds\,,
\end{equation}
for every smooth function with compact support
$\varphi:\Omega\times[0,T)\to\R$ and $t\in[0,T)$.\\
This means also that the time derivative at the left
member has to exist. The right member does not give any problem
since the triods are $C^2$, at least.

We will say that a Brakke flow is {\em smooth} if all the triods are
geometrically smooth.
\end{dfnz}

\begin{rem}
It is straightforward to check that a solution of Problem~\ref{problema}
is also a smooth Brakke flow with equality.

Actually, the original definition of Brakke flow stated
in~\cite[Section~3.3]{brakke}
allows equality~\eqref{brakkeqqq} to be an inequality (and triods 
$\TTT_t$ to be one--dimensional countably 
rectifiable subsets of $\R^2$ with a distributional notion of
curvature, called {\em varifolds}, see~\cite{simon}),
precisely,
\begin{equation*}
\frac{\overline{d\,}}{dt}\int_{\TTT_t}\varphi(x,t)\,d\HH^1(x)\leq
-\int_{\TTT_t}\varphi(x,t) k^2\,d\HH^1(x)
+\int_{\TTT_t}\langle\nabla\varphi(x,t)
\,\vert\,\underline{k}\rangle\,d\HH^1(x)
+ \int_{\TTT_t}\varphi_t(x,t)\,d\HH^1(x)\,,
\end{equation*}
must hold for every {\em positive} smooth function with compact support
$\varphi:\Omega\times[0,T)\to\R$ and $t\in[0,T)$, where 
$\frac{\overline{d\,}}{dt}$ is the {\em upper} derivative (the
$\limup$ of the incremental ratios) and $\HH^1$ 
is the Hausdorff one--dimensional measure in $\R^2$ (we will use this 
notation through all the paper).\\
This weaker condition was introduced by Brakke
in order to prove an existence result~\cite[Section~4.13]{brakke} for
a family of initial sets much wider than the networks of curves, but, on the
other hand, it let open the possibility of instantaneous vanishing of
some parts of the set.\\
Since for our triods (probably, this can be done also for a general
network) we are able to show the existence of a Brakke flow 
via a different method, and this flow is composed of smooth triods and
satisfies the equality, for sake of simplicity, we included such extra
properties in the definition.
\end{rem}

A big difference between Brakke flows and the evolutions obtained as
solutions of Problem~\eqref{problema} is that the former triods are
simply considered as {\em sets} without any mention to their
parametrization (that clearly is not unique). This means that actually
a Brakke flow can be a family of triods given by the maps $\gamma^i(x,t)$
which are $C^2$ in space, but possibly do not have absolutely any
regularity with respect to the time variable $t$.

If we consider two different smooth triods $\TTT_0^1$ and $\TTT_0^2$
which are the same subset of $\R^2$, 
giving two different solutions $\TTT_t^1$ and
$\TTT^2_t$ of Problem~\eqref{problema} on some common interval,
then the two associated Brakke flows coincide
(Proposition~\ref{geouniq}), that is, as subsets of $\R^2$, forgetting
the parametrization, actually $\TTT_t^1=\TTT^2_t$ for every time $t$.\\
This means that the geometric evolution problem has a satisfactory
uniqueness property if the initial triod is smooth.

{\em In general, when we will speak of {\em geometric problem} we will
  mean that we are thinking of the triods as subsets of $\R^2$,
  independently of the parametrizations of their curves.} 

An open question is whether any smooth Brakke flow with equality
admits a parametrization of the initial triod with an associate
solution of Problem~\eqref{problema} representing it at least for
some time.\\
A positive answer would imply the uniqueness in the class of these
special Brakke flows and the coincidence of the two formulations, from
the geometric point of view.

\smallskip

Finally we state precisely the conjecture which is the main topic of
the second part of the paper.\\
Here and in the following, we 
denote with $L^i$ the lengths of the three curves 
and with $L=L^1+L^2+L^3$ the total length of the triod. 

\begin{conge}\label{congettura} Let $\TTT_t$ be a smooth evolution of 
embedded triods on a maximal time interval $[0,T)$.\\
If $\limnf_{t\to T} L^i\not=0$ for every $i\in\{1, 2, 3\}$, then
$T=+\infty$ and $\TTT_t$ converges, as $t\to+\infty$, to the minimal
connection between the three points $P^i$.
\end{conge}

\section{Small Time Existence and A Priori Estimates}
\label{kestimates}

The first small time existence result for a flow
very similar to Problem~\eqref{problema} is due to Bronsard and
Reitich~\cite{bronsard}. In their paper it is shown the existence of a
{\em unique} solution  $\gamma^i\in C^{2+2\alpha,1+\alpha}([0,1]\times[0,T])$
of the same parabolic system, for an initial
triod composed of three curves $\sigma^i\in C^{2+2\alpha}([0,1])$ and
satisfying the natural compatibility conditions. The only
difference is that they impose the 
{\em Neumann} boundary condition of orthogonal intersection with
$\partial\Omega$, instead of keeping the end points
$P^i\in\partial\Omega$ fixed as we do.

The same technique works also in our case and
gives the small time existence of a {\em unique} solution
$\gamma^i\in C^{2+2\alpha,1+\alpha}([0,1]\times[0,T])$ (notice that
this means that the curves are $C^{2+2\alpha}$ till the 3--point and
their end points $P^i$) of the following parabolic system
\begin{equation}\label{prob0}
\begin{cases}
\begin{array}{ll}
\gamma_x^i(x,t)\not=0\qquad &\qquad\text{ regularity}\\
\gamma^i(0,t)=\gamma^j(0,t)\qquad&\qquad\text{ concurrence at the
  3--point}\\
\sum_{i=1}^3\frac{\gamma_x^i(0,t)}{\vert{\gamma_x^i(0,t)}\vert}=0\qquad&\qquad\text{
  angles of $120$ degrees at the 3--point}\\
\gamma^i(1,t)=P^i\qquad &\qquad\text{ fixed end points condition}\\
\gamma^i(x,0)=\sigma^i(x)\qquad &\qquad\text{ initial data}\\
\gamma^i_t(x,t)=\frac{\gamma_{xx}^i(x,t)}{{\vert{\gamma_x^i(x,t)}\vert}^2}\qquad
&\qquad\text{ motion by curvature}
\end{array}
\end{cases}
\end{equation}
given any initial $C^{2+2\alpha}$ triod $\TTT_0=\cup_{i=1}^3\sigma^i$,
with $\alpha\in(0,1/2)$,  satisfying the {\em compatibility conditions of
  order $2$}, that is,
$$
\text{$\frac{\sigma_{xx}^i(1)}{{\vert{\sigma_x^i(1)}\vert}^2}=0$
{ for every} $i\in\{1,2,3\}$ { and}
$\frac{\sigma_{xx}^i(0)}{{\vert{\sigma_x^i(0)}\vert}^2}=
\frac{\sigma_{xx}^j(0)}{{\vert{\sigma_x^j(0)}\vert}^2}$
for every $i,j\in\{1,2,3\}$}
$$
(the compatibility conditions of order $0$ and $1$ are automatically
satisfied, since they are equivalent to say  that the three curves
$\sigma^i$ form a triod with angles of $120$ degrees).

Now we look for a higher regularity result.

\begin{teo}\label{smoothexist}
  For any initial smooth triod $\TTT_0$ there exists a unique
  smooth solution of Problem~\eqref{problema} on a maximal time
  interval $[0,T)$.
\end{teo}
\begin{proof} Since the initial triod $\TTT_0$ satisfies the
  compatibility conditions at every order,
  the method of Bronsard and Reitich actually
  provides a way, for every $n\in\NN$,  to get a unique solution in
  $C^{2n+2\alpha,n+\alpha}([0,1]\times[0,T_n])$, satisfying the compatibility
  conditions of order $0, \dots, n$ at every time.\\
  Then, by standard methods of one--dimensional parabolic equations we
  can obtain a solution which belongs to
  $C^{\infty}([0,1]\times[0,T))$ for some small positive time $T>0$
  and consider its maximal time of existence.\\
  We give just the line of the proof and we indicate the relevant
  references for the details.

Let us consider a solution $\gamma^i\in
C^{2n+2\alpha,n+\alpha}([0,1]\times[0,T_n])$ for $n\geq 2$, then the 
functions $\gamma_x^i(x,t)$  belong to $C^{2n-1+2\alpha,n-1/2
  +\alpha}([0,1]\times[0,T_n])$ (see for
instance~\cite[Lemma~5.1.1]{lunardi1}), then we look at
the parabolic system satisfied by $v^i(x,t)=\gamma_t^i(x,t)$,
\begin{equation*}
\begin{cases}
v^i_t(x,t)=\frac{v_{xx}^i(x,t)}{{\vert{\gamma_x^i(x,t)}\vert}^2} -
2\frac{\langle v^i_x(x,t)\,\vert\, \gamma_x^i(x,t)\rangle
\gamma_{xx}^i(x,t)}{{\vert{\gamma_x^i(x,t)}\vert}^4}\\
v^i(0,t)=v^j(0,t)\\
\sum_{i=1}^3\frac{v_x^i(0,t)}{\vert{\gamma_x^i(0,t)}\vert}-
\frac{\langle v^i_x(0,t)\,\vert\, \gamma_x^i(0,t)\rangle
  \gamma_x^i(0,t)}{\vert{\gamma_x^i(0,t)}\vert^3}=0\\
v^i(1,t)=0\\
v^i(x,0)=\frac{\sigma^i_{xx}(x)}{\vert\sigma_x^i(x)\vert^2}
\end{cases}
\end{equation*}
for every $i, j\in\{1, 2, 3\}$.\\
This system can be rewritten as
\begin{equation*}
\begin{cases}
v^i_t(x,t)=v_{xx}^i(x,t) f^i(x,t) + \langle v^i_x(x,t)\,\vert\,
g^i(x,t)\rangle\\
v^i(0,t)=v^j(0,t) \qquad &\qquad\\
\sum_{i=1}^3 v_x^i(0,t)p^i(t) + \langle v^i_x(0,t)\,\vert\,
q^i(t)\rangle r^i(t)=0\\
v^i(1,t)=0\\
v^i(x,0)=h^i(x)
\end{cases}
\end{equation*}
with coefficients $f^i, g^i\in C^{2n-2+2\alpha,n-1
  +\alpha}([0,1]\times[0,T_n])$, $p^i, q^i, r^i\in C^{2n-1+\alpha,n-1/2
  +\alpha}([0,T_n])$ and $h^i\in C^{2n+2\alpha}([0,1])$, since the
  initial triod is smooth.\\
By Solonnikov~\cite{solonnikov1} results, $v^i=\gamma^i_t$ belongs to
$C^{2n+2\alpha,n+\alpha}([0,1]\times[0,T_n])$
and since $\gamma^i_{xx}=\gamma^i_t\vert \gamma^i_x\vert^2$ with
$ \vert \gamma^i_x\vert^2\in C^{2n-1+2\alpha,n-1/2+\alpha}([0,1]\times[0,T_n])$, we
get also $\gamma^i_{xx}\in C^{2n-1+2\alpha,n-1/2+\alpha}([0,1]\times[0,T_n])$.\\
Following~\cite{lusiw}, we can then conclude that
$\gamma^i\in C^{2n+1+2\alpha,n+1/2+\alpha}([0,1]\times[0,T_n])$.\\
Iterating this argument, we see that $\gamma^i\in
C^{\infty}([0,1]\times[0,T_n])$, moreover, since for every $n\in\NN$
the solution obtained via the method of Bronsard and Reitich is
unique, it must coincide with $\gamma^i$ and we can choose all the
$T_n$ to be the same positive value $T$. Finally, by the same reason, 
all the compatibility conditions are satisfied at every time, that is,
the evolving triods are smooth.\\
The facts that these triods actually stay in the convex set $\Omega$ and
that they do not develop self--intersections during the flow will
follow by the results of Section~\ref{disuellesec}.
\end{proof}

\begin{prop}\label{equality1000} Any solution of Problem~\eqref{problema} is a smooth
  Brakke flow with equality.\\
Moreover,  for every curve $\gamma^i(\cdot,t)$ and for every time
$t\in[0,T)$ we have
$$
\frac{dL^i(t)}{dt}=-\lambda^i(0,t)-\int_{\gamma^i(\cdot,t)}k^2\,ds
$$
and
$$
\frac{dL(t)}{dt}=-\int_{\TTT_t}k^2\,ds\,.
$$
Hence, the total length $L(t)$ is decreasing in time and uniformly
bounded by the length of the initial triod $\TTT_0$.
\end{prop}
\begin{proof}
The geometrical smoothness of the flow is clear.\\
The time derivative of the measure $ds$ on the curve $\gamma^i$ is
given by $(\lambda_s^i-(k^i)^2)\,ds$, considering a smooth function
with compact support $\varphi:\Omega\times[0,T)\to\R$, we compute
\begin{align*}
\frac{d\,}{dt}\int_{\gamma^i(\cdot,t)}\varphi(\gamma^i,t)\,ds
= &\, \int_{\gamma^i(\cdot,t)}\varphi(\gamma^i,t) (\lambda^i_s-(k^i)^2)\,ds
+\int_{\gamma^i(\cdot,t)}\langle\nabla\varphi(\gamma^i,t)\,\vert\,\underline{v}^i\rangle\,ds
+ \int_{\gamma^i(\cdot,t)}\varphi_t(\gamma^i,t)\,ds\\
=&\,\int_{\gamma^i(\cdot,t)}\partial_s(\lambda^i\varphi(\gamma^i,t))
-\lambda^i\langle\nabla\varphi(\gamma^i,t)\,\vert\,\tau^i\rangle-\varphi(\gamma^i,t) (k^i)^2\,ds\\
&\, +\int_{\gamma^i(\cdot,t)}\langle\nabla\varphi(\gamma^i,t)\,\vert\,\underline{v}^i\rangle\,ds
+ \int_{\gamma^i(\cdot,t)}\varphi_t(\gamma^i,t)\,ds\\
=&\,\int_{\gamma^i(\cdot,t)}\partial_s(\lambda^i\varphi(\gamma^i,t))
-\varphi(\gamma^i,t) (k^i)^2\,ds\\
&\, +\int_{\gamma^i(\cdot,t)}\langle\nabla\varphi(\gamma^i,t)\,\vert\,\underline{k}^i\rangle\,ds
+ \int_{\gamma^i(\cdot,t)}\varphi_t(\gamma^i,t)\,ds\\
=&\,-\int_{\gamma^i(\cdot,t)}\varphi(\gamma^i,t) (k^i)^2\,ds
+\int_{\gamma^i(\cdot,t)}\langle\nabla\varphi(\gamma^i,t)\,\vert\,\underline{k}^i\rangle\,ds\\
&\, + \int_{\gamma^i(\cdot,t)}\varphi_t(\gamma^i,t)\,ds
+\lambda^i(1,t)\varphi(P^i,t)-\lambda^i(0,t)\varphi(O(t),t)\,.
\end{align*}
Since $\lambda^i(1,t)$ is zero, being zero the velocity $v$ at the end
points $P^i$, and since the sum of $\lambda^i$ is zero at the
3--point, adding these three equalities for $i\in\{1, 2, 3\}$ we obtain
formula~\eqref{brakkeqqq}.\\
The formulas for the lengths are given by the same computation with
$\varphi\equiv 1$.
\end{proof}

Seeing the initial triod $\TTT_0$ simply as a subset of $\R^2$, it can admit
more than a single parametrization of its curves making it a smooth
triod, so there are various flows arising by this theorem, associated to different
parametrizations. The following proposition shows that they
must coincide geometrically.

\begin{prop}\label{geouniq} If $\TTT^1_0=\cup_{i=1}^3\sigma^i_0$ and
  $\TTT^2_0=\cup_{i=1}^3\xi^i_0$ are two smooth triods which
  coincide as subset of $\R^2$ and $\TTT_t^1$, $\TTT_t^2$ are the
  relative flows given by Theorem~\ref{smoothexist} on a common time
  interval $[0,T)$, then at every time $t$ the triods $\TTT_t^1$ and
  $\TTT_t^2$ coincide as sets.
\end{prop}
\begin{proof}
Let $\gamma^i(x,t)$ and $\eta^i(x,t)$ be the two smooth flows associated to
$\sigma^i=\gamma^i(\cdot,0)$ and 
$\xi^i=\eta^i(\cdot,0)$, which parametrize  the same triod, seen as
a subset of $\R^2$.\\
We fix an index $i\in\{1, 2, 3\}$. Since $\sigma^i$ and $\xi^i$
are smooth regular parametrization of the same curve of the initial
triod, the map
$\varphi^i=(\sigma^i)^{-1}\composed\xi^i:[0,1]\to[0,1]$ is an orientation
preserving, smooth diffeomorphisms of the unit interval with 
itself, hence $\varphi^i(0)=0$ and $\varphi^i(1)=1$. 
Moreover, by the compatibility conditions we have
$\sigma^i_{xx}(1)=\xi^i_{xx}(1)=(\sigma^i\composed\varphi^i)_{xx}(1)=0$, hence
$$
0=\sigma^i_{xx}(\varphi^i(1))\vert\varphi^i_x(1)\vert^2+\sigma^i_{x}(\varphi^i(1))\varphi^i_{xx}(1)=
\sigma^i_{x}(\varphi^i(1))\varphi^i_{xx}(1)=\sigma^i_{x}(1)\varphi^i_{xx}(1)
$$
which implies that $\varphi^i_{xx}(1)=0$. At the 3--point, we have
\begin{align*}
\frac{\xi_{xx}^i(0)}{{\vert{\xi_x^i(0)}\vert}^2}-\frac{\sigma_{xx}^i(0)}{{\vert{\sigma_x^i(0)}\vert}^2}
=&\, \frac{{\sigma^i_{xx}(\varphi^i(0))\vert\varphi^i_x(0)\vert^2
+\sigma^i_{x}(\varphi^i(0))\varphi^i_{xx}(0)}}
{\vert\sigma^i_{x}(\varphi^i(0))\varphi^i_{x}(0)\vert^2}
-\frac{\sigma_{xx}^i(0)}{{\vert{\sigma_x^i(0)}\vert}^2}\\
=&\,\frac{{\sigma^i_{xx}(0)\vert\varphi^i_x(0)\vert^2+\sigma^i_{x}(0)\varphi^i_{xx}(0)}
-\sigma_{xx}^i(0)\vert\varphi^i_x(0)\vert^2}{\vert\sigma^i_{x}(0)\varphi^i_{x}(0)\vert^2}\\
=&\,\frac{\sigma^i_{x}(0)\varphi^i_{xx}(0)}
{\vert\sigma^i_{x}(0)\varphi^i_{x}(0)\vert^2}\,,
\end{align*}
hence, $\frac{\sigma^i_{x}(0)\varphi^i_{xx}(0)}
{\vert\sigma^i_{x}(0)\varphi^i_{x}(0)\vert^2}=\frac{\sigma^j_{x}(0)\varphi^j_{xx}(0)}
{\vert\sigma^j_{x}(0)\varphi^j_{x}(0)\vert^2}$ for every pair
$i, j\in\{1, 2, 3\}$. This means that $\varphi^i_{xx}(0)=0$ since the
tangents to the curves $\sigma^i$ and $\sigma^j$ are not parallel.\\
Now we look for three smooth functions
$\psi^i:[0,1]\times[0,T)\to[0,1]$ 
satisfying the following parabolic system
\begin{equation}
\begin{cases}
\begin{array}{l}
\psi^i_t(x,t)=\frac{\psi_{xx}^i(x,t)}{{\vert
    \gamma^i_x(\psi^i(x,t),t)\vert}^2{\vert\psi_x^i(x,t)\vert}^2}\\
\psi^i(0,t)=0\\
\psi^i(1,t)=1\\
\psi^i(x,0)=\varphi^i(x)
\end{array}
\end{cases}
\end{equation}
for every $(x,t)\in[0,1]\times[0,T)$ and every index $i\in\{1, 2, 3\}$.\\
We see that the compatibility conditions of order $2$ are satisfied by the
initial data, indeed, here these reduce only to
$\varphi^i_{xx}(0)=\varphi^i_{xx}(1)=0$.\\
By standard methods (now the problem is {\em scalar}, 
see~\cite{lunardi1,solonnikov1}), being $\varphi^i$
smooth and regular  ($\varphi^i_x(x)\not=0$ for every $x\in[0,1]$ since it
is a diffeomorphisms), the functions $\gamma^i_x$ bounded from above
and away from zero, and holding the compatibility conditions of order $0$, $1$ and
$2$, this quasilinear problem has a solution on some maximal time interval
$[0,T^\prime)$, with $T^\prime\leq T$, belonging to 
$C^{2+2\alpha,1+\alpha}([0,1]\times[0,T^{\prime\prime}])$ 
(for some $\alpha\in(0,1/2)$) for every time $T^{\prime\prime}\in(0,T^\prime)$.\\
It is now straightforward to see that the functions 
$\theta^i(x,t)=\gamma^i(\psi^i(x,t),t)$ coincide with $\eta^i(x,t)$ at
time $t=0$ (indeed, $\xi^i=\sigma^i\composed\varphi^i$) and 
\begin{align*}
\frac{\partial \theta^i(x,t)}{\partial t}
=\frac{\partial \gamma^i(\psi^i(x,t),t)}{\partial t}
=&\,\gamma^i_t(\psi^i(x,t),t)+\gamma_x^i(\psi^i(x,t),t)\psi^i_t(x,t)\\
=&\,\frac{\gamma^i_{xx}(\psi^i(x,t),t)}{\vert
  \gamma^i_x(\psi^i(x,t),t)\vert^2}
+\gamma_x^i(\psi^i(x,t),t)\frac{\psi_{xx}^i(x,t)}{{\vert
    \gamma^i_x(\psi^i(x,t),t)\vert}^2{\vert\psi_x^i(x,t)\vert}^2}\\
=&\,\frac{\gamma^i_{xx}(\psi^i(x,t),t)\vert\psi^i_x(x,t)\vert^2
+\gamma_x^i(\psi^i(x,t),t)\psi_{xx}^i(x,t)}{{\vert
    \gamma^i_x(\psi(x,t),t)\vert}^2{\vert\psi_x^i(x,t)\vert}^2}\\
=&\,\frac{(\gamma^i(\psi^i(x,t),t))_{xx}}{\vert(\gamma^i(\psi^i(x,t),t))_{x}\vert^2}
=\frac{\theta^i_{xx}(x,t)}{\vert\theta^i_x(x,t)\vert^2}\,.
\end{align*}
Then, since such functions $\theta^i$ and $\eta^i$ satisfy both
system~\eqref{prob0} and they take the same initial data at time
$t=0$, they must coincide on $[0,T^\prime)$, by the uniqueness
of the solution in $C^{2+2\alpha,1+\alpha}$ proved by Bronsard and
Reitich in~\cite{bronsard}.\\
At the maximal time $T^\prime$ it has to happen that 
$\psi^i_x$ is no more bounded away from zero or the $C^{2+2\alpha}$ norm
of $\psi^i(\cdot,t)$ is not bounded from above, but a simple computation
shows that then the same holds for $\theta^i$ hence also for $\eta^i$, which
is smooth and regular on $[0,T)$.\\
This clearly implies that $T^\prime=T$ and that 
the triods $\TTT^1_t$ and $\TTT^2_t$ are the
same subset of $\R^2$ for every time $t$ in the interval $[0,T)$.
\end{proof}

\begin{rem} 
This proposition clearly sets positively the question about the {\em
  geometric} uniqueness of the flow of a smooth triod.\\
Actually, we do not know if, at least in this special initial case,
uniqueness holds also in the class of smooth Brakke flows with
equality.

Clearly, if a triod is geometrically smooth but not 
smooth, we can reparametrize it and apply
Theorem~\ref{smoothexist}  in order to get a smooth flow (which is a
smooth Brakke flow with equality).
\end{rem}

Now in order to improve these results and to study the global
existence and regularity of the evolution in the next sections, 
we work out a priori estimates for $k^i$, $\lambda^i$ and their
derivatives.

\begin{rem}
Sometimes we will consider time depending functions, defined as the
maximum of a geometric quantity over the triods, in order to get
estimates by means of ODE's and maximum principle arguments. 
Even if the evolution is smooth, such functions will be typically only
Lipschitz, hence they can fail to be differentiable at some times, 
so there will be a little misuse of notation in writing a derivative
that possibly does not exist at {\em every} time. However, the 
arguments used, which are pointwise and apparently affected by the 
lack of differentiability, still work also in this situation, 
as explained in details by Hamilton 
in~\cite[Sections~3 and~4]{hamilton2}.
\end{rem}

We fix some non standard notation for the computations in the
sequel.\\
We denote with $\pol_\sigma(\ders^j\lambda,\ders^h k)$ a 
polynomial in $\lambda,\dots,\ders^j\lambda$ and
$k,\dots,\ders^h k$ with constant coefficients, 
such that every monomial it contains is of the
form
$$
C\prod_{l=0}^j (\ders^l\lambda)^{\alpha_l} \cdot \prod_{l=0}^h
(\ders^lk)^{\beta_l}\,
\text{ { with} $\,\, \sum_{l=0}^j  (l+1)\alpha_l + \sum_{l=0}^h
  (l+1)\beta_l = \sigma$,}
$$
we will call $\sigma$ the {\em geometric order} of $\pol_\sigma$.\\
Moreover, if one of the two arguments of $\pol_\sigma$ 
does not appear, it means that the polynomial does not contain it, for
instance, $\pol_\sigma(\ders^h k)$
does not contain neither $\lambda$ nor its derivatives.\\
We denote with
$\qol_\sigma(\dert^j\lambda,\ders^h k)$ a polynomial as before in
$\lambda,\dots,\dert^j\lambda$ and $k,\dots,\ders^h k$ such that all its
monomials are of the form
$$
C\prod_{l=0}^j (\dert^l\lambda)^{\alpha_l} \cdot \prod_{l=0}^h
(\ders^lk)^{\beta_l}\,
\text{ { with} $\,\, \sum_{l=0}^j  (2l+1)\alpha_l + \sum_{l=0}^h
  (l+1)\beta_l = \sigma$.}
$$
Finally, when we will write $\pol_\sigma(\vert\ders^j\lambda\vert,
\vert\ders^h k\vert)$ (or $\qol_\sigma(\vert\dert^j\lambda\vert,
\vert\ders^h k\vert)$) we will mean a finite sum of terms like
$$
C\prod_{l=0}^j \vert\ders^l\lambda\vert^{\alpha_l} \cdot \prod_{l=0}^h
\vert\ders^lk\vert^{\beta_l}\,
\text{ { with} $\,\, \sum_{l=0}^j  (l+1)\alpha_l + \sum_{l=0}^h
  (l+1)\beta_l = \sigma$,}
$$
where $C$ is a positive constant and the exponents $\alpha_l,\beta_l$
are non negative {\em real} values (analogously for $\qol_\sigma$).\\
Clearly we have $\pol_\sigma(\ders^j\lambda,\ders^h k)\leq
\pol_\sigma(\vert\ders^j\lambda\vert,\vert\ders^h k\vert)$.  

\begin{rem}
We advise the reader that in the following computations these
polynomials can vary from one line to another,
by addition of similar terms, what has to be kept in mind is that
the coefficients and the number of monomials they contains are
independent of $k$, $\lambda$ and their derivatives, since they arise 
by the algorithmic construction of the polynomials.\\
We will often denote with $C$ a generic constant which also can vary
from one passage to another.
\end{rem}
We will make extensive use of Young inequality in the following form
$$
ab\leq \varepsilon a^p+ C(\varepsilon,p,q)b^q\qquad \text{
    { for} $a,b,\varepsilon<0$, $p,q\in(1,\infty)$ { and} $1/p+1/q=1$.}
$$

\begin{lemma}\label{kexpr}
The following formulas hold
\begin{equation}
\begin{array}{ll}\dert\ders^jk=\ders^{j+2}k+\lambda\ders^{j+1}k +
\pol_{j+3}(\ders^{j}k)\qquad& \text{{ for every} $j\in \NN$,}\label{kappas}\\
\ders^jk=\dert^{j/2}k+\qol_{j+1}(\dert^{j/2-1}\lambda,\ders^{j-1}k)\qquad
&\text{{ if}  $j\geq2$ { is even,}}\\
\ders^jk=\dert^{(j-1)/2}k_s+\qol_{j+1}(\dert^{(j-3)/2}\lambda,\ders^{j-1}k)\qquad
&\text{{ if}  $j\geq1$ { is odd.}}
\end{array}
\end{equation}
\end{lemma}
\begin{proof} The case $j=0$ of the first formula is
  equation~\eqref{dertdik}. Suppose that the formula holds for
  $(j-1)$, using the commutation rule~\eqref{commut} we have
\begin{align*}
\dert\ders^j k=&\, \ders\dert\ders^{j-1}k + (k^2-\lambda_s)\ders^j k\\
=&\, \ders [\ders^{j+1}k + \lambda\ders^jk+ \pol_{j+2}(\ders^{j-1}k)]
-\lambda_s\ders^j k +\pol_{j+3}(\ders^jk)\\
=&\, \ders^{j+2}k + \lambda_s\ders^jk+ \lambda\ders^{j+1}k+
\pol_{j+3}(\ders^{j}k)-\lambda_s\ders^j k +\pol_{j+3}(\ders^jk)\\
=&\, \ders^{j+2}k + \lambda\ders^{j+1}k+\pol_{j+3}(\ders^{j}k)
\end{align*}
which gives the inductive step.\\
The second formula also follows by induction. The case $j=2$ is again
equation~\eqref{dertdik}. If the case $(j-2)$ holds, then by the first
formula,
\begin{align*}
\ders^jk=&\, \dert\ders^{j-2}k -\lambda\ders^{j-1}k + \pol_{j+1}(\ders^{j-2}k)=
\dert[\dert^{j/2-1}k+\qol_{j-1}(\dert^{j/2-2}\lambda,\ders^{j-3}k)]\\
=&\, \dert^{j/2}k+\dert\qol_{j-1}(\dert^{j/2-2}\lambda,\ders^{j-3}k)\,.
\end{align*}
Now, when we differentiate in $t$ the term
$\qol_{j-1}(\dert^{j/2-2}\lambda,\ders^{j-3}k)$
we will get a polynomial in
$\lambda,\dots,\dert^{j/2-1}\lambda, k,\dots\ders^{j-1}k$ and time derivatives
of space derivatives of $k$. Using the first formula we can express
these latter as polynomials in $\lambda$ and space derivatives of $k$, up to the
order $\ders^{j-1}k$. Moreover, it is easy to check that the resulting
polynomial is of the form
$\qol_{j+1}(\dert^{j/2-1}\lambda,\ders^{j-1}k)$, hence the formula for
$j$ is proved.\\
The odd case is analogous.
\end{proof}

\begin{lemma}\label{lambexpr}
The following formulas hold
\begin{equation*}
\begin{array}{ll}
\dert\ders^j\lambda=\ders^{j+2}\lambda-\lambda\ders^{j+1}\lambda -2k\ders^{j+1}k+
\pol_{j+3}(\ders^j\lambda,\ders^jk)
&\qquad\text{{ for every} $j\in \NN$,}\\
\ders^j\lambda=\dert^{j/2}\lambda+\pol_{j+1}(\ders^{j-1}\lambda,\ders^{j-1}k)
&\qquad\text{{ if}  $j\geq2$ { is even,}}\\
\ders^j\lambda=\dert^{(j-1)/2}\lambda_s+\pol_{j+1}(\ders^{j-1}\lambda,\ders^{j-1}k)
&\qquad\text{{ if}  $j\geq1$ { is odd.}}
\end{array}
\end{equation*}
\end{lemma}
\begin{proof} The case $j=0$ of the first formula is
  equation~\eqref{dertdilamb}, then the proof follows as
  for $k$ in the previous lemma.
\end{proof}

\begin{rem}\label{qolpol}
We state the following {\em calculus rules} which will be used
extensively in the sequel,
\begin{align*}
\pol_\alpha(\ders^j\lambda,\ders^hk)\cdot
\pol_\beta(\ders^l\lambda,\ders^mk)=&\,
\pol_{\alpha+\beta}(\ders^{\max\{j,l\}}\lambda,\ders^{\max\{h,m\}}k)\,,\\
\qol_\alpha(\dert^j\lambda,\ders^hk)\cdot
\qol_\beta(\dert^l\lambda,\ders^mk)=&\,
\qol_{\alpha+\beta}(\dert^{\max\{j,l\}}\lambda,\ders^{\max\{h,m\}}k)\,.
\end{align*}
Since the time derivatives of $k$ and $\lambda$ and their space
derivatives can be expressed in terms of these latter, by means of 
Lemmas~\ref{kexpr} and~\ref{lambexpr}, we have 
\begin{align*}
\ders^l\pol_\alpha(\ders^j\lambda,\ders^hk)=
\pol_{\alpha+l}(\ders^{j+l}\lambda,\ders^{h+l}k)\,,\qquad&
\dert^l\pol_\alpha(\ders^j\lambda,\ders^hk)=
\pol_{\alpha+2l}(\ders^{j+2l}\lambda,\ders^{h+2l}k)\\
\dert^l\qol_\alpha(\dert^j\lambda,\ders^hk)=
\qol_{\alpha+2l}(\dert^{j+l}\lambda,\ders^{h+2l}k)\,,\qquad&
\qol_\alpha(\dert^j\lambda,\ders^hk)=\pol_{\alpha}
(\ders^{2j}\lambda,\ders^{\max\{h,2j-1\}}k)\,.
\end{align*}
Moreover, by relations~\eqref{lambdakappa}, 
at the 3--point $\dert^j\lambda^i=({\mathrm S}\dert^j {\mathrm K})^i$,
that is, the time derivatives of $\lambda^i$ are expressible as time
derivatives of the functions $k^i$. Then, by using repeatedly such
relation and the first formula of Lemma~\ref{kexpr}, we can express
these latter as space derivatives of $k^i$. Hence, we have the
relation
$$
\sum_{i=1}^3\qol_\sigma(\dert^j\lambda^i,\ders^h
k^i)\,\biggr\vert_{\text{{ at    the 3--point}}}=
\pol_\sigma(\ders^{\max\{2j, h\}}
{\mathrm K})\,\biggr\vert_{\text{{ at
    the 3--point}}}
$$
with the meaning that this last polynomial contains also product of
derivatives of different $k^i$'s, because of the action of the linear
operator ${\mathrm S}:\R^3\to\R^3$.\\
We will often make use of this identity in the computations in the
sequel in the following form,
$$
\sum_{i=1}^3\qol_\sigma(\dert^j\lambda^i,\ders^h
k^i)\,\biggr\vert_{\text{{ at    the 3--point}}}\leq
\Vert \pol_\sigma(\vert\ders^{\max\{2j, h\}}k\vert)\Vert_{L^\infty}\,.
$$
\end{rem}

Before proceeding we prove also a relation holding at the end
points. 

\begin{lemma}\label{evenly}
At the three end points $P^i$ there holds
$\ders^jk^i=\ders^j\lambda^i=0$, for every even $j\in\NN$.
\end{lemma}
\begin{proof}
The first case $j=0$ simply follows from the fact that the velocity
$\underline{v}=\lambda\tau+k\nu$ is always zero at the three fixed end
points $P^i$.\\
We argue by induction, we suppose that
for every even natural $l\leq j-2$ we have
$\ders^lk^i=\ders^l\lambda^i=0$, then, by using the first equation in
Lemma~\ref{kexpr}, we get
$$
\ders^{j}k^i=\dert\ders^{j-2}k^i-\lambda^i\ders^{j-1}k^i -\pol_{j+1}(\ders^{j-2}k^i)
$$
at the points $P^i$.\\
We already know that $\lambda^i=0$ and by the
inductive hypothesis $\ders^{j-2}k^i=0$, thus
$\dert\ders^{j-2}k^i=0$. 
Since $\pol_{j+1}(\ders^{j-2}k^i)$ is a sum of terms like $C\prod_{l=0}^{j-2}
(\ders^lk^i)^{\alpha_l}$ with $\sum_{l=0}^{j-2}(l+1)\alpha_l=j+1$ which is
odd, at least one of the terms of this sum has to be odd, hence at
least for one index $l$, the product $(l+1)\alpha_l$ is odd. It
follows that at least for one even $l$ the exponent $\alpha_l$ is
nonzero. Hence, at least one even derivatives is present in every monomial
of $\pol_{j+1}(\ders^{j-2}k^i)$, which contains only derivatives up to the
order $(j-2)$.\\
Again, by the inductive hypothesis we then conclude that at
the end points $\ders^{j}k^i=0$.\\
We can deal with $\lambda^i$ similarly, by
means of the first equation of Lemma~\ref{lambexpr}.
\end{proof}

Taking into account that the time derivative of
the measure $ds$ is given by $(\lambda_s-k^2)\,ds$ and using the first
relation of Lemma~\ref{kexpr}, we compute for $j\in\NN$
\begin{align}\label{evolint000}
\frac{d\,}{dt} \int_{{\TTT_t}} |\ders^j k|^2\,ds=
&\, 2\int_{{\TTT_t}} \ders^j k\, \dert\ders^jk\,ds + \int_{{\TTT_t}}
|\ders^j k|^2(\lambda_s -k^2)\,ds\\
=&\, 2\int_{{\TTT_t}} \ders^j k\, \ders^{j+2}k +\lambda \ders^{j+1}
k\,\ders^{j} k +\pol_{j+3}(\ders^jk)\,\ders^jk\,ds 
+\int_{{\TTT_t}}
|\ders^j k|^2(\lambda_s -k^2)\,ds\nonumber\\
=&\, -2\int_{{\TTT_t}} \vert\ders^{j+1} k\vert^2\,ds + \int_{{\TTT_t}}
\ders(\lambda|\ders^j k|^2)\,ds +\int_{{\TTT_t}}\pol_{2j+4}(\ders^jk)\,ds\nonumber\\
&\, - 2\sum_{i=1}^3 \ders^j k^i\,\ders^{j+1}k^i\,
\biggr\vert_{\text{{ at the 3--point}}}
+ 2\sum_{i=1}^3 \ders^j k^i\,\ders^{j+1}k^i\,
\biggr\vert_{\text{{ at the point} $P^i$}}\nonumber\\
=&\, -2\int_{{\TTT_t}} \vert\ders^{j+1}k\vert^2\,ds  + \int_{\TTT_t}
\pol_{2j+4}(\ders^{j}k)\,ds\nonumber\\
&\,- \sum_{i=1}^3 2\ders^j
k^i\,\ders^{j+1}k^i+\lambda^i\vert\ders^{j}k^i\vert^2\,
\biggr\vert_{\text{{ at      the 3--point}}}\nonumber
\end{align}
where we integrated by parts a couple of times and we eliminated the
contributions given by the end points $P^i$ by means of
Lemma~\ref{evenly}.

In the very special (and important as we will see) case $j=0$ we get
explicitly
$$
\frac{d\,}{dt} \int_{{\TTT_t}} k^2\,ds
= -2\int_{{\TTT_t}} \vert k_s\vert^2\,ds  + \int_{\TTT_t} k^4\,ds 
- \sum_{i=1}^3 2 k^ik^i_s+\lambda^i(k^i)^2\,
\biggr\vert_{\text{{ at      the 3--point}}}\,.
$$
Then, recalling relation~\eqref{orto} with $l, m=0$, we have 
$\sum_{i=1}^3 k^ik^i_s+\lambda^i(k^i)^2\,\bigr\vert_{\text{{ at      the
      3--point}}}=0$, and substituting in
the last term above,
\begin{equation}\label{ksoltanto}
\frac{d\,}{dt} \int_{{\TTT_t}} k^2\,ds
= -2\int_{{\TTT_t}} \vert k_s\vert^2\,ds  + \int_{\TTT_t} k^4\,ds 
+\sum_{i=1}^3 \lambda^i(k^i)^2\,\biggr\vert_{\text{{ at      the 3--point}}}
\end{equation}
hence, we lowered the maximum order of the space derivatives of the
curvature in the 3--point term, particular now it is lower than the
one of the ``nice'' negative integral.

Now we are going to do the same for the general case, when $j\geq2$ is
even.\\
By means of formulas~\eqref{kappas} we have
\begin{align*}
2\ders^jk &\,\,\ders^{j+1} k + \lambda |\ders^j
k|^2\\
=&\,2[\dert^{j/2}k+\qol_{j+1}(\dert^{j/2-1}\lambda,\ders^{j-1}k)]\cdot
[\dert^{j/2}k_s+\qol_{j+2}(\dert^{j/2-1}\lambda,\ders^{j}k)]
+\qol_{2j+3}(\lambda,\ders^jk)\\
=&\,2[\dert^{j/2}k+\qol_{j+1}(\dert^{j/2-1}\lambda,\ders^{j-1}k)]\cdot
[\dert^{j/2}(k_s+k\lambda)+\qol_{j+2}(\dert^{j/2}\lambda,\ders^{j}k)]
+\qol_{2j+3}(\lambda,\ders^jk)\\
=&\,2\dert^{j/2}k\cdot\dert^{j/2}(k_s+k\lambda)
+ \dert^{j/2}k\cdot\qol_{j+2}(\dert^{j/2}\lambda,\ders^{j}k)
+\qol_{j+1}(\dert^{j/2-1}\lambda,\ders^{j-1}k)\cdot\dert^{j/2}(k_s+k\lambda)\\
&\,+\qol_{j+1}(\dert^{j/2-1}\lambda,\ders^{j-1}k)
\cdot\qol_{j+2}(\dert^{j/2}\lambda,\ders^{j}k)
+\qol_{2j+3}(\lambda,\ders^jk)\\
=&\,2\dert^{j/2}k\cdot\dert^{j/2}(k_s+k\lambda)+
\qol_{j+1}(\dert^{j/2-1}\lambda,\ders^{j-1}k)\cdot\dert^{j/2}k_s
+ \qol_{2j+3}(\dert^{j/2}\lambda,\ders^{j}k)\,.
\end{align*}
We now examine the term
$\qol_{j+1}(\dert^{j/2-1}\lambda,\ders^{j-1}k)\cdot\dert^{j/2}k_s$,
which contains $(j+1)$--th space derivatives of $k$
(after expansion of the $j/2$--th time derivative of $k_s$).\\
By using the third relation of Lemma~\ref{kexpr}, it can be written as
\begin{align*}
\qol_{j+1}(\dert^{j/2-1}\lambda,\ders^{j-1}k)&\,\cdot
\dert[\ders^{j-1}k+\qol_{j}(\dert^{j/2-2}\lambda,\ders^{j-2}k)]\\
&=\,\qol_{j+1}(\dert^{j/2-1}\lambda,\ders^{j-1}k)\cdot
[\dert\ders^{j-1}k+\qol_{j+2}(\dert^{j/2-1}\lambda,\ders^{j}k)]\\
&=\,\qol_{j+1}(\dert^{j/2-1}\lambda,\ders^{j-1}k)\cdot
\dert\ders^{j-1}k+ \qol_{2j+3}(\dert^{j/2-1}\lambda,\ders^jk)\,,
\end{align*}
moreover, if we look
at the polynomial $\qol_{j+1}(\dert^{j/2-1}\lambda,\ders^{j-1}k)$, we
can see that among its monomials, only those of the 
form $A\lambda\ders^{j-1}k$ or $Bk\ders^{j-1}k$ can contain the
derivative $\ders^{j-1}k$ (because of the {\em geometric order} of
$\pol_{2k+1}$). Hence,
\begin{align*}
\qol_{j+1}&(\dert^{j/2-1}\lambda,\ders^{j-1}k)\cdot\dert^{j/2}k_s\\
=&\, [\qol_{j+1}(\dert^{j/2-1}\lambda,\ders^{j-2}k)+A\lambda\ders^{j-1}k +
Bk\ders^{j-1}k]\cdot\dert\ders^{j-1}k+\qol_{2j+3}(\dert^{j/2-1}\lambda,\ders^jk)\\
=&\,(A\lambda\ders^{j-1}k + Bk\ders^{j-1}k)\cdot\dert\ders^{j-1}k +
\qol_{2j+3}(\dert^{j/2-1}\lambda,\ders^{j}k)\\
&\,+\dert[\qol_{j+1}(\dert^{j/2-1}\lambda,\ders^{j-2}k)\ders^{j-1}k]
-\dert\qol_{j+1}(\dert^{j/2-1}\lambda,\ders^{j-2}k)\ders^{j-1}k\\
=&\,(A\lambda+Bk)\dert(\ders^{j-1}k/2)^2
+\dert\qol_{2j+1}(\dert^{j/2-1}\lambda,\ders^{j-1}k)
+ \qol_{2j+3}(\dert^{j/2}\lambda,\ders^{j}k)\\
=&\,\dert[(A\lambda+Bk)(\ders^{j-1}k/2)^2] - (A\dert\lambda+B\dert k)(\ders^{j-1}k/2)^2\\
&\,+\dert\qol_{2j+1}(\dert^{j/2-1}\lambda,\ders^{j-1}k)
+ \qol_{2j+3}(\dert^{j/2}\lambda,\ders^{j}k)\\
=&\,\dert\qol_{2j+1}(\dert^{j/2-1}\lambda,\ders^{j-1}k)
+ \qol_{2j+3}(\dert^{j/2}\lambda,\ders^{j}k)\, .
\end{align*}
It follows that
\begin{align*}
\sum_{i=1}^3 2\ders^j k^i\,\ders^{j+1}k^i&\,+\lambda^i|\ders^j
k^i|^2\lambda\, \biggr\vert_{\text{{ at the 3--point}}}\\
=&\,\sum_{i=1}^3 2\dert^{j/2}k^i\cdot\dert^{j/2}(k^i_s+k^i\lambda^i)\\
&\, \phantom{\sum_{i=1}^3}+\dert\qol_{2j+1}(\dert^{j/2-1}\lambda^i,\ders^{j-1}k^i)
+\qol_{2j+3}(\dert^{j/2}\lambda^i,\ders^{j}k^i)
\,\biggr\vert_{\text{{ at the 3--point}}}\\
=&\,\sum_{i=1}^3 \dert\qol_{2j+1}(\dert^{j/2-1}\lambda^i,\ders^{j-1}k^i)
+\qol_{2j+3}(\dert^{j/2}\lambda^i,\ders^{j}k^i)
\,\biggr\vert_{\text{{ at the 3--point}}}
\end{align*}
by relations~\eqref{orto}.\\
Resuming, if $j\geq2$ is even, we have
\begin{align}\label{pippo100}
\frac{d\,}{dt} \int_{{\TTT_t}} |\ders^j k|^2\,ds=&\,
-\int_{{\TTT_t}} \vert\ders^{j+1}k\vert^2\,ds  +
\int_{\TTT_t}\pol_{2j+4}(\ders^{j}k)\,ds\\
&\, + \sum_{i=1}^3 \dert\qol_{2j+1}(\dert^{j/2-1}\lambda^i,\ders^{j-1}k^i)
+\qol_{2j+3}(\dert^{j/2}\lambda^i,\ders^{j}k^i)
\,\biggr\vert_{\text{{ at the 3--point}}}\,.\nonumber
\end{align}
Now, the key tool to estimate the terms 
$\int_{\TTT_t}\pol_{2j+4}(\ders^{j}k)\,ds$ and
$\sum_{i=1}^3\qol_{2j+3}(\dert^{j/2}\lambda^i,\ders^{j}k^i)\,\bigr
\vert_{\text{{ at the 3--point}}}$ are
the following Gagliardo--Nirenberg interpolation inequalities
(see~\cite{adams,aubin0}, for instance).
\begin{prop} Let $\gamma$ be a smooth regular curve in
  $\R^2$ with finite length ${\mathrm L}$. If $u$ is a smooth function defined on
  $\gamma$ and $m\geq1$, $p\in[2,+\infty]$, we have the estimates
  \begin{equation}\label{int1}
    {\Vert\partial_s^n u\Vert}_{L^p}
    \leq C_{n,m,p}
      {\Vert\partial_s^m  u\Vert}_{L^2}^{\sigma}
      {\Vert u\Vert}_{L^2}^{1-\sigma}+
      \frac{B_{n,m,p}}{{\mathrm L}^{m\sigma}}{\Vert u\Vert}_{L^2}
  \end{equation}
  for every $n\in\{0,\dots, m-1\}$ where
$$
\sigma=\frac{n+1/2-1/p}{m}
$$
and the constants $C_{n,m,p}$ and $B_{n,m,p}$ are independent of $\gamma$.
\end{prop}

\begin{rem} We put in evidence the particular case $p=+\infty$,
\begin{equation}\label{int2}
    {\Vert\partial_s^n u\Vert}_{L^\infty}
    \leq C_{n,m}
      {\Vert\partial_s^m  u\Vert}_{L^2}^{\sigma}
      {\Vert u\Vert}_{L^2}^{1-\sigma}+
      \frac{B_{n,m}}{{\mathrm L}^{m\sigma}}{\Vert
        u\Vert}_{L^2}\qquad\text{ { with} }\quad \text{ $\sigma=\frac{n+1/2}{m}$.}
  \end{equation}
It clearly follows that for a family of curves with lengths
equibounded from below by some positive value, these inequalities
hold with uniform constants.
\end{rem}

Every monomial of $\pol_{2j+4}(\ders^{j}k)$ is of the form
$C\prod_{l=0}^j(\ders^lk)^{\alpha_l}$ with
$\sum_{l=0}^j(l+1)\alpha_l=2j+4$, then we estimate its integral by
means of H\"older inequality,
$$
C\int_{{\TTT_t}}\prod_{l=0}^j(\ders^lk)^{\alpha_l}\,ds
\leq C\prod_{l=0}^j\left(\int_{\TTT_t}\vert\ders^lk\vert^{\alpha_l\beta_l}\right)^{1/\beta_l}\,ds
=C\prod_{l=0}^j\Vert\ders^lk\Vert^{\alpha_l}_{L^{\alpha_l\beta_l}}
$$
where the exponents $\beta_l$ satisfy $\sum 1/\beta_l=1$
and $\alpha_l\beta_l>2$ for every $l\in\{0,\dots,j\}$ such that
$\alpha_l\not=0$. These conditions can be fulfilled choosing
$\beta_l=\frac{2j+4}{(l+1)\alpha_l}$, 
then $\alpha_l\beta_l=(2j+4)/(l+1)>2$
since $l\leq j$ and
$\sum 1/\beta_l=\sum_{l=0}^j(l+1)\alpha_l/(2j+4)=1$.\\
Notice that the constant $C$ depends only on the
{\em structure} of the polynomial $\pol_{2j+4}(\ders^jk)$, that is,
only on $j\in\NN$.\\
Putting $n=l$, $m=j+1$, $p=\alpha_l\beta_l$ and
$u=k$ in  inequality~\eqref{int1} we get
\begin{equation*}
  {\Vert\partial_s^l k\Vert}_{L^{\alpha_l\beta_l}}
  \leq C_l\left(
    {\Vert\partial_s^{j+1}  k\Vert}_{L^2}^{\sigma_l}
    {\Vert k\Vert}_{L^2}^{1-\sigma_l}+
    {\Vert k\Vert}_{L^2}\right)
\end{equation*}
with $\sigma_l=\frac{l+1/2-1/(\alpha_l\beta_l)}{j+1}$ for every
$l\in\{0,\dots,j\}$ and the constants $C_l$ depend only on the lengths
of the curves.\\
Hence, since the number of monomials of $\pol_{2j+4}(\ders^jk)$ 
depends only on $j\in\NN$,
\begin{align*}
\int_{{\TTT_t}}\pol_{2j+4}(\ders^jk)\,ds
\leq&\, C \prod_{l=0}^j\left(
    {\Vert\partial_s^{j+1}  k\Vert}_{L^2}
    + {\Vert k\Vert}_{L^2}\right)^{\sigma_l\alpha_l}
{\Vert k\Vert}_{L^2}^{(1-\sigma_l)\alpha_l}\\
\leq&\, C \left(    {\Vert\partial_s^{j+1}  k\Vert}_{L^2}
    + {\Vert k\Vert}_{L^2}\right)^{\sum_{l=0}^j\sigma_l\alpha_l}
{\Vert k\Vert}_{L^2}^{\sum_{l=0}^j(1-\sigma_l)\alpha_l}\, .
\end{align*}
Now we have
\begin{align*}
\sum_{l=0}^j\sigma_l\alpha_l&=\,\sum_{l=0}^j\alpha_l
\frac{l+1/2-1/(\alpha_l\beta_l)}{j+1}=\sum_{l=0}^j
\frac{(l+1/2)\alpha_l-1/\beta_l}{j+1}\\
&=\,\frac{-1+\sum_{l=0}^j(l+1)\alpha_l - 1/2\alpha_l}{j+1}
=\frac{2j+3-1/2\sum_{l=0}^j\alpha_l}{j+1}\\
&\leq\, \frac{2j+3-1/2\sum_{l=0}^j\alpha_l(l+1)/(j+1)}{j+1}
=\frac{2j+3-(2j+4)/2(j+1)}{j+1}\\
&=\, \frac{2j+3 -1 -1/(j+1)}{j+1}=
2 - \frac{1}{(j+1)^2}<2
\end{align*}
then, by Young inequality,
$$
\left(    {\Vert\partial_s^{j+1}  k\Vert}_{L^2}
    + {\Vert k\Vert}_{L^2}\right)^{\sum_{l=0}^j\sigma_l\alpha_l}
{\Vert k\Vert}_{L^2}^{\sum_{l=0}^j(1-\sigma_l)\alpha_l}
\leq \varepsilon \left(    {\Vert\partial_s^{j+1}  k\Vert}_{L^2}
    + {\Vert k\Vert}_{L^2}\right)^2
+C {\Vert  k\Vert}_{L^2}^{2\frac{\sum_{l=0}^j(1-\sigma_l)\alpha_l}
{2-\sum_{l=0}^j\sigma_l\alpha_l}}
$$
and this last exponent is equal to
\begin{align*}
{2\frac{\sum_{l=0}^j(1-\sigma_l)\alpha_l}
{2-\sum_{l=0}^j\sigma_l\alpha_l}}=&\,
{2\frac{\sum_{l=0}^j\alpha_l
-\frac{2j+3-1/2\sum_{l=0}^j\alpha_l}{j+1}}
{2- \frac{2j+3-1/2\sum_{l=0}^j\alpha_l}{j+1}}}\\
=&\, {2\frac{(j+1)\sum_{l=0}^j\alpha_l - 2j-3+1/2\sum_{l=0}^j\alpha_l}
{2j+2- 2j-3+1/2\sum_{l=0}^j\alpha_l}}\\
=&\, {2\frac{-2j-3 +(j+3/2)\sum_{l=0}^j\alpha_l}
{-1 +1/2\sum_{l=0}^j\alpha_l}}\\
=&\, 2(2j+3)\,.
\end{align*}
Choosing a value $\varepsilon>0$ small enough and controlling, via
interpolation again, the term $\Vert k\Vert_{L^2}^2$, we conclude
$$
\int_{{\TTT_t}} \pol_{2j+4}(\ders^jk)\,ds\leq 
1/4\int_{{\TTT_t}} \vert\ders^{j+1}k\vert^2\,ds
+ C\left(\int_{{\TTT_t}} k^2\,ds\right)^{2j+3} + C
$$
where the constant $C$ depends only on $j\in\NN$ and the lengths of the
curves of the triod.\\
The term $\sum_{i=1}^3\qol_{2j+3}(\dert^{j/2}\lambda^i,\ders^{j}k^i)\,\bigr
\vert_{\text{{ at the 3--point}}}$ can be estimated similarly. Taking into account
Remark~\ref{qolpol}, we have
$\sum_{i=1}^3\qol_{2j+3}(\dert^{j/2}\lambda^i,\ders^{j}k^i)\,\bigr
\vert_{\text{{ at the
      3--point}}}\leq\Vert\pol_{2j+3}(\vert\ders^{j}k\vert)\Vert_{L^\infty}$ and
this latter can be controlled with a sum of terms like 
$C\prod_{l=0}^{j}\Vert\ders^{l}k\Vert_{L^\infty}^{\alpha_l}$
with $\sum_{l=0}^j(l+1)\alpha_l=2j+3$.\\
Then, we use interpolation inequalities with $p=+\infty$,
$$
\Vert\ders^lk\Vert_{L^\infty}\leq C_l\left(
    {\Vert\partial_s^{j+1}  k\Vert}_{L^2}^{\sigma_l}
    {\Vert k\Vert}_{L^2}^{1-\sigma_l}+
    {\Vert k\Vert}_{L^2}\right)
$$
with $\sigma_l=\frac{l+1/2}{j+1}$, hence
\begin{align*}
\sum_{i=1}^3\qol_{2j+3}(\dert^{j/2}\lambda^i,\ders^{j}k^i)\,\biggr
\vert_{\text{{ at the 3--point}}}
\leq&\, C \prod_{l=0}^j\left(
    {\Vert\partial_s^{j+1}  k\Vert}_{L^2}
    + {\Vert k\Vert}_{L^2}\right)^{\sigma_l\alpha_l}
{\Vert k\Vert}_{L^2}^{(1-\sigma_l)\alpha_l}\\
\leq&\, C \left(    {\Vert\partial_s^{j+1}  k\Vert}_{L^2}
    + {\Vert k\Vert}_{L^2}\right)^{\sum_{l=0}^j\sigma_l\alpha_l}
{\Vert k\Vert}_{L^2}^{\sum_{l=0}^j(1-\sigma_l)\alpha_l}
\end{align*}
and
\begin{align*}
\sum_{l=0}^j\sigma_l\alpha_l&=\,\sum_{l=0}^j\alpha_l
\frac{l+1/2}{j+1}=\sum_{l=0}^j
\frac{(l+1/2)\alpha_l}{j+1}\\
&=\,\frac{\sum_{l=0}^j(l+1)\alpha_l - 1/2\alpha_l}{j+1}
=\frac{2j+3-1/2\sum_{l=0}^j\alpha_l}{j+1}\\
&\leq\, \frac{2j+3-1/2\sum_{l=0}^j\alpha_l(l+1)/(j+1)}{j+1}
=\frac{2j+3-(2j+3)/2(j+1)}{j+1}\\
&=\, \frac{2j+3 -1 -1/2(j+1)}{j+1}=
2 - \frac{1}{2(j+1)^2}<2\,.
\end{align*}
As before, by Young inequality,
$$
\left(    {\Vert\partial_s^{j+1}  k\Vert}_{L^2}
    + {\Vert k\Vert}_{L^2}\right)^{\sum_{l=0}^j\sigma_l\alpha_l}
{\Vert k\Vert}_{L^2}^{\sum_{l=0}^j(1-\sigma_l)\alpha_l}
\leq \varepsilon \left(    {\Vert\partial_s^{j+1}  k\Vert}_{L^2}
    + {\Vert k\Vert}_{L^2}\right)^2
+ C {\Vert  k\Vert}_{L^2}^{2\frac{\sum_{l=0}^j(1-\sigma_l)\alpha_l}
{2-\sum_{l=0}^j\sigma_l\alpha_l}}
$$
and the last exponent is again equal to $2(2j+3)$. Choosing here also
a value $\varepsilon>0$ small enough, we get an estimate analogous to
the previous one.

Hence, for every even $j\geq2$ we can finally write
\begin{align}
\label{final1} \frac{d\,}{dt} \int_{{\TTT_t}} |\ders^j k|^2\,ds\leq
&\,-1/2\int_{{\TTT_t}} \vert\ders^{j+1}k\vert^2\,ds
+ C\left(\int_{{\TTT_t}} k^2\,ds\right)^{2j+3}  + C\\
&\, + \dert\sum_{i=1}^3
\qol_{2j+1}(\dert^{j/2-1}\lambda^i,\ders^{j-1}k^i)
\,\biggr\vert_{\text{{ at the 3--point}}}\nonumber\\
\leq &\,C\left(\int_{{\TTT_t}} k^2\,ds\right)^{2j+3}  + C+ \dert\sum_{i=1}^3
\qol_{2j+1}(\dert^{j/2-1}\lambda^i,\ders^{j-1}k^i)
\,\biggr\vert_{\text{{ at the 3--point}}}\,.\nonumber
\end{align}
Recalling the computation in the special case special case $j=0$, this
argument gives the same final estimate
without the the last term.
\begin{equation}\label{evolint999}
\left\vert\frac{d\,}{dt} \int_{{\TTT_t}} k^2\,ds\right\vert\leq 
C\left(\int_{{\TTT_t}} k^2\,ds\right)^{3} + C\,.
\end{equation}

Integrating~\eqref{final1} in time on $[0,t]$ and estimating we get
\begin{align*}
\int_{\TTT_t} |\ders^j k|^2\,ds\leq&\,
\int_{\TTT_0} |\ders^j k|^2\,ds +
C\int_{0}^t\left(\int_{\TTT_\xi} k^2\,ds\right)^{2j+3}\,d\xi
+ Ct\\
&\, + \sum_{i=1}^3
\qol_{2j+1}(\dert^{j/2-1}\lambda^i(0,t),\ders^{j-1}k^i(0,t))
-\qol_{2j+1}(\dert^{j/2-1}\lambda^i(0,0),\ders^{j-1}k^i(0,0)) \\
\leq&\, C\int_{0}^t\left(\int_{\TTT_\xi}
  k^2\,ds\right)^{2j+3}\,d\xi + \Vert\pol_{2j+1}(\vert\ders^{j-1}k\vert)\Vert_{L^\infty}
+ Ct + C
\end{align*}
where in the last passage we used as before Remark~\ref{qolpol}. 
The constant $C$ depends only on $j\in\NN$ and the triod $\TTT_0$.\\
Interpolating again by means of inequalities~\eqref{int2} (we
leave the details to the reader),  one gets
$$
\Vert\pol_{2j+1}(\vert\ders^{j-1}k\vert)\Vert_{L^\infty}\leq
1/2\Vert\partial_s^{j}k\Vert^2_{L^2}+
      C\Vert k \Vert_{L^2}^{4j+2}\,.
$$
Hence, putting all together, for every even $j\in\NN$.
$$
\int_{\TTT_t} |\ders^j k|^2\,ds
\leq C\int_{0}^t\left(\int_{\TTT_\xi}
  k^2\,ds\right)^{2j+3}\,d\xi + C\left(\int_{\TTT_t}
  k^2\,ds\right)^{2j+1} + Ct + C\,.
$$
Passing from integral to $L^\infty$ estimates by using
inequalities~\eqref{int2}, we have the following proposition.

\begin{prop}\label{pluto1000} If the lengths of the three curves are
  positively bounded from below and the $L^2$ norm of
  $k$ is bounded, uniformly on $[0,T)$, then the curvature of $\TTT_t$ 
  and all its space derivatives of are uniformly bounded in the same
  time interval by some constants depending only on the $L^2$ integrals of
  the space derivatives of $k$ on the initial triod $\TTT_0$.
\end{prop}

Now in the hypotheses of this proposition we deal with $\lambda$ and
its derivatives.\\
At the 3--point
$\sum_{i=1}^3(\lambda^i)^2=\sum_{i=1}^3(k^i)^2$, hence
the squared modulus of the velocity
$v^2=\vert \underline{v}\vert^2$ is uniformly bounded at $O$.\\
Then, since $\gamma^i_t(1,t)=\underline{v}^i(1,t)=0$
for every index $i\in\{1,2,3\}$, by the maximum principle applied to the equation for $v^2$,
\begin{equation*}
\dert v^2 =(v^2)_{ss} -2\lambda^2_s -2k^2_s -\lambda(v^2)_s + 2v^2 k^2
\end{equation*}
which follows from equation~\eqref{dertdilamb}
\begin{equation*}
\dert\lambda =\lambda_{ss} -\lambda\lambda_s - 2kk_s +\lambda k^2
\end{equation*}
and equation~\eqref{dertdik}, we see that if $v^2$ gets larger than
its value at the 3--point, then its maximum is taken in the interior of
some curve of the triod so, as $k^2$ is uniformly bounded,
\begin{equation*}
\dert v^2_{\max} \leq 2v^2_{\max} k^2_{\max}\leq Cv^2_{\max}\,.
\end{equation*}
Integrating this {\em linear} differential inequality, we obtain that
$\underline{v}$ and hence $\lambda$ are also uniformly bounded as $k$
and its derivatives in the time interval $[0,T)$.

By Lemma~\ref{lambexpr} and computing like for $k$, we get
\begin{align}\label{levole}
\frac{d\,}{dt} \int_{{\TTT_t}} |\ders^j \lambda|^2\,ds
=&\, -2\int_{{\TTT_t}} \vert\ders^{j+1}\lambda\vert^2\,ds  -2
\int_{\TTT_t} \lambda\ders^j\lambda\ders^{j+1}\lambda
+2 k\ders^j\lambda \ders^{j+1}k\,ds\\
&\,+ \int_{\TTT_t}\pol_{2j+4}(\ders^{j}\lambda,\ders^{j}k)\,ds\nonumber\\
&\, - 2\sum_{i=1}^3 \ders^j \lambda^i\,\ders^{j+1}\lambda^i
\, \biggr\vert_{\text{{ at the 3--point}}}
+ 2\sum_{i=1}^3 \ders^j \lambda^i\,\ders^{j+1}\lambda^i
\, \biggr\vert_{\text{{ at the point} $P^i$}}\nonumber\\
\leq&\, -\int_{{\TTT_t}} \vert\ders^{j+1}\lambda\vert^2\,ds  +
\int_{\TTT_t}\pol_{2j+4}(\ders^{j}\lambda,\ders^j k)
\,ds + \int_{\TTT_t}\vert\ders^{j+1}k\vert^2\,ds\nonumber\\
&\, - 2\sum_{i=1}^3 \ders^j \lambda^i\,\ders^{j+1}\lambda^i
\, \biggr\vert_{\text{{ at the 3--point}}}\nonumber\\
\leq&\, -\int_{{\TTT_t}} \vert\ders^{j+1}\lambda\vert^2\,ds  +
\int_{\TTT_t}\rol_{2j+4}(\vert\ders^{j}\lambda\vert) + C 
- 2\sum_{i=1}^3 \ders^j \lambda^i\,\ders^{j+1}\lambda^i
\, \biggr\vert_{\text{{ at the 3--point}}}\nonumber
\end{align}
where in the first passage we used Peter--Paul inequality 
$ab\leq \varepsilon a^2+ b^2/4\varepsilon$ with $\varepsilon=1/2$,
$a=\ders^{j+1}\lambda$ and
$b=\lambda \ders^j\lambda$ on the first term of the second integral,
and with $\varepsilon=1/2$,
$a=\ders^{j+1}k$ and
$b=k\ders^j\lambda$ on the second term of the second integral. Then,
we summed and absorbed the terms without $(j+1)$--th derivatives 
into $\pol_{2j+4}(\ders^j\lambda,\ders^{j}k)$.\\
In the second passage,  using Young inequality, we ``separated'' in all the
monomials of $\pol_{2j+4}(\ders^j\lambda,\ders^jk)$ 
the derivatives of $\lambda$ and $k$, controlling them 
with $\rol_{2j+4}(\vert\ders^j\lambda\vert)\
+\rol_{2j+4}(\vert\ders^jk\vert)$, where $\rol_{2j+4}$ denotes
a ``polynomial'' (similar to $\pol$ and $\qol$) with {\em real
  exponents} all greater or equal than 1.\\
Then we estimated the term $\rol_{2j+4}(\vert\ders^jk\vert)$ with
some constant, as  we know that $k$ and its space derivatives
are bounded. Notice that the number of monomials of
$\pol_{2j+4}(\ders^j\lambda,\ders^jk)$ depends only on $j\in\NN$.\\
Finally, the contributions
of the end points $P^i$ vanish by Lemma~\ref{evenly}, since 
at least one of the two derivatives of $\lambda$ is even.

We estimate by interpolation, exactly like for $k$, the term 
$\int_{\TTT_t}\rol_{2j+4}(\vert\ders^{j}\lambda\vert)\,ds$
with a small fraction of the ``good'' term $\int_{{\TTT_t}}
\vert\ders^{j+1}\lambda\vert^2\,ds$ and a possibly large multiple of
$\left(\int_{{\TTT_t}}\lambda^2\,ds\right)^{2j+3}$, 
which is bounded by the argument above.\\
Hence, it only remains to control the 3--point term
$$
-2\sum_{i=1}^3 \ders^j \lambda^i\,\ders^{j+1}\lambda^i
\, \biggr\vert_{\text{{ at the 3--point}}}\leq
C\sum_{i=1}^3 \vert\ders^j
\lambda^i\vert\,\vert\ders^{j+1}\lambda^i\vert\biggr\vert_{\text{{ at
      the 3--point}}}\,.
$$
Now, if $j\in\NN$ is odd, by the second formula of Lemma~\ref{lambexpr},
$$
\ders^{j+1}\lambda^i=\dert^{(j+1)/2}\lambda^i+\pol_{j+2}(\ders^{j}\lambda^i,\ders^{j}k^i)
$$
thus, recalling that
$\vert\dert^{(j+1)/2}\lambda^i\vert\,\bigr\vert_{\text{{ at the
      3--point}}}\leq
\Vert\pol_{j+2}(\vert\ders^{j+1}k\vert)\Vert_{L^\infty}$ by Remark~\ref{qolpol}, 
\begin{align*}
\vert\ders^j
\lambda^i\vert\,\vert\ders^{j+1}\lambda^i\vert
\,\biggr\vert_{\text{{ at      the 3--point}}}
\leq&\,\vert\ders^j\lambda^i\vert\,\left(\vert\dert^{(j+1)/2}\lambda^i\vert+
\vert\pol_{j+2}(\ders^{j}\lambda^i,\ders^{j}k^i)\vert\right)
\,\biggr\vert_{\text{{ at      the 3--point}}}\\
\leq&\, \Vert\ders^j\lambda\Vert_{L^\infty}\,\Vert\pol_{j+2}(\vert\ders^{j+1}k\vert)\Vert_{L^\infty}
+\Vert\pol_{2j+3}(\vert\ders^{j}\lambda\vert,\vert\ders^{j}k\vert)\Vert_{L^\infty}\\
\leq&\, \Vert\ders^j\lambda\Vert_{L^\infty}^{\frac{2j+3}{j+1}}+
\Vert\pol_{j+2}(\vert\ders^{j+1}k\vert)\Vert_{L^\infty}^{\frac{2j+3}{j+2}} +
\Vert\pol_{2j+3}(\vert\ders^{j}\lambda\vert,\vert\ders^{j}k\vert)\Vert_{L^\infty}\\
\leq&\,\Vert\rol_{2j+3}(\vert\ders^{j}\lambda\vert)\Vert_{L^\infty}
+C+\Vert\pol_{2j+3}(\vert\ders^{j}\lambda\vert,\vert\ders^{j}k\vert)\Vert_{L^\infty}\\
\leq&\,\Vert\rol_{2j+3}(\vert\ders^{j}\lambda\vert)\Vert_{L^\infty}
+C+\Vert\rol_{2j+3}(\vert\ders^{j}\lambda\vert)+\rol_{2j+3}(\vert\ders^{j}k\vert)\Vert_{L^\infty}\\
\leq&\,\Vert\rol_{2j+3}(\vert\ders^{j}\lambda\vert)\Vert_{L^\infty}
+C
\end{align*}
where we used Young inequality and the fact that
$\Vert\pol_{j+2}(\vert\ders^{j+1}k\vert)\Vert_{L^\infty}$ is uniformly bounded.\\
Moreover, we separated, as before, the derivatives of $\lambda$ and $k$ in
every monomial of $\pol_{2j+3}(\ders^{j}\lambda,\ders^{j}k)$, hence
estimating them with $\rol_{2j+3}(\vert\ders^j\lambda\vert)\
+\rol_{2j+3}(\vert\ders^jk\vert)$. Finally, we controlled the $k$--terms with
some constants and we can now interpolate the $\lambda$--terms like we did for
$k$, since these latter do not contain $(j+1)$--th space derivatives of $\lambda$.

Hence, coming back to computation~\eqref{levole}, we conclude that for
every odd $j\in\NN$
$$
\frac{d\,}{dt} \int_{{\TTT_t}} |\ders^j \lambda|^2\,ds\leq C<+\infty
$$
for a constant $C$ depending only on $j\in\NN$ and $\TTT_0$.

Like for $k$, since we know that $\Vert\lambda\Vert_{L^2}$ is bounded,
passing from integral to $L^\infty$ estimates by means of
inequalities~\eqref{int2}, we obtain that also all the space
derivatives of $\lambda$ are uniformly
bounded in $[0,T)$.

Then, we can bound from above and positively from below the term
$\vert\gamma_x\vert$ at the denominator in the evolution
equation~\eqref{evoleq}.
$$
\dert\log\vert\gamma_x(x,t)\vert=
\frac{\langle\gamma_x\,\vert\,\gamma_{xt}\rangle}{\vert\gamma_x\vert^2}
=\langle\tau\,\vert\,\ders (\lambda\tau+k\nu)\rangle=\lambda_s-k^2\leq C<+\infty
$$
for a constant $C$ independent of $x\in[0,1]$ and $t\in[0,T)$.
This implies that $\vert\gamma_x\vert$ is bounded from above and
away from zero, uniformly in space and time as
$$
\Bigl\vert\log\vert\gamma_x(x,t)\vert\Bigr\vert\leq \int_0^t
\Bigl\vert\dert\log\vert\gamma_x(x,\xi)\vert\Bigr\vert\,d\xi \leq
Ct\leq CT<+\infty\,.
$$
Since $\vert\gamma_x\vert$ is uniformly bounded and 
$\partial_x=\vert\gamma_x\vert\ders$, by using the evolution
equation~\eqref{evoleq}, it follows that
all the mixed derivatives in $x$ and $t$ of $\gamma^i$ for every
$i\in\{1,2,3\}$ are uniformly bounded in $[0,1]\times[0,T)$.

\begin{prop}\label{unif222} If $\TTT_t$ is a smooth evolution of the
  initial triod $\TTT_0=\cup_{i=1}^3\sigma^i$ such that the lengths 
  of the three curves are uniformly bounded away from zero and the
  $L^2$ norm of the curvature is uniformly bounded by some constants
  in the time interval $[0,T)$, then
\begin{itemize}
\item all the derivatives in space and time of $k$ and $\lambda$ are uniformly
bounded in $[0,1]\times[0,T)$,
\item all the derivatives in space and time of the curves 
$\gamma^i(x,t)$ are uniformly bounded in $[0,1]\times[0,T)$,
\item the quantities $\vert\gamma^i_x(x,t)\vert$ are uniformly bounded
  from above and away from zero in $[0,1]\times[0,T)$.
\end{itemize}
All the bounds depend only on the uniform controls on $k$ and the
lengths of the curves, and on the $L^\infty$ norms of the derivatives
of the maps $\sigma^i$ composing the initial triod $\TTT_0$.
\end{prop}

Now, we work out a second family of estimates 
where everything is controlled only by the $L^2$ norm of the curvature
and the inverses of the lengths of the three curves at time zero.\\
As before we consider the smooth evolution $\TTT_t$ of a smooth triod $\TTT_0$
in the time interval $[0,T)$.

\begin{prop}\label{stimaL}
For every $M>0$ there exists a time $T_M\in (0,T)$ such that if the 
$L^2$ norm of the curvature and
the inverses of the lengths of the three curves of $\TTT_0$ are
bounded by $M$, then the $L^2$ norm of $k$ and the inverses of the
lengths of the curves of $\TTT_t$ are smaller than $2M^2+6M$, for every
time $t\in[0,T_M]$.
\end{prop}
\begin{proof}
The evolution equations for the lengths of the three curves are given by
$\frac{dL^i(t)}{dt}=-\lambda^i(0,t)-\int_{\gamma^i(\cdot,t)}k^2\,ds$
(Proposition~\ref{equality1000}), 
then recalling computation~\eqref{ksoltanto}, we have
\begin{align*}
\frac{d\,}{dt} \left(\int_{\TTT_t} k^2\,ds + \right. \,& \left.
\sum_{i=1}^3\frac{1}{L^i}\right)\leq -2\int_{{\TTT_t}}
k_{s}^2\,ds + \int_{{\TTT_t}} k^4\,ds +C\Vert k\Vert_{L^\infty}^3
- \sum_{i=1}^3\frac{1}{(L^i)^2}\frac{dL}{dt}\\
=&\, -2\int_{{\TTT_t}} k_{s}^2\,ds + \int_{{\TTT_t}} k^4\,ds
+C\Vert k\Vert^3_{L^\infty}
+\sum_{i=1}^3\frac{\lambda^i(0,t)+\int_{\gamma^i(\cdot,t)}k^2\,ds}{(L^i)^2}\\
\leq&\,-2\int_{{\TTT_t}} k_{s}^2\,ds + \int_{{\TTT_t}} k^4\,ds
+ C\Vert k\Vert_{L^\infty}^3
+C\sum_{i=1}^3\frac{\Vert k\Vert_{L^\infty}}{(L^i)^2}
+\sum_{i=1}^3\frac{\int_{\TTT_t}k^2\,ds}{(L^i)^2}\\
\leq&\,-2\int_{{\TTT_t}} k_{s}^2\,ds + \int_{{\TTT_t}} k^4\,ds
+ C\Vert k\Vert_{L^\infty}^3
+ C\left(\int_{\TTT_t}k^2\,ds\right)^3 + C\sum_{i=1}^3\frac{1}{(L^i)^3}
\end{align*}
where we used Young inequality in the last passage.\\
Interpolating as before (and applying again Young inequality) 
but keeping now in evidence the terms depending on $L^i$
in inequalities~\eqref{int1} and~\eqref{int2}, we obtain
\begin{align*}
\frac{d\,}{dt} \left(\int_{\TTT_t} k^2\,ds +
  \sum_{i=1}^3\frac{1}{L^i}\right)
\leq&\,-\int_{\TTT_t} k_{s}^2\,ds
+ C\left(\int_{\TTT_t}k^2\,ds\right)^3 +
C\sum_{i=1}^3\frac{\left(\int_{\TTT_t}k^2\,ds\right)^2}
{L^i}\\
&\,+C\sum_{i=1}^3\frac{\left(\int_{\TTT_t}k^2\,ds\right)^{3/2}}
{(L^i)^{3/2}}+ C\sum_{i=1}^3\frac{1}{(L^i)^3}\\
\leq&\, C\left(\int_{\TTT_t}k^2\,ds\right)^3
+ C\sum_{i=1}^3\frac{1}{(L^i)^3}\\
\leq&\, C\left(\int_{\TTT_t}k^2\,ds +
  \sum_{i=1}^3\frac{1}{L^i}\right)^3
\end{align*}
with a {\em universal} constant $C$ independent of the triods.\\
This means that the function $f(t)=\int_{\TTT_t} k^2\,ds +
  \sum_{i=1}^3\frac{1}{L^i(t)}$ satisfies the differential inequality
$f^\prime\leq Cf^3$, hence, after integration the thesis follows.
\end{proof}

By means of this proposition we can strengthen the conclusion of
Proposition~\ref{unif222}.

\begin{cor} In the hypothesis of the previous proposition, in the time
  interval $[0,T_M]$ all the bounds in Proposition~\ref{unif222} depends only
  on the $L^2$ norm of the curvature and the lengths of the curves of
  $\TTT_0=\cup_{i=1}^3\sigma^i$ and on the $L^\infty$ norms of
  the derivatives of the maps $\sigma^i$.
\end{cor}

Then, from now on we assume that the $L^2$ norm of the curvature
and the inverses of the lengths of the three curves are bounded 
in the interval $[0,T_M]$.

Before dealing with the general estimate, we work out a special
case as an example.\\
By means of computations~\eqref{evolint000}, we obtain
\begin{align*}
\frac{d\,}{dt} \int_{{\TTT_t}} k^2+ tk_s^2+\frac{t^2 k_{ss}^2}{2} \,ds=
&\, -2\int_{{\TTT_t}} k_s^2 + t k_{ss}^2 + \frac{t^2 k_{sss}^2}{2}
\,ds + \int_{{\TTT_t}} k_s^2+ t k_{ss}^2 \,ds\\
&\, + \int_{\TTT_t}
k^4 + t \pol_{6}(k_s)
+t^2 \pol_{8}(k_{ss})\,ds\\
&\, + \sum_{i=1}^3  \lambda^i(k^i)^2
- t[2k^i_s k^i_{ss}+\lambda^i(k_s^i)^2]\\
&\, \phantom{+ \sum_{i=1}^3 }  - t^2[2k^i_{ss}
k_{sss}^i+\lambda^i(k^i_{ss})^2]/2
\, \biggr\vert_{\text{{ at the 3--point}}}\\
=&\, -\int_{{\TTT_t}} k_s^2 + t k_{ss}^2 + t^2 k_{sss}^2
\,ds\\
&\, + \int_{\TTT_t}
k^4 + t \pol_{6}(k_s)
+t^2 \pol_{8}(k_{ss})\,ds\\
&\, + \sum_{i=1}^3  \lambda^i(k^i)^2
- t[2k^i_s k^i_{ss}+\lambda^i(k_s^i)^2]\\
&\, \phantom{+ \sum_{i=1}^3 }  - t^2[2k^i_{ss}
k_{sss}^i+\lambda^i(k^i_{ss})^2]/2
\, \biggr\vert_{\text{{ at the 3--point}}}\,.
\end{align*}
Estimating, as we did before in order to get inequality~\eqref{final1},
the terms coming from the integrals of $k^2$ and $t^2 k_{ss}^2$ (the {\em even}
terms) we get
\begin{align*}
\frac{d\,}{dt} \int_{{\TTT_t}} k^2+ tk_s^2+t^2 k_{ss}^2/2 \,ds\leq
&\, -1/2\int_{{\TTT_t}} k_s^2 + t k_{ss}^2 + t^2 k_{sss}^2
\,ds\\
&\, + \int_{\TTT_t}
t \pol_{6}(k_s)\,ds
- t\sum_{i=1}^3 2k^i_s k^i_{ss}+\lambda^i(k_s^i)^2
\, \biggr\vert_{\text{{ at the 3--point}}}\\
&\, +C\left(\int_{{\TTT_t}} k^2\,ds\right)^3
+Ct^2\left(\int_{{\TTT_t}} k^2\,ds\right)^7+C(1+t^2)\\
&\,+t^2\dert\sum_{i=1}^3\qol_5(\lambda^i, k^i_s)
\,\biggr\vert_{\text{{ at the 3--point}}}\,.
\end{align*}
The integral term $\int_{\TTT_t}t\pol_6(k_s)\,ds$ and the
term $-t\sum_{i=1}^3 \lambda^i(k_s^i)^2\, \biggr\vert_{\text{{ at the
    3--point}}}$, which is of type $t\qol_5(\lambda^i,k^i_s)$,
can be estimated by means of interpolation
inequalities like in the even derivative case with a small fraction of the
term $t\int_{\TTT_t}k_{ss}^2$ and a possibly large multiple 
of $t\left(\int_{\TTT_t}k^2\,ds\right)^5$. The same holds  for the
term $-2t\sum_{i=1}^3\qol_5(\lambda^i, k_s^i)\,
\biggr\vert_{\text{{ at the
    3--point}}}$ arising by putting $t^2$
inside the time derivative of the last sum. Finally
since the $L^2$ norm of $k$ is bounded by some constant and $t\leq
T_M$, we conclude
\begin{align*}
\frac{d\,}{dt} \int_{{\TTT_t}} k^2+ tk_s^2+t^2 k_{ss}^2/2 \,ds\leq
&\, -1/4\int_{{\TTT_t}} k_s^2 + t k_{ss}^2 + t^2 k_{sss}^2
\,ds + C \\
&\, + \dert\sum_{i=1}^3 t^2\qol_5(\lambda^i, k^i_s)
\,\biggr\vert_{\text{{ at the 3--point}}}
- 2t\sum_{i=1}^3 k^i_sk^i_{ss}
\, \biggr\vert_{\text{{ at the 3--point}}}\,.
\end{align*}
The last term is the only one which needs a special analysis, so we deal
with it in the general case.

Considering now $j\in\NN$ even and following exactly the same line,
if we differentiate the function
$$
\int_{{\TTT_t}} k^2+ tk_s^2+ \frac{t^2 k_{ss}^2}{2!}+\dots+
\frac{t^j\vert \partial_s^jk\vert^2}{j!}\,ds\,,
$$
and we estimate as above, we obtain
\begin{align}\label{ciaociao1}
\frac{d\,}{dt}\int_{{\TTT_t}}&\, k^2+ tk_s^2+ \frac{t^2 k_{ss}^2}{2!}+\dots+
\frac{t^j\vert \partial_s^jk\vert^2}{j!}\,ds\\
\leq&\,-\varepsilon\int_{{\TTT_t}} k_s^2 + t k_{ss}^2 + t^2 k_{sss}^2+\dots+
t^j\vert \partial_s^{j+1}k\vert^2\,ds + C\nonumber\\
&\, + \dert\sum_{i=1}^3 t^2\qol_5(\lambda^i, k_s^i)
+ t^4\qol_9(\dert\lambda^i,k_{sss}^i)
+\dots+ t^j\qol_{2j+1}(\dert^{j/2-1}\lambda^i,\ders^{j-1}k^i)
\,\biggr\vert_{\text{{ at the 3--point}}}\nonumber\\
&\, +C\sum_{i=1}^3 t k_s^i k^i_{ss}+
t^3 k^i_{sss} k_{ssss}^i +\dots +
t^{j-1}\ders^{j-1} k^i\,\ders^{j}k^i\, \biggr\vert_{\text{{ at
    the 3--point}}}\nonumber
\end{align}
in the time interval $[0,T_M]$, where $\varepsilon>0$ and 
$C$ are two constants depending only on the
$L^2$ norm of the curvature  and the inverses of the lengths of the
three curves of $\TTT_0$.

We deal now with the terms $\sum_{i=1}^3t^{h-1}\ders^{h-1}
k^i\ders^{h}k^i\, \biggr\vert_{\text{{ at
    the 3--point}}}$ with $h\in\{2,\dots,j\}$ even (notice that this
family contains the term $tk^i_sk^i_{ss}$ of the case $j=2$ above).
By formulas~\eqref{kappas} we have
\begin{align*}
\ders^{h-1} k\cdot&\,\ders^{h}k\\
=&\,[\dert^{h/2-1}k_s+\qol_{h}(\dert^{h/2-2}\lambda,\ders^{h-2}k)]\cdot
[\dert^{h/2}k+\qol_{h+1}(\dert^{h/2-1}\lambda,\ders^{h-1}k)]\\
=&\,[\dert^{h/2-1}(k_s+\lambda k)+\qol_{h}(\dert^{h/2-1}\lambda,\ders^{h-2}k)]\cdot
[\dert^{h/2}k+\qol_{h+1}(\dert^{h/2-1}\lambda,\ders^{h-1}k)]\\
=&\, \dert^{h/2-1}(k_s+\lambda k)\cdot\dert^{h/2}k
+ \dert^{h/2-1}(k_s+\lambda k)\cdot
\qol_{h+1}(\dert^{h/2-1}\lambda,\ders^{h-1}k)\\
&\,+\qol_{h}(\dert^{h/2-1}\lambda,\ders^{h-2}k)
\cdot\qol_{h+1}(\dert^{h/2-1}\lambda,\ders^{h-1}k)
+ \dert^{h/2}k\cdot  \qol_{h}(\dert^{h/2-1}\lambda,\ders^{h-2}k)\\
=&\, \dert^{h/2-1}(k_s+\lambda k)\cdot\dert^{h/2}k
+\qol_{2h+1}(\dert^{h/2-1}\lambda,\ders^{h-1}k)
+ \dert^{h/2}k\cdot  \qol_{h}(\dert^{h/2-1}\lambda,\ders^{h-2}k)\\
=&\, \dert^{h/2-1}(k_s+\lambda k)\cdot\dert^{h/2}k
+\qol_{2h+1}(\dert^{h/2-1}\lambda,\ders^{h-1}k)\\
&\,+ [\ders^{h}k+\qol_{h+1}(\dert^{h/2-1}\lambda,\ders^{h-1}k)]
\cdot  \qol_{h}(\dert^{h/2-1}\lambda,\ders^{h-2}k)\\
=&\, \dert^{h/2-1}(k_s+\lambda k)\cdot\dert^{h/2}k
+\qol_{2h+1}(\dert^{h/2-1}\lambda,\ders^{h-1}k)
+ \ders^{h}k\cdot  \qol_{h}(\dert^{h/2-1}\lambda,\ders^{h-2}k)
\end{align*}
and since we are summing at the 3--point, the first product is zero by
relations~\eqref{orto} and, by Remark~\ref{qolpol}, we get
\begin{align*}
\sum_{i=1}^3t^{h-1}\ders^{h-1} k^i\ders^{h}k^i\,
\biggr\vert_{\text{{ at
    the 3--point}}}
= &\,\sum_{i=1}^3t^{h-1}\qol_{2h+1}(\dert^{h/2-1}\lambda^i,\ders^{h-1}k^i)\\
 &\,\phantom{\sum_{i=1}^3} + t^{h-1}
\ders^{h}k^i\cdot  \qol_{h}(\dert^{h/2-1}\lambda^i,\ders^{h-2}k^i)\,
\biggr\vert_{\text{{ at
    the 3--point}}}\\
\leq &\,t^{h-1}\Vert \pol_{2h+1}(\vert\ders^{h-1}k\vert)\Vert_{L^\infty}
+t^{h-1}\Vert\ders^{h}k\Vert_{L^\infty}\Vert\pol_{h}(\vert\ders^{h-2}k\vert)\Vert_{L^\infty}\,.
\end{align*}
The term $t^{h-1}\Vert \pol_{2h+1}(\vert\ders^{h-1}k\vert)\Vert_{L^\infty}$
is controlled as before by a small fraction of the term
$t^{h-1}\int_{\TTT_t}\vert \ders^h k\vert^2\,ds$ and a possibly large
multiple of $t^{h-1}$ times some power of the $L^2$ norm of $k$ 
(which is bounded), whereas 
$t^{h-1}\Vert\ders^{h}k\Vert_{L^\infty}\Vert\pol_{h}(\vert\ders^{h-2}k\vert)\Vert_{L^\infty}$
is the critical term.

Again by means of interpolation inequalities~\eqref{int2} we have
\begin{align*}
\Vert \ders^h k\Vert_{L^\infty}\leq &\,C\Vert
\ders^{h+1}k\Vert_{L^2}^{1/2}\Vert \ders^hk\Vert_{L^2}^{1/2}
+C\Vert \ders^h k \Vert_{L^2}\leq
C\left(\Vert
\ders^{h+1}k\Vert_{L^2}^{1/2}+ C\Vert \ders^hk\Vert_{L^2}^{1/2}\right)
\Vert \ders^h k \Vert_{L^2}^{1/2}\\
\Vert \pol_{h}(\vert\ders^{h-2}k\vert)\Vert_{L^\infty}
\leq &\,C\Vert\ders^{h+1}k\Vert_{L^2}^{\sigma}\Vert
k\Vert_{L^2}^{\sum_{l=0}^{h-2}\alpha_l-\sigma}+C\Vert k\Vert_{L^2}^h\leq
C\Vert\ders^{h+1}k\Vert_{L^2}^{\sigma}+ C\leq
C\left(\Vert
\ders^{h+1}k\Vert_{L^2}^{2}+ C\right)^{\sigma/2}\\
\Vert \ders^h k\Vert_{L^2}\leq &\,C\Vert
\ders^{h+1}k\Vert_{L^2}^{h/(h+1)}\Vert k\Vert_{L^2}^{1/(h+1)}
+C\Vert k \Vert_{L^2}\leq
C\Vert \ders^{h+1}k\Vert_{L^2} + C
\end{align*}
with $\sigma=\frac{\sum_{l=0}^{h-2}(l+1/2)\alpha_l}{h+1}$ as 
every monomial of $\pol_h(\vert\ders^{h-2}k\vert)$ is less than 
$C\prod_{l=0}^{h-2}\vert\ders^l k\vert^{\alpha_l}$.\\
Hence, putting together the first and third inequalities above,
$$
\Vert \ders^h k\Vert_{L^\infty}\leq C\left(\Vert
\ders^{h+1}k\Vert_{L^2}^{2}+ C\Vert \ders^hk\Vert_{L^2}^{2}\right)^{1/4}
\Vert \ders^h k \Vert_{L^2}^{1/2}\leq
C\left(\Vert \ders^{h+1}k\Vert_{L^2}^{2}+ C\right)^{1/4}
\Vert \ders^h k \Vert_{L^2}^{1/2}
$$
and multiplying this last with the second we get
\begin{align*}
t^{h-1}\Vert\ders^{h}&k\Vert_{L^\infty}\Vert\pol_{h}(\vert\ders^{h-2}k\vert)\Vert_{L^\infty}\\
&\,\leq Ct^{h-1}\left(\Vert \ders^{h+1}k\Vert_{L^2}^{2}+ C\right)^{\sigma/2+1/4}
\Vert \ders^h k \Vert_{L^2}^{1/2}\\
&\,= Ct^{h-1}\left(\int_{\TTT_t} \vert\ders^{h+1}k\vert^2\,ds
+ C\right)^{\frac{\sum_{l=0}^{h-2}(l+1/2)\alpha_l}{2(h+1)}+\frac{1}{4}}
\left(\int_{\TTT_t} \vert\ders^{h}k\vert^2\,ds\right)^{1/4}\\
&\, = C\frac{\left(t^h\int_{\TTT_t} \vert\ders^{h+1}k\vert^2\,ds
  +Ct^h\right)^{\frac{\sum_{l=0}^{h-2}(l+1/2)\alpha_l}{2(h+1)}+\frac{1}{4}}
\left(t^{h-1}\int_{\TTT_t} \vert\ders^{h}k\vert^2\,ds\right)^{1/4}}{t^{1-h}\cdot
t^{h\frac{\sum_{l=0}^{h-2}(l+1/2)\alpha_l}{2(h+1)}+\frac{h}{4}}
\cdot t^{\frac{h-1}{4}}}\\
&\, = C\frac{\left(t^h\int_{\TTT_t} \vert\ders^{h+1}k\vert^2\,ds
  + Ct^h\right)^{\frac{\sum_{l=0}^{h-2}(l+1/2)\alpha_l}{2(h+1)}+\frac{1}{4}}
\left(t^{h-1}\int_{\TTT_t} \vert\ders^{h}k\vert^2\,ds\right)^{1/4}}
{t^{3/4-h/2+h\frac{\sum_{l=0}^{h-2}(l+1/2)\alpha_l}{2(h+1)}}}\,.
\end{align*}
Now, applying Young inequality, if we elevate to
$\biggl(\frac{\sum_{l=0}^{h-2}(l+1/2)\alpha_l}{2(h+1)}+\frac{1}{4}\biggr)^{-1}$
the first term at the numerator, to $4$ the second, as
$\sum_{l=0}^{h-2}(l+1)\alpha_l=h$, it follows that the exponent
which, at the end, goes on $t$ at the denominator is
\begin{align*}
\theta_h=&\,\left(3/4-h/2+h\frac{\sum_{l=0}^{h-2}(l+1/2)\alpha_l}{2(h+1)}\right)
\left(1-1/4- \frac{\sum_{l=0}^{h-2}(l+1/2)\alpha_l}{2(h+1)}-1/4
\right)^{-1}\\
&\,= \left(3/4-h/2+h\frac{h-1/2\sum_{l=0}^{h-2}\alpha_l}{2(h+1)}\right)
\left(1/2- \frac{h-1/2\sum_{l=0}^{h-2}\alpha_l}{2(h+1)}\right)^{-1}\\
&\,= \frac{3h+3-2h^2-2h+2h^2-h\sum_{l=0}^{h-2}\alpha_l}{4(h+1)}\cdot
\frac{4(h+1)}{2h+2 -2h +\sum_{l=0}^{h-2}\alpha_l}\\
&\,= \frac{h+3-h\sum_{l=0}^{h-2}\alpha_l}{2 +\sum_{l=0}^{h-2}\alpha_l}
= \frac{3h+3}{2 +\sum_{l=0}^{h-2}\alpha_l}-h
\end{align*}
and since
$\sum_{l=0}^{h-2}\alpha_l\geq\sum_{l=0}^{h-2}\alpha_l\frac{l+1}{h-1}=\frac{h}{h-1}=1+1/(h-1)$,
we have
$$
\theta_h\leq \frac{2h-3}{3h-2}=1-\frac{h+1}{3h-2}<1\qquad\text{ { for
  every} even $h\geq2$.}
$$
Thus,
$$
\sum_{i=1}^3t^{h-1}\ders^{h-1} k^i\ders^{h}k^i\,
\biggr\vert_{\text{{ at
    the 3--point}}}
\leq \varepsilon_h/2 \left (t^h\int_{\TTT_t} \vert\ders^{h+1}k\vert^2\,ds
+ t^{h-1}\int_{\TTT_t} \vert\ders^{h}k\vert^2\,ds + Ct^h\right) +C/t^{\theta_h}
$$
with $\theta_h<1$ and $\varepsilon_h>0$ which can be chosen
arbitrarily small.

We apply this argument for every even $h$ from 2 to 
$j$, choosing accurately the values $\varepsilon_j>0$.\\
Hence, we can continue estimate~\eqref{ciaociao1} as follows,
\begin{align*}
\frac{d\,}{dt}\int_{{\TTT_t}}&\, k^2+ tk_s^2+ \frac{t^2 k_{ss}^2}{2!}+\dots+
\frac{t^j\vert \partial_s^jk\vert^2}{j!}\,ds\\
\leq&\,-\varepsilon/2\int_{{\TTT_t}} k_s^2 + t k_{ss}^2 + t^2 k_{sss}^2+\dots+
t^j\vert \partial_s^{j+1}k\vert^2\,ds + C + C/t^{\theta_2}+\dots+C/t^{\theta_j}\\
&\, + \dert\sum_{i=1}^3 t^2\qol_5(\lambda^i, k_s^i)
+ t^4\qol_9(\dert\lambda^i,k_{sss}^i)
+\dots+ t^j\qol_{2j+1}(\dert^{j/2-1}\lambda^i,\ders^{j-1}k^i)
\,\biggr\vert_{\text{{ at the 3--point}}}\\
\leq &\, C + C/t^\theta + \dert\sum_{i=1}^3 t^2\qol_5(\lambda^i, k_s^i)
+ t^4\qol_9(\dert\lambda^i,k_{sss}^i)
+\dots+ t^j\qol_{2j+1}(\dert^{j/2-1}\lambda^i,\ders^{j-1}k^i)
\,\biggr\vert_{\text{{ at the 3--point}}}
\end{align*}
for some $\theta<1$.\\
Integrating this inequality in time on $[0,t]$ with $t\leq T_M$ and taking
into account Remark~\ref{qolpol}, we get
\begin{align*}
\int_{{\TTT_t}} k^2&\,+ tk_s^2+ \frac{t^2 k_{ss}^2}{2!}+\dots+
\frac{t^j\vert \partial_s^jk\vert^2}{j!}\,ds\\
\leq&\, \int_{{\TTT_0}} k^2\,ds + CT_M+ CT_M^{(1-\theta)}\\
&\,+ \sum_{i=1}^3 t^2\qol_5(\lambda^i, k_s^i)
+ t^4\qol_9(\dert\lambda^i,k_{sss}^i)
+\dots+ t^j\qol_{2j+1}(\dert^{j/2-1}\lambda^i,\ders^{j-1}k^i)
\,\biggr\vert_{\text{{ at the 3--point}}}\\
\leq&\, \int_{{\TTT_0}} k^2\,ds + C
+ t^2\Vert \pol_5(\vert k_s\vert)\Vert_{L^\infty}
+ t^4\Vert \pol_9(\vert k_{sss}\vert)\Vert_{L^\infty}
+\dots+ t^j\Vert \pol_{2j+1}(\vert \ders^{j-1}k\vert) \Vert_{L^\infty}\,.
\end{align*}
Now we absorb all the polynomial terms, after interpolating each
one of them between the corresponding ``good'' integral in the left
member and some power of the $L^2$ norm of $k$, 
as we did in showing Proposition~\ref{pluto1000}, hence we finally
obtain for every even $j\in\NN$,
$$
\int_{{\TTT_t}} k^2+ tk_s^2+ \frac{t^2 k_{ss}^2}{2!}+\dots+
\frac{t^j\vert \partial_s^jk\vert^2}{j!}\,ds \leq C_j
$$
with $t\in[0,T_M]$ and a constant $C_j$ depending
only on $\int_{\TTT_0}k^2\,ds$ and the inverses of the lengths of the
three curves at time zero.

This family of inequalities clearly implies
$$
\int_{{\TTT_t}} \vert \partial_s^jk\vert^2 \,ds\leq
\frac{C_j j!}{t^j}\qquad\text{ { for every even} $j\in\NN$.}
$$

Then, passing as before from integral to $L^\infty$ estimates by means of
inequalities~\eqref{int2}, we have the following proposition.

\begin{prop}\label{topolino5}
For every $\mu>0$ the curvature and all
  its space derivatives of $\TTT_t$ are uniformly bounded in the time
  interval $[\mu,T_M]$ (where $T_M$ is given by
  Proposition~\ref{stimaL}) by some constants depending only on $\mu$,
  the $L^2$ norm of $k$ of $\TTT_0$ and the inverses of the lengths of the three
  curves at time zero.
\end{prop}

By means of these a priori estimates we can now work out some
results about the flow and improve Theorem~\ref{smoothexist}.

\begin{teo}\label{curvexplod}
If $[0,T)$ is the maximal time interval of existence of a smooth
solution $\TTT_t$ with $T<+\infty$ of Problem~\eqref{problema},
then
\begin{enumerate}
\item either the inferior limit of the length of at least one curve of
  $\TTT_t$ goes to zero when $t\to T$,
\item or $\limup_{t\to T}\int_{\TTT_t}k^2\,ds=+\infty$.
\end{enumerate}
Moreover, if the lengths of the three curves are uniformly positively
bounded from below, then the superior limit in $(2)$ is a limit and 
there exists a positive constant $C$ such that $\int_{{\TTT_t}} k^2\,ds \geq
C/\sqrt{T-t}\to+\infty$ for every $t\in[0, T)$.
\end{teo}
\begin{proof}
If the three lengths are uniformly bounded away from zero and the
$L^2$ norm of $k$ is bounded, by Proposition~\ref{unif222} and
Ascoli--Arzel\`a Theorem, the triods
$\TTT_t$ converge in $C^\infty$ to a smooth triod $\TTT_T$ as $t\to
T$. Then, applying Theorem~\ref{smoothexist} to $\TTT_T$
we could restart the flow obtaining a smooth evolution on a longer
time interval, hence  contradicting the maximality of the interval
$[0,T)$.\\
By means of differential inequality~\eqref{evolint999}, we have
$$
\frac{d\,}{dt} \int_{{\TTT_t}} k^2\,ds
\leq C \left(\int_{{\TTT_t}} k^2\,ds\right)^{3} + C
\leq C \left(1+\int_{{\TTT_t}} k^2\,ds\right)^{3}\,,
$$
which, after integration between $t,r\in[0,T)$ with $t<r$, gives
$$
\frac{1}{\left(1+\int_{{\TTT_t}} k^2\,ds\right)^{2}}
-\frac{1}{\left(1+\int_{{\TTT_r}} k^2\,ds\right)^{2}}\leq C(r-t)\,.
$$
Then, if case $(1)$ does not hold, we can choose 
a sequence of times $r_j\to T$ such that 
$\int_{\TTT_{r_j}} k^2\,ds\to+\infty$. Putting $r=r_j$ in the inequality above and passing
to the limit we get
$$
\frac{1}{\left(1+\int_{{\TTT_t}} k^2\,ds\right)^{2}}\leq C(T-t)
$$
hence,
$$
\int_{{\TTT_t}} k^2\,ds \geq \frac{{C}}{\sqrt{T-t}}-1\geq
\frac{C}{\sqrt{T-t}}\to+\infty\,,
$$
for some positive constant $C$.
\end{proof}
This theorem obviously implies the following corollary.
\begin{cor}\label{kexplod} If $[0,T)$ is the maximal time
interval of existence of a smooth solution $\TTT_t$ with $T<+\infty$ 
and the three lengths are uniformly bounded away from zero, then
\begin{equation}\label{krate}
\max_{\TTT_t}k^2\geq\frac{C}{\sqrt{T-t}}\to+\infty\,,
\end{equation}
as $t\to T$.
\end{cor}

\begin{rem}\label{kratrem}
In the case of the evolution of a closed curve in the plane 
there exist a constant $C>0$ such that if at time $T>0$ a 
singularity develops, then 
$$
\max_{{\TTT}_t} k^2\geq\frac{C}{{T-t}}
$$
for every $t\in[0,T)$ (see~\cite{huisk3}).

If this lower bound on the rate of blowing up of the curvature (which
is clearly stronger than the one in inequality~\eqref{krate}) holds
also in the case of the evolution of a triod is an open problem.
\end{rem}

\begin{prop}\label{unif333} For every $M>0$ there exists a
  positive time $T_M$ such that if the $L^2$ norm of the curvature and
  the inverses of the lengths of the smooth triod $\TTT_0$ are bounded
  by $M$, then the  maximal time of existence $T>0$ of the associated solution of
Problem~\eqref{problema} with initial data $\TTT_0$ is larger than
$T_M$.
\end{prop}
\begin{proof}
By Proposition~\ref{stimaL} in the interval $[0,\min\{T_M,T\})$ the
$L^2$ norm of $k$ and the inverses of the lengths of the three curves
of $\TTT_t$ are bounded by $2M^2+6M$.\\
Then, by Theorem~\ref{curvexplod}, the value $\min\{T_M,T\}$ cannot
coincide with the maximal time of existence, hence $T>T_M$.
\end{proof}

By means of Proposition~\ref{topolino5} we can now work out an
existence result for an initial triod $\TTT_0$ which 
is neither smooth nor geometrically smooth, but it is only
$C^2$ and satisfies the $120$ degrees condition.

\begin{teo}\label{brakkeevolution} If $\TTT_0$ is a $C^2$ initial
  triod (not necessarily geometrically smooth) 
  then there exists a Brakke flow with equality
  $\TTT_t$ of the initial triod $\TTT_0$ for some positive time interval
  $[0,T)$.\\
  Moreover, the triods $\TTT_t$ are geometrically smooth for every time $t>0$ and
  the curvatures $k^i$ belong to $C^\infty([0,1]\times(0,T))$, hence
  the flow is a smooth Brakke flow with equality for every positive time.\\
  Finally, the unit tangents 
  $\tau^i$ are continuous in $[0,1]\times[0,T)$ and the function 
  $\int_{\TTT_t}k^2\,ds$ is continuous on $[0,T)$.
\end{teo}
\begin{proof}
We can approximate in $W^{2,2}(0,1)$ (hence in $C^1([0,1])$)
the triod $\TTT_0=\cup_{i=1}^3\sigma^i$ with a family of smooth triods
$\TTT_j$, composed of $C^\infty$ curves $\sigma_j^i\to\sigma^i$, as
$j\to\infty$ with
$\frac{d\sigma_j^i(0)}{dx}=\frac{d\sigma^i(0)}{dx}$ and
$\frac{d^2\sigma_j^i(0)}{dx^2}=0$.\\
By the convergence in $W^{2,2}$ and in $C^1$, the inverses of the
lengths of the initial curves, the integrals $\int_{\TTT_j} k^2
+\lambda^2\,ds$ and $\vert\partial_x\sigma_j^i(x)\vert$ (from
above and away from zero) for all the approximating triods are
equibounded, thus Proposition~\ref{unif333} assures the existence
of a uniform interval $[0,T)$ of existence of smooth
evolutions given by the curves
$\gamma_j^i(x,t):[0,1]\times[0,T)\to\overline{\Omega}$.\\
Now, by the same reason, Proposition~\ref{topolino5} gives uniform
estimates on the $L^\infty$ norms of the curvature and of all its
derivatives in every rectangle $[0,1]\times[\mu,T)$, with
$\mu>0$.\\
We can then select a subsequence (not relabelled) such that
the curves $\gamma_j^i$, after 
reparametrization proportional to arclength,
converge to some
$\gamma^i(x,t):[0,1]\times[0,T)\to\overline{\Omega}$ 
(composing the triods $\TTT_t$),
\begin{itemize}
\item uniformly in $[0,1]\times[0,T)$,
\item in $C^\infty$ in every rectangle $[0,1]\times[\mu,T)$,
  with $\mu>0$.
\end{itemize}
Moreover,  since all the approximating flows are composed
of smooth triods and the curvatures converge smoothly, when $t>0$ the triods
$\TTT_t$ are geometrically smooth (see Remark~\ref{geocomprem}).\\
It is then an exercise to see that the unit tangents $\tau^i$ are
continuous functions also at $t=0$, that is, on all the rectangle
$[0,1]\times[0,T)$ (by the uniform control on $\Vert k\Vert_{L^2}$ and
Sobolev embedding theorem).
Notice that, the continuity of $\gamma^i$ also implies that the measures
$\HH^1\res\TTT_t$ weakly$^{{\displaystyle{\star}}}$ converge to $\HH^1\res\TTT_0$, where
$\HH^1$ is the Hausdorff one--dimensional measure.\\
Now we show that $\TTT_t=\cup_{i=1}^3\gamma^i$ is a 
Brakke flow with equality. By the smoothness of the flow 
for every positive time, we have only to check the derivative
$\frac{d\,}{dt}\int_{\TTT_t}\varphi\,ds$ at $t=0$.\\
For every smooth positive test function
$\varphi:\overline{\Omega}\times[0,T)\to\R$ 
the functional $\int_{\TTT}\varphi k^2\,ds$ is lower semicontinuous in
the convergence of the triods with their unit tangents (see~\cite{simon}, 
moreover,  integrating on $[0,t)$
inequality~\eqref{evolint999} (forgetting the absolute value), 
for the approximating flows $\gamma_j^i$, and passing to the limit, we see that
the function $\int_{\TTT_t}k^2\,ds$ is actually {\em continuous} at
$t=0$. Then, by a standard argument, it follows that also the functions
$\int_{\TTT_t}\varphi\,k^2\,ds$ are continuous at $t=0$, 
for every positive $\varphi$, hence {\em for
  every} smooth function $\varphi$ with compact support in $\Omega$.\\
Analogously, also the terms
$\int_{\TTT_t}\langle\nabla\varphi\,\vert\,\underline{k}\rangle\,ds$
and $\int_{\TTT_t}\varphi_t\,ds$ are continuous at $t=0$, hence
integrating equation~\eqref{brakkeqqq}, satisfied by the approximating
flows, on $[0,t)$ and then passing to the limit we get
$$
\int_{\TTT_t}\varphi\,ds-\int_{\TTT_0}\varphi\,ds=
-\int_{0}^t\int_{\TTT_\xi}\varphi\,k^2\,ds\,d\xi
+\int_{0}^t\int_{\TTT_\xi}\langle\nabla\varphi\,\vert\,\underline{k}\rangle\,ds\,d\xi
+ \int_{0}^t\int_{\TTT_\xi}\varphi_t\,ds\,d\xi
$$
which clearly says, by the fundamental theorem of calculus, that the
derivative $\frac{d\,}{dt}\int_{\TTT_t}\varphi\,ds$ exists at $t=0$ and
that $\TTT_t$ is a Brakke flow with equality.
\end{proof}

\begin{rem}\label{oprob1} \hspace{1truecm}
  \begin{enumerate}
  \item The relevance of this theorem is that the initial
    triod is not required to satisfy any compatibility
    condition, but only to have angles of $120$ degrees,
    in particular, it is not necessary that the sum of the three
    curvatures at its 3--point is zero.
  \item It should be noticed that if the three initial curves
    are $C^\infty$, the flow $\gamma^i$ is smooth till $t=0$ far from
    the 3--point, that is, in a closed rectangle included in
    $[0,1]\times[0,T)\setminus\{(0,0)\}$ we can locally reparametrize
    the curves to get a smooth flow also at $t=0$.\\
    This follows from the local estimates for the motion by curvature
    (see~\cite{eckhui2}, for instance).
  \item As we said in the introduction, the next important question is what
    can be said if the initial triod does not satisfy the $120$ degrees
    condition. One would hope to have a
    suitable definition of evolution (possibly {\em weak}) such that the $120$
    degrees condition is satisfied instantaneously, that is, at every
    positive time, like it happens here for the geometrical smoothness.
  \item The uniqueness of the limit $\gamma^i$ is an open problem as
    well as its dependence on the approximating procedure.\\ 
    Even more important is the geometric uniqueness of such a Brakke 
    flow with continuous unit tangents, forgetting the
    parametrizations and looking at the triods as subsets of $\R^2$.\\
    Finally, if the initial triod is smooth (or
    geometrically smooth) this flow should be a reparametrization
    of the smooth evolution given by Theorem~\ref{smoothexist} (see
    Proposition~\ref{geouniq}).
  \item We do not know if every Brakke flow with
    equality starting from a $C^2$ initial triod $\TTT_0$, which
    becomes immediately smooth (possibly 
    requiring also the continuity of the unit tangents), 
    can be obtained in this way, even when the initial triod 
    $\TTT_0$ is smooth.\\
    This problem is clearly related to the uniqueness of the 
    smooth Brakke flows with equality (maybe further restricting the
    candidates to a special class with extra geometric properties). A 
    positive answer would also allow us to extend to them 
    the analysis of the singularities carried on in the next sections.
  \end{enumerate}
\end{rem}

\begin{rem}\label{manytrip}
We point out that all this section can be extended to
  networks of curves with many 3--points.\\
  Indeed, an analog of the result of Bronsard and Reitich for such situation 
  (they also remark that) can be obtained generalizing the
  ``algebraic'' analysis in their paper~\cite{bronsard} in order to
  show the {\em complementary conditions} for the 
  system~\ref{prob0} associated to the network. 
  Then, all the estimates can be generalized simply adding the
  contributions of all the 3--points and of every end point, since each one of
  them has to satisfy the relations between $k$, $\lambda$ and 
  their derivatives computed in Section~\ref{due} for a single triod.
\end{rem}

In all the discussion of this section, we did not take
care of the fact that the triods have to remain in the domain
$\Omega$ and of the condition of embeddedness,
which are required in the formulation of Problem~\eqref{problema},
actually, we only concentrated on the analytic properties of the
solution of the parabolic system~\eqref{prob0}.\\
It will follow by the geometric results of the
next section that since the initial triod is embedded, if the lengths
of three curves stay away from zero, then, during the evolution,
the triods do not develop self--intersections and ``touch'' the
boundary of $\Omega$ only with their end points.

{\em In the rest of the paper we restrict ourselves only to the smooth
  flows given by Theorem~\ref{smoothexist} and we will 
analyse the possible formation of singularities}.

\section{Geometric Properties of the Flow}
\label{disuellesec}

Let us consider a smooth evolution $\TTT_t$ in the time 
interval $[0,T)$ of an initially embedded smooth triod $\TTT_0$ in
the convex $\Omega$.

The first thing we want to show is that the triods cannot get out of
$\Omega$.

\begin{prop}\label{omegaok} The triods $\TTT_t$  
intersect the boundary of $\Omega$ only at the end points.\\
Moreover, for every positive time such intersections are transversal.
\end{prop}
\begin{proof}
Even if some of the three curves of the initial triod are tangent to
$\partial\Omega$ at the end points $P^i$, by the strong maximum
principle, as $\Omega$ is convex, the intersections become immediately
transversal and stay so for every time.\\
By continuity, the 3--point cannot hit the boundary at least 
for some time $T^\prime>0$.  Then, fixing a time $t\in[0,T^\prime)$,
even if the curve $\gamma^i(\cdot,t)$ intersects 
the boundary at some of its inner points, that is,
$\gamma^i(x,t)\in\partial\Omega$ with $x\not=\{0,1\}$, again by the
strong maximum principle and the convexity of $\Omega$ we can conclude
that for every subsequent time in the interval $(t,T^\prime)$ there are no
more intersections.\\
This argument clearly implies that if $t_0>0$ is the ``first time'' when
the triods intersect the boundary with an inner point, this latter has
to be the 3--point $O$. The minimality of $t_0$ is then easily
contradicted by the convexity of $\Omega$, the $120$ degrees condition
and the non zero length of the three curves of $\TTT_{t_0}$.
\end{proof}

Now we concentrate on the condition of embeddedness.

Given the smooth flow $\TTT_t=F(\TTT,t)$, we take two points $p=F(x,t)$ and
$q=F(y,t)$ belonging to $\TTT_t$ 
and we define $\Gamma_{p,q}$ to be the {\em geodesic}
curve contained in $\TTT_t$ connecting $p$ and $q$. Then we let 
$A_{p,q}$ to be the area of the open region ${\mathcal A}_{p,q}$ in $\R^2$ enclosed
by the segment  $[p,q]$ and the curve $\Gamma_{p,q}$.
When the ${\mathcal A}_{p,q}$ is not connected, we let
$A_{p,q}$ to be the sum of the areas of its connected components.

We consider the function 
$\Phi_t: \TTT\times\TTT\to \R\cup\{+\infty\}$ as 
\begin{equation*}
\Phi_t(x,y)=
\begin{cases}
\frac{\vert p-q\vert^2}{A_{p,q}}\qquad &\text{ { if} $x\not=y$,}\\
4\sqrt{3}\qquad  & \text{ { if} $x$ and $y$ coincide with the 3--point
  $O$ of $\TTT$,}\\
+\infty\qquad  & \text{ { if}  $x=y\not=O$}
\end{cases}
\end{equation*}
where $p=F(x,t)$ and $q=F(y,t)$.\\
Since $\TTT_t$ is smooth and the $120$ degrees condition holds, it is
easy to check that $\Phi_t$ is a lower semicontinuous function. Hence,
by the compactness of $\TTT$, the following infimum is actually a
minimum
\begin{equation}\label{equant}
E(t) = \inf_{x,y\in\TTT}\Phi_t(x,y)
\end{equation}
for every $t\in[0,T)$.\\
Similar geometric quantities have already been applied to similar problems 
in~\cite{hamilton3}, ~\cite{chzh} and~\cite{huisk2}.

If the triod $\TTT_t$ has no self--intersections we have
$E(t)>0$, the converse is clearly also true.\\
Moreover, $E(t)\leq\Phi_t(0,0)=4\sqrt{3}$ always holds, thus when
$E(t)>0$ the two points $(p,q)$ of a minimizing pair $(x,y)$ can coincide if
and only if $p=q=O$.\\
Finally, since the evolution is smooth it is easy to see that the
function $E:[0,T)\to\R$ is continuous.

These properties set the question of the possible self--intersections of
solutions of Problem~\eqref{problema} that we let open in the previous
section at the end of the proof of Theorem~\ref{smoothexist} 
(the fact that the triods remain in $\Omega$ is shown by
Proposition~\ref{omegaok}). Indeed, if the initial triod $\TTT_0$ is
embedded, we have $E(0)>0$, hence $E(t)>0$ for some time.

Now we want to show something more, that is, in {\em all} the
maximal interval of existence of a smooth flow self--intersections cannot
happen, in other words, if $E(0)>0$ then $E(t)>0$ for every
$t\in[0,T)$.\\
Since we are dealing with embedded triods, with a little abuse of
notation, we consider the function $\Phi_t$ defined on
$\TTT_t\times\TTT_t$ and we speak of a minimizing pair for the couple
of points $(p,q)\in\TTT_t\times\TTT_t$ instead of
$(x,y)\in\TTT\times\TTT$.

\begin{lemma}\label{lemet1}
Assume that $0<E(t)<4\sqrt{3}$, then for any minimizing pair 
$(p,q)$ we have $p\ne q$ and
neither $p$ nor $q$ coincides with the 3--point $O$ of $\TTT_t$.\\
Moreover, it is always possible to find a minimizing pair such that 
${\mathcal A}_{p,q}$ is a single connected region and 
the segment $[p,q]$ meet the (one or two) curves containing $p$ 
and $q$ only at these points and with transversal intersection.
\end{lemma}
\begin{proof}
We already saw that by the very definition of $\Phi_t$, the inequality
$0<E(t)<4\sqrt{3}$  implies $p\ne q$. 
The proof that $p$ and $q$ can be chosen in order that the region
${\mathcal A}_{p,q}$ is connected and the segment $[p,q]$ does not
intersect the curves which contain the two points, 
goes like in~\cite[Lemma~2.1]{chzh}).\\
We prove the first claim, assuming by contradiction that 
$p=O\in\TTT_t$ and $\Phi_t(O,q)=E(t)$.
By the above, we can suppose that the segment $[O,q]$ is contained in
in the sector between the curves $\gamma^1$ and
$\gamma^2$, that $q\in\gamma^1$ and that the region ${\mathcal
  A}_{O,q}$ is bounded by such segment and the curve $\gamma^1$.\\
If the angle $\alpha\geq0$ formed by $[O,q]$ and $\tau^2(O,t)$ is
smaller than $90$ degrees it is easy to see that moving a little the
point $p$ along $\gamma^2$, the distance $\vert p-q\vert$ decreases
while the area $A_{p,q}$ increases, hence the ratio $\vert
O-q\vert^2/A_{O,q}$ cannot be minimal.\\
Thus the width of the angle $\alpha$ has to be greater or equal than 
$90$ degrees.\\
We consider then the points $p(s)=\gamma^1(x(s),t)$ with arclength
parameter $s\in[0,\varepsilon)$ (then $p(0)=O$ and
$\frac{dp(0)}{ds}=\tau^1(O,t)$)  and we compute 
the right derivative at $s=0$ of $\Phi_t(p(s),q)$ 
(see the proof of Proposition~\ref{lemet2}),
\begin{align*}
\frac{d\,}{ds}\Phi_t(p(s),q)\,\biggl\vert_{s=0}=&\,
\frac{A_{p(s),q}\frac{d\,}{ds}\vert p(s)-q\vert^2- \vert
  p(s)-q\vert^2\frac{d\,}{ds}A_{p(s),q}}
{A_{p(s),q}^2}\,\biggl\vert_{s=0}\\
=&\,\frac{-2A_{p(0),q}\vert p(0)-q\vert\cos{(120-\alpha)}+1/2\vert
  p(0)-q\vert^3 \sin{(120-\alpha)}}{A_{p(0),q}^2}\\
=&\,\frac{\vert O-q\vert}{2A_{O,q}^2}\biggl[-4A_{O,q}\cos{(120-\alpha)}
+\vert O-q\vert^2 \sin{(120-\alpha)}\biggr]
\end{align*}
which has to be non negative by minimality.\\
Hence, it follows that 
$ 4A_{O,q}\cos{(120-\alpha)}\leq \vert O-q\vert^2 \sin{(120-\alpha)}$ and 
$$
4\cot{(120-\alpha)}\leq \frac{\vert O-q\vert^2}{A_{O,q}}=\Phi_t(O,q)=E(t)\,.
$$
Since $\pi/2\leq\alpha\leq2\pi/3$ we have
$\cot{(120-\alpha)}\geq\sqrt{3}$, hence 
we conclude $E(t)\geq 4\sqrt{3}$ which is in contradiction with the
initial hypothesis.\\
Finally, as neither $p$ nor $q$ 
coincides with the 3--point, if $[p,q]$ intersects $\TTT_t$
tangentially at $p$ or $q$, moving such point a little in the
direction which decreases the distance $\vert p-q\vert$, the area $A_{p,q}$
is little changed, so a variation as above gives a contradiction 
(see~\cite[Lemma~2.1]{chzh}).
\end{proof}

\begin{rem}\label{goodpoints} Looking at the proof of the connectedness 
  in~\cite{chzh}, it can be proved also that if
  there is at least one minimizing pair such that both its 
  points are distinct from the end points $P^i$, then the same holds also
for the pair $(p,q)$ with ${\mathcal A}_{p,q}$ connected whose
existence is assured by this lemma.
\end{rem}

\begin{prop}\label{lemet2} The function $E(t)$ is monotone increasing 
  in every time interval where $0<E(t)<4\sqrt{3}$ and 
  for at least one minimizing pair $(p,q)$ of $\Phi_t$ 
  neither $p$ nor $q$ coincides with one of the end points $P^i$.
\end{prop}
\begin{proof} 
We assume that $0<E(t)<4\sqrt{3}$ and that there exists a minimizing
pair $(p,q)$ for $\Phi_t$ such that the two points are both distinct
from the end points $P^i$, for every $t$ in some interval of time. 
Since $E(t)$ is a locally Lipschitz function, 
to prove the statement it is then enough to show that
$\frac{dE(t)}{dt}>0$ for every time $t$ such that this derivative
exists (which happens almost everywhere in the interval).\\
Fixed a minimizing pair $(p,q)$ at time $t$, satisfying the
conclusions of Lemma~\ref{lemet1} and Remark~\ref{goodpoints}, 
we choose a value  $\eps>0$ smaller than the geodesic distances of $p$ and $q$ from
the 3--point $O$ of $\TTT_t$ and between them, moreover if $p$ and $q$
both belong to the same curve we can also suppose that $q$ is the
closest to $O$.\\
By simplicity, we discuss the situation where the points $p$, $q$ are like in 
Figure~\ref{figura1}, the computations in the other cases are analogous.\\
\begin{figure}[h]
\begin{picture}(0,0)%
\includegraphics{trip1.pstex}%
\end{picture}%
\setlength{\unitlength}{3947sp}%
\begingroup\makeatletter\ifx\SetFigFont\undefined%
\gdef\SetFigFont#1#2#3#4#5{%
  \reset@font\fontsize{#1}{#2pt}%
  \fontfamily{#3}\fontseries{#4}\fontshape{#5}%
  \selectfont}%
\fi\endgroup%
\begin{picture}(4449,3932)(64,-3185)
\put(2251,614){\makebox(0,0)[lb]{\smash{\SetFigFont{12}{14.4}{\familydefault}{\mddefault}{\updefault}{$P^3$}%
}}}
\put(151,-3136){\makebox(0,0)[lb]{\smash{\SetFigFont{12}{14.4}{\familydefault}{\mddefault}{\updefault}{$P^1$}%
}}}
\put(4051,-3136){\makebox(0,0)[lb]{\smash{\SetFigFont{12}{14.4}{\familydefault}{\mddefault}{\updefault}{$P^2$}%
}}}
\put(2326,-1486){\makebox(0,0)[lb]{\smash{\SetFigFont{12}{14.4}{\familydefault}{\mddefault}{\updefault}{$O$}%
}}}
\put(3151,-1786){\makebox(0,0)[lb]{\smash{\SetFigFont{12}{14.4}{\familydefault}{\mddefault}{\updefault}{$q$}%
}}}
\put(3451,-1036){\makebox(0,0)[lb]{\smash{\SetFigFont{12}{14.4}{\familydefault}{\mddefault}{\updefault}{${\mathcal A}_{p,q}$}%
}}}
\put(548,-2618){\makebox(0,0)[lb]{\smash{\SetFigFont{12}{14.4}{\familydefault}{\mddefault}{\updefault}{$\tau(p)$}%
}}}
\put(1196,-2567){\makebox(0,0)[lb]{\smash{\SetFigFont{12}{14.4}{\familydefault}{\mddefault}{\updefault}{$\alpha(p)$}%
}}}
\put(1095,-2348){\makebox(0,0)[lb]{\smash{\SetFigFont{12}{14.4}{\familydefault}{\mddefault}{\updefault}{$p$}%
}}}
\put(3526,-361){\makebox(0,0)[lb]{\smash{\SetFigFont{29}{34.8}{\familydefault}{\mddefault}{\updefault}{$\Omega$}%
}}}
\end{picture}
 \caption{\label{figura1}}
\end{figure}
Possibly taking a smaller $\eps>0$, we fix an arclength
coordinate $s\in (-\eps,\eps)$ and a local parametrization $p(s)$ of
the curve containing in a neighborhood of $p=p(0)$, with the 
same orientation of the original one.
Let $\eta(s)=\vert p(s)-q\vert$ and $A(s)=A_{p(s),q}$, since  
\begin{equation*}
E(t)=\min_{s\in(-\eps,\eps)}\frac{\eta^2(s)}{A(s)}=
\frac{\eta^2(0)}{A(0)}\,,
\end{equation*}
if we differentiate in $s$ we obtain 
\begin{equation}\label{eqdsul}
\frac{d\eta^2(0)}{ds}A(0)=
\frac{dA(0)}{ds}\eta^2(0)\,.
\end{equation}
As the intersection of the segment $[p,q]$ with the triod is
transversal, we have an angle $\alpha(p)\in(0,\pi)$ 
determined by the unit tangent $\tau(p)$ and the vector
$q-p$. We compute
\begin{align*} 
\frac{{ d} \eta^2(0)}{{ d} s} &=\,
-2 \langle \tau(p)\,\vert\, q-p\rangle = -2 \vert
p-q\vert \cos\alpha(p)\\ 
\frac{{ d} A(0)}{{ d} s} &=\, 
\frac 12 \vert \tau(p) \wedge (q-p)\vert = 
\frac 12 \langle \nu(p)\,\vert\, q-p\rangle = 
\frac 12  \vert p-q\vert \sin\alpha(p)\\
\end{align*}
Putting these derivatives in equation~\eqref{eqdsul} and
recalling that $\eta^2(0)/A(0)=E(t)$,  we get
\begin{equation}\label{topolino2}
\cot\alpha(p) = -\frac{\vert p-q\vert^2}{4A_{p,q}}= -\frac{E(t)}{4}\,.
\end{equation}
Since $0<E(t)<4\sqrt{3}$ we get $-\sqrt{3}<\cot\alpha(p)<0$ which
implies
\begin{equation*}
\frac{\pi}{2} <\alpha(p) < \frac 56 \pi\,.
\end{equation*}
The same argument clearly holds for the point $q$, hence 
defining $\alpha(q)\in(0,\pi)$ to be the angle determined by the
unit tangent $\tau(q)$ and the vector $p-q$, by equation~\eqref{topolino2}
it follows that $\alpha(p)=\alpha(q)$ and we simply write $\alpha$
for both.

We consider now a different variation, moving at the same time the
points $p$ and $q$, in a way that $\frac{dp(s)}{ds}=\tau(p(s))$ and 
$\frac{dq(s)}{ds}=\tau(q(s))$.\\
As above, letting $\eta(s)=\vert p(s)-q(s)\vert$ and 
$A(s)=A_{p(s),q(s)}$, by minimality we have
\begin{equation}\label{eqdersec}
\frac{{ d} \eta^2(0)}{{ d} s}A(0)=
\frac{{ d} A(0)}{{ d} s}\eta^2(0) \qquad\text{ { and} } \qquad
\frac{{ d}^2 \eta^2(0)}{{ d} s^2}A(0)\ge 
\frac{{ d}^2 A(0)}{{ d} s^2}\eta^2(0)\,.
\end{equation}
Computing as before,
\begin{align*}
\frac{{ d} \eta^2(0)}{{ d} s} 
&=\,2 \langle p-q \,\vert\,\tau(p)-\tau(q)\rangle\\
\frac{{ d} A(0)}{{ d} s} = -4\vert p-q\vert\cos\alpha
&=\, -\frac12\langle p-q\,\vert\,\nu(p)+\nu(q)\rangle\\
\frac{{ d}^2 \eta^2(0)}{{ d} s^2} = \vert p-q\vert\sin\alpha
&=\,2 \langle \tau(p)-\tau(q) \,\vert\,\tau(p)-\tau(q)\rangle +
2 \langle p-q \,\vert\, k(p)\nu(p) -k(q)\nu(q)\rangle\\
&=\, 2\vert \tau(p)-\tau(q)\vert^2+
2 \langle p-q \,\vert\, k(p)\nu(p) -k(q)\nu(q)\rangle\\
&=\, 8\cos^2\alpha+
2 \langle p-q \,\vert\, k(p)\nu(p) -k(q)\nu(q)\rangle\\
\frac{{ d}^2 A(0)}{{ d} s^2} 
&=\,\frac12\langle \tau(p)-\tau(q)\,\vert\,\nu(p)+\nu(q)\rangle
+\frac12\langle p-q\,\vert\,k(p)\tau(p)+k(q)\tau(q)\rangle\\
&=\,\frac12\langle\tau(p)\,\vert\,\nu(q)\rangle
-\frac12\langle \tau(q)\,\vert\,\nu(p)\rangle
+\frac12\langle p-q\,\vert\,k(p)\tau(p)+k(q)\tau(q)\rangle\\
&=\,-2\sin\alpha\cos\alpha
-1/2\vert p-q\vert(k(p)-k(q))\cos\alpha\,.
\end{align*}
Putting the last two relations in the second inequality 
of~\eqref{eqdersec}, we get 
\begin{align*}
(8\cos^2\alpha+&\, 
2 \langle p-q \,\vert\, k(p)\nu(p) -k(q)\nu(q)\rangle)A_{p,q}\\
\geq&\,(-2\sin\alpha\cos\alpha-1/2\vert p-q\vert(k(p)-k(q))\cos\alpha
)\vert p-q\vert^2
\end{align*}
hence, keeping in mind that $\tan\alpha=-4/E(t)$, by
equation~\eqref{topolino2}, we obtain
\begin{align}\label{eqfin}
2A_{p,q}\langle p-q \,\vert\, &\,k(p)\nu(p) -k(q)\nu(q)\rangle
+1/2\vert p-q\vert^3(k(p)-k(q))\cos\alpha\\
\geq &\,-2\sin\alpha\cos\alpha\vert p-q\vert^2
-8A_{p,q}\cos^2\alpha\nonumber\\
=&\,-2A_{p,q}\cos^2\alpha\left(\tan\alpha\frac{\vert
  p-q\vert^2}{A_{p,q}} + 4\right)\nonumber\\
=&\,-2A_{p,q}\cos^2\alpha\left(-\frac{4}{E(t)}\,E(t) + 4\right)=0\,.\nonumber
\end{align}

We consider now a time $t_0$ such that the derivative
$\frac{dE(t_0)}{dt}$ exists and we compute it with the following
standard trick,
$$
\frac{dE(t_0)}{dt} = \left.\frac{\partial}{\partial
    t}\Phi_t(p,q)\right\vert_{t=t_0}
$$
for {\em any} pair $(p,q)$ such that $p, q\in\TTT_{t_0}$ and $\frac{\vert
  p-q\vert^2}{A_{p,q}}=E(t_0)$.\\
Considering then a minimizing pair $(p,q)$ for $\Phi_{t_0}$ with
all the previous properties, by minimality, we are free to choose the 
``motion'' of the points $p(t)$, $q(t)$ ``inside'' the triods $\TTT_t$
in computing such partial derivative.\\
Since locally the triods are moving by curvature and we know that
neither  $p$ nor $q$ coincides with the 3--point or the end points, 
we can find $\eps>0$ and two  smooth curves $p(t), q(t)\in\TTT_t$  for
every $t\in(t_0-\eps,t_0+\eps)$ such that
\begin{align*}
p(t_0) &=\, p \qquad \text{ and }\qquad \frac{dp(t)}{dt} = k(p(t), t)
~\nu(p(t),t)\,, \\
q(t_0) &=\, q \qquad \text{ and }\qquad \frac{dq(t)}{dt}= k(q(t),t)
~\nu(q(t),t)\,. 
\end{align*}
Then,
\begin{equation}\label{eqderE}
\frac{dE(t_0)}{dt} = 
\left.\frac{\partial}{\partial t}\Phi_{t}(p,q)\right\vert_{t=t_0}=\frac{1}{A^2_{p,q}}\left.\left( A_{p,q}
\frac{d\vert p(t)-q(t)\vert^2}{dt} - 
\vert p-q\vert^2 \frac{dA_{p(t),q(t)}}{dt}\right)\right\vert_{t=t_0}\,.
\end{equation}

With a straightforward computation we get the following equalities, 
\begin{align*}
\left.\frac{{ d} \vert p(t)-q(t)\vert^2}{{ d} t} \right\vert_{t=t_0}
&\,= 2 \langle p-q\,\vert\, k(p)\nu(p) -k(q)\nu(q)\rangle\\
\left.\frac{{ d} A_{p(t),q(t)}}{{ d} t}\right\vert_{t=t_0} &\,=
\int_{\Gamma_{p,q}} \langle\underline{k}(s)\,\vert\nu_{\Gamma_{p,q}}\rangle\,ds
+ 1/2\vert p-q\vert\langle {\nu_{[p,q]}}\,\vert\, k(p)\nu(p)
  +k(q)\nu(q)\rangle\\
&\,= 2\alpha -5\pi/3-1/2\vert p-q\vert(k(p)-k(q))\cos\alpha
\end{align*}
where we wrote $\nu_{\Gamma_{p,q}}$ and $\nu_{[p,q]}$ for the exterior
unit normal to the region ${\mathcal A}_{p,q}$, respectively at the
points of the geodesic $\Gamma_{p,q}$ and of the segment $[p,q]$.\\
Substituting these derivatives in equation~\eqref{eqderE} we get
\begin{align*}
\frac{dE(t_0)}{dt} 
=&\,\frac{1}{A^2_{p,q}}
\Bigl(2A_{p,q}\langle p-q\,\vert\, k(p)\nu(p)
-k(q)\nu(q)\rangle +1/2 \cos\alpha\vert p-q\vert^3(k(p)-k(q))\Bigr)\\
&\,-\frac{\vert p-q\vert^2}{A_{p,q}^2}\Bigl(2\alpha
-\frac{5\pi}{3}\Bigr)\,,
\end{align*}
and, by equation~\eqref{eqfin}, the first term in parentheses is non
negative, hence
$$
\frac{dE(t_0)}{dt} 
\geq -\frac{\vert p-q\vert^2}{A^2_{p,q}}\Bigl(2\alpha
-\frac{5\pi}{3}\Bigr)\,.
$$
By equation~\eqref{topolino2} we have $\alpha=\arctan (-4/E(t_0))$,
hence the following inequality holds
$$
\frac{dE(t_0)}{dt} \geq
\frac{2E(t_0)}{A_{p,q}}\Bigl(\arctan(4/E(t_0))-\arctan(1/\sqrt{3})\Bigr)>0\,.
$$
As the area $A_{p,q}$ is bounded by the area of $\Omega$, we conclude
that for every $t$ in an interval such that the minimum of $\Phi_t$ is
taken by at least one pair of inner points and the derivative of $E(t)$ exists, 
$$
\frac{dE(t_0)}{dt} \geq 
2C E(t_0)\Bigl(\arctan(4/E(t_0))-\arctan(1/\sqrt{3})\Bigr)
$$
for the time independent constant $C=1/{\mathrm {Area}}(\Omega)>0$.
\end{proof}

\begin{rem}\label{netdl} 
  All this analysis can be extended step by step to a network
  with many 3--points but without {\em loops}, that is, a {\em
    tree}. The presence of loops complicates the analysis of the
  minimality properties, because of the possible presence of more than one
  geodesic in $\TTT_t$ between the two points $p$ and $q$.\\
\end{rem}

We now show that on all $[0,T)$ we have $E(t)>0$.

\begin{teo}\label{dlteo} If $\Omega$ is bounded and {\em strictly} convex, 
there exists  a constant $C>0$ depending only on $\TTT_0$ 
  such that $E(t)>C>0$ for every $t\in[0,T)$.\\
Hence, the triods $\TTT_t$ remain embedded in all the maximal interval of 
existence of the flow.
\end{teo}
\begin{proof}
We define three flows of networks of curves 
${\mathbb  H}^1_t$,  ${\mathbb  H}^2_t$, ${\mathbb   H}^3_t$ in the 
interval $[0,T)$. The network ${\mathbb   H}^i_t$ is
obtained as the set theoretic union of $\TTT_t$ with its symmetric
image $\TTT_t^i$ with respect to the point $P^i$.\\
As the triods $\TTT_t$ are contained in the convex
$\Omega$ which is strictly convex, 
this operation does not introduces self--intersections and
since, by Lemma~\ref{evenly} all the even derivatives of $k$ and
$\lambda$ are zero at the end points $P^i$, each one of 
${\mathbb  H}^1_t$,  ${\mathbb  H}^2_t$, ${\mathbb   H}^3_t$
is a smooth flow by curvature of a centrally symmetric network, which
is a tree and it is composed of five curves, two 3--points and four
fixed end points.\\
We define for these networks the functions $E^1, E^2, E^3:[0,T)\to\R$, analogous
to the function $E:[0,T)\to\R$ of $\TTT_t$ and we set
$\Pi(t)=\min\{E^1(t), E^2(t), E^3(t)\}$ which clearly turns out to be
a locally Lipschitz function on $[0,T)$ satisfying $\Pi(t)\leq E^i(t)\leq
E(t)$ for every time $t$ and index $i\in\{1,2,3\}$, since every
${\mathbb  H}^i_t$ contains a copy of $\TTT_t$ (actually two copies). 
Moreover, as there are no self--intersections by construction, $\Pi(0)>0$.\\
Showing that for every time $\Pi(t)>0$, we prove the theorem.\\
We consider a time $t\in[0,T)$ such that the time derivatives of
$\Pi$ and of all the $E^i$ exist (almost everywhere), 
then for every index $i\in\{1, 2, 3\}$ such that $E^i(t)=\Pi(t)$ we
must have $\frac{dE^i(t)}{dt}=\frac{d\Pi(t)}{dt}$.\\
Extending the
previous analysis to the flow ${\mathbb H}^i_t$, which is a tree hence
Remark~\ref{netdl} applies, if the minimum  $E^i(t)=\Pi(t)$ 
is taken by at least one pair of inner points (notice that $P^i$
became an inner point for ${\mathbb H}^i_t$), then this derivative is
positive.\\
If every minimizing pair $(p,q)$ is constituted by two end
points then, by our construction, 
the squared distance $\vert p-q\vert^2$ is bounded
from below and $E^i(t)=\Pi(t)>C>0$ for some uniform constant $C>0$
independent of $t\in[0,T)$, determined by $\TTT_0$.\\
The same holds if one is an end point of ${\TTT}_t$, different from
$P^i$, and the other is a point of $\TTT_t^i$ and {\em viceversa} 
(${\mathbb  H}^i_t$ is centrally symmetric).\\
In the last case, when $(p,q)$ is composed of an inner point and an 
end point of ${\mathbb H}^i_t$, both in the same copy of $\TTT_t$, if
we consider the same pair and associated region, contained in the
other two networks, in at least one of them the two points are both
inner, for instance this happens in the network, say ${\mathbb
  H}^j_t$, if the end point of ${\mathbb H}^i_t$ in the pair is a copy
of $P^j$. Hence, we found an inner minimizing pair of ${\mathbb
  H}^j_t$, that is $E^j(t)=E^i(t)=\Pi(t)$, which implies, by the
previous discussion, that $\frac{dE^i(t)}{dt}>0$ so $\frac{d\Pi(t)}{dt}>0$.

We conclude that if $\Pi(t)$ is under some
constant $C>0$ on $[0,T)$ then it is increasing. Since
$\Pi(0)>0$ this argument gives a uniform bound from below on $\Pi(t)$
on $[0,T)$, hence on $E(t)$.
\end{proof}

\begin{rem} The reason why we put the strict convexity of $\Omega$ in
  the hypothesis is that, in the very special situation such that 
the three end points $P^i$ stay on a line and one of them is
the middle point of the segment determined by the other two, then the
symmetry operation with respect to this middle point produces an
intersection between $\TTT_t$ and its image at the other two end
points, hence the argument above cannot be applied since two loops
have formed.

All the results of this section, in particular the previous one can be
  generalized to networks of curves with many 3--points which are trees,
  that is, without loops in a strictly convex set $\Omega$.\\
We will deal with the case of a general
  network, not necessarily a tree, in a domain $\Omega$ only convex, 
  in a forthcoming paper.
\end{rem}

\section{Blow Up and Self--Similar Solutions}
\label{blw}

As before we suppose to have a smooth embedded solution ${\TTT}_t$ of
Problem~\eqref{problema} in a bounded and strictly convex
$\Omega\subset\R^2$  on a maximal time interval $[0,T)$. Moreover, we
assume that there is a constant $\delta>0$ which uniformly bounds from
below the lengths $L^i(t)$ of the three curves of $\TTT_t$.

By Theorem~\ref{curvexplod} and Corollary~\ref{kexplod} 
the maximum of the modulus of the curvature and its $L^2$ norm go to 
$+\infty$, as $t\to T$.\\
As it is standard, we divide the possible singularities in two cases
(recall Remark~\ref{kratrem}) according to the rate of blow up of the
curvature.\\
We say that we have a {\em Type I} singularity (or a {\em fast}
singularity) if there exists a constant $C$ such that
\begin{equation}\label{typeIeq}
\max_{{\TTT}_t} k^2\leq\frac{C}{{T-t}}
\end{equation}
for every $t\in[0,T)$.\\
In the case inequality~\eqref{typeIeq} does not hold for any constant
$C$ we say that the singularity is of {\em Type II} (or a {\em slow}
singularity).

Blowing up in a proper way the evolving triods around a {\em Type I}
singularity, we obtain as possible limits an unbounded and
embedded triod without end points or an unbounded and embedded curve
with at most one end point,  moving by curvature simply shrinking
homothetically.

In the standard case of the evolution of a closed curve in the plane, 
it is possible to do a blow up of a {\em Type II} singularity obtaining a
{\em translating} solution of the motion by curvature which has to be
a straight line or the {\em grim reaper}  (see~\cite{altsch,hamilton4}). In our
situation dealing with triods, we are not able at the moment to get
the same conclusion, but we can anyway 
hope to exclude the presence of {\em Type II} singularities, 
analysing the possible blow up.

In the rest of this section we classify the embedded 
triods without end points (see Remark~\ref{noendp}) 
shrinking homothetically  by curvature and, since they
could become relevant, we classify also the translating ones.

\begin{dfnz}\label{smttt} 
  We say that a family of unbounded triods without end
  points is a smooth flow by curvature in $\R^2$ 
  if locally, in space and time, it can be parametrized by 
  some maps giving a smooth flow as in Definition~\ref{smoothflow}
  (but with the end points free to move).
\end{dfnz}

\subsection{Classification of Homothetic Flows of Triods}

We can clearly suppose that the origin of $\R^2$ is the center of
homothety.\\
By a straightforward computation, it can be seen that the curves 
$\gamma$ in $\R^2$ which are {\em homothetically shrinking}
around the origin of $\R^2$ under the curvature flow satisfy the
equation $k=-\lambda\langle x\,\vert\,\nu\rangle$, for some {\em
  positive}  constant $\lambda$, at every point $x\in\gamma$.\\
Then, by a rescaling, we can suppose that 
$k=-\langle x\,\vert\,\nu\rangle$, that is $\lambda=1$.

\begin{lemma}\label{lemomt}
If a closed, unbounded and embedded triod $\TTT$ in $\R^2$ without end
points  satisfies $k=-\langle x\,\vert\,\nu\rangle$ at every point
$x\in\TTT$, then the 3--point coincides with the origin of $\R^2$ and
the three curves of the triod are halflines forming angles of $120$ degrees.
\end{lemma}
\begin{proof}
By the work of Abresch and Langer~\cite{ablang1}, since the three
curves of the triod satisfy the equation above, they can be only halflines, 
pieces of circles or pieces of a special family of curves 
({\em  curves of Abresch and Langer}) described in~\cite{ablang1},
which are bounded, periodic and with transversal
self--intersections. Since the triod has no end points, if an 
edge is a piece of one of these curves or of a circle,
following this edge in the direction opposite to the 3--point one
would get a self--intersection which is not present, 
hence, these two possibilities have to be excluded.\\
Then, the three curves are halflines meeting with angles of $120$
degrees and,  by means of the equation above, it follows immediately
that the 3--point has to coincide with the origin of $\R^2$.
\end{proof}

By the same argument used in this proof we can also get the following
two lemmas that we will need in the next section.

\begin{lemma}\label{lemomt2}
If a closed, unbounded and embedded curve $\gamma$ without end points 
satisfies $k=-\langle x\,\vert\,\nu\rangle$ at every point
$x\in\gamma$, then the curve is a straight line through the origin of $\R^2$.
\end{lemma}

\begin{lemma}\label{lemomt3}
If a closed, unbounded and embedded curve $\gamma$ with only one end
satisfies $k=-\langle x\,\vert\,\nu\rangle$ at every point
$x\in\gamma$, then the curve is a halfline.
\end{lemma}

These lemmas implies the following classification result.

\begin{prop}\label{propomt} If $\TTT_t$ is a flow of unbounded
  triods (curves) without end points (with one or without end points)
  shrinking homothetically during the motion by curvature, then 
  $\TTT_t$ is composed of three halflines forming three angles of $120$
  degrees (one halfline or a straight line), hence it is not moving at
  all ($k\equiv0$).
\end{prop}

\subsection{Classification of Translating Flows of Triods}

The curves $\gamma$ in $\R^2$ which move by translation with constant 
velocity $w\in\R^2$, for the curvature flow satisfy the equation
$k=\langle w\,\vert\,\nu\rangle$ at every point 
(observe that this equation is {\em translation invariant}).

\begin{dfnz} The {\em grim reaper} relative to the vector $e_1$ is the
  graph of the function $x=-\log(\cos y)$ in
  $\R^2$ when $y$ varies in the interval $(-\pi/2,\pi/2)$.\\
The {\em grim reaper} relative to a non zero vector $w\in\R^2$ 
is obtained rotating and dilating the {\em grim reaper} relative to $e_1$, in a
way to make this latter coincide with $w$.
\end{dfnz}

\begin{rem} Notice that the {\em grim reaper} relative to $w$ is
  a smooth convex curve asymptotic to two straight lines in $\R^2$ 
  parallel to such vector.
\end{rem}
\begin{figure}[h]
\setlength{\unitlength}{0.240900pt}
\ifx\plotpoint\undefined\newsavebox{\plotpoint}\fi
\begin{picture}(1500,450)(0,0)
\font\gnuplot=cmr10 at 10pt
\gnuplot
\sbox{\plotpoint}{\rule[-0.200pt]{0.400pt}{0.400pt}}%
\put(451,249){\makebox(0,0)[l]{$e_1$}}
\put(1324,366){\makebox(0,0)[l]{$y={\pi}/2$}}
\put(1324,87){\makebox(0,0)[l]{$y=-{\pi}/2$}}
\put(405,225){\vector(1,0){184}}
\put(1216,41){\usebox{\plotpoint}}
\multiput(1085.24,41.60)(-45.810,0.468){5}{\rule{31.500pt}{0.113pt}}
\multiput(1150.62,40.17)(-248.620,4.000){2}{\rule{15.750pt}{0.400pt}}
\multiput(848.45,45.60)(-18.613,0.468){5}{\rule{12.900pt}{0.113pt}}
\multiput(875.23,44.17)(-101.225,4.000){2}{\rule{6.450pt}{0.400pt}}
\multiput(728.75,49.61)(-17.877,0.447){3}{\rule{10.900pt}{0.108pt}}
\multiput(751.38,48.17)(-58.377,3.000){2}{\rule{5.450pt}{0.400pt}}
\multiput(667.68,52.60)(-8.670,0.468){5}{\rule{6.100pt}{0.113pt}}
\multiput(680.34,51.17)(-47.339,4.000){2}{\rule{3.050pt}{0.400pt}}
\multiput(612.66,56.60)(-6.915,0.468){5}{\rule{4.900pt}{0.113pt}}
\multiput(622.83,55.17)(-37.830,4.000){2}{\rule{2.450pt}{0.400pt}}
\multiput(568.40,60.60)(-5.599,0.468){5}{\rule{4.000pt}{0.113pt}}
\multiput(576.70,59.17)(-30.698,4.000){2}{\rule{2.000pt}{0.400pt}}
\multiput(527.32,64.61)(-7.160,0.447){3}{\rule{4.500pt}{0.108pt}}
\multiput(536.66,63.17)(-23.660,3.000){2}{\rule{2.250pt}{0.400pt}}
\multiput(500.96,67.60)(-3.990,0.468){5}{\rule{2.900pt}{0.113pt}}
\multiput(506.98,66.17)(-21.981,4.000){2}{\rule{1.450pt}{0.400pt}}
\multiput(473.79,71.60)(-3.698,0.468){5}{\rule{2.700pt}{0.113pt}}
\multiput(479.40,70.17)(-20.396,4.000){2}{\rule{1.350pt}{0.400pt}}
\multiput(446.41,75.61)(-4.704,0.447){3}{\rule{3.033pt}{0.108pt}}
\multiput(452.70,74.17)(-15.704,3.000){2}{\rule{1.517pt}{0.400pt}}
\multiput(427.87,78.60)(-2.967,0.468){5}{\rule{2.200pt}{0.113pt}}
\multiput(432.43,77.17)(-16.434,4.000){2}{\rule{1.100pt}{0.400pt}}
\multiput(408.11,82.60)(-2.528,0.468){5}{\rule{1.900pt}{0.113pt}}
\multiput(412.06,81.17)(-14.056,4.000){2}{\rule{0.950pt}{0.400pt}}
\multiput(390.53,86.60)(-2.382,0.468){5}{\rule{1.800pt}{0.113pt}}
\multiput(394.26,85.17)(-13.264,4.000){2}{\rule{0.900pt}{0.400pt}}
\multiput(371.73,90.61)(-3.365,0.447){3}{\rule{2.233pt}{0.108pt}}
\multiput(376.36,89.17)(-11.365,3.000){2}{\rule{1.117pt}{0.400pt}}
\multiput(358.77,93.60)(-1.943,0.468){5}{\rule{1.500pt}{0.113pt}}
\multiput(361.89,92.17)(-10.887,4.000){2}{\rule{0.750pt}{0.400pt}}
\multiput(345.19,97.60)(-1.797,0.468){5}{\rule{1.400pt}{0.113pt}}
\multiput(348.09,96.17)(-10.094,4.000){2}{\rule{0.700pt}{0.400pt}}
\multiput(330.39,101.61)(-2.695,0.447){3}{\rule{1.833pt}{0.108pt}}
\multiput(334.19,100.17)(-9.195,3.000){2}{\rule{0.917pt}{0.400pt}}
\multiput(320.02,104.60)(-1.505,0.468){5}{\rule{1.200pt}{0.113pt}}
\multiput(322.51,103.17)(-8.509,4.000){2}{\rule{0.600pt}{0.400pt}}
\multiput(309.02,108.60)(-1.505,0.468){5}{\rule{1.200pt}{0.113pt}}
\multiput(311.51,107.17)(-8.509,4.000){2}{\rule{0.600pt}{0.400pt}}
\multiput(298.43,112.60)(-1.358,0.468){5}{\rule{1.100pt}{0.113pt}}
\multiput(300.72,111.17)(-7.717,4.000){2}{\rule{0.550pt}{0.400pt}}
\multiput(287.60,116.61)(-1.802,0.447){3}{\rule{1.300pt}{0.108pt}}
\multiput(290.30,115.17)(-6.302,3.000){2}{\rule{0.650pt}{0.400pt}}
\multiput(279.85,119.60)(-1.212,0.468){5}{\rule{1.000pt}{0.113pt}}
\multiput(281.92,118.17)(-6.924,4.000){2}{\rule{0.500pt}{0.400pt}}
\multiput(271.26,123.60)(-1.066,0.468){5}{\rule{0.900pt}{0.113pt}}
\multiput(273.13,122.17)(-6.132,4.000){2}{\rule{0.450pt}{0.400pt}}
\multiput(262.16,127.61)(-1.579,0.447){3}{\rule{1.167pt}{0.108pt}}
\multiput(264.58,126.17)(-5.579,3.000){2}{\rule{0.583pt}{0.400pt}}
\multiput(255.68,130.60)(-0.920,0.468){5}{\rule{0.800pt}{0.113pt}}
\multiput(257.34,129.17)(-5.340,4.000){2}{\rule{0.400pt}{0.400pt}}
\multiput(248.68,134.60)(-0.920,0.468){5}{\rule{0.800pt}{0.113pt}}
\multiput(250.34,133.17)(-5.340,4.000){2}{\rule{0.400pt}{0.400pt}}
\multiput(240.71,138.61)(-1.355,0.447){3}{\rule{1.033pt}{0.108pt}}
\multiput(242.86,137.17)(-4.855,3.000){2}{\rule{0.517pt}{0.400pt}}
\multiput(235.09,141.60)(-0.774,0.468){5}{\rule{0.700pt}{0.113pt}}
\multiput(236.55,140.17)(-4.547,4.000){2}{\rule{0.350pt}{0.400pt}}
\multiput(229.51,145.60)(-0.627,0.468){5}{\rule{0.600pt}{0.113pt}}
\multiput(230.75,144.17)(-3.755,4.000){2}{\rule{0.300pt}{0.400pt}}
\multiput(224.09,149.60)(-0.774,0.468){5}{\rule{0.700pt}{0.113pt}}
\multiput(225.55,148.17)(-4.547,4.000){2}{\rule{0.350pt}{0.400pt}}
\multiput(217.82,153.61)(-0.909,0.447){3}{\rule{0.767pt}{0.108pt}}
\multiput(219.41,152.17)(-3.409,3.000){2}{\rule{0.383pt}{0.400pt}}
\multiput(213.92,156.60)(-0.481,0.468){5}{\rule{0.500pt}{0.113pt}}
\multiput(214.96,155.17)(-2.962,4.000){2}{\rule{0.250pt}{0.400pt}}
\multiput(209.92,160.60)(-0.481,0.468){5}{\rule{0.500pt}{0.113pt}}
\multiput(210.96,159.17)(-2.962,4.000){2}{\rule{0.250pt}{0.400pt}}
\multiput(205.37,164.61)(-0.685,0.447){3}{\rule{0.633pt}{0.108pt}}
\multiput(206.69,163.17)(-2.685,3.000){2}{\rule{0.317pt}{0.400pt}}
\multiput(201.92,167.60)(-0.481,0.468){5}{\rule{0.500pt}{0.113pt}}
\multiput(202.96,166.17)(-2.962,4.000){2}{\rule{0.250pt}{0.400pt}}
\multiput(197.92,171.60)(-0.481,0.468){5}{\rule{0.500pt}{0.113pt}}
\multiput(198.96,170.17)(-2.962,4.000){2}{\rule{0.250pt}{0.400pt}}
\multiput(194.95,175.00)(-0.447,0.685){3}{\rule{0.108pt}{0.633pt}}
\multiput(195.17,175.00)(-3.000,2.685){2}{\rule{0.400pt}{0.317pt}}
\multiput(190.92,179.61)(-0.462,0.447){3}{\rule{0.500pt}{0.108pt}}
\multiput(191.96,178.17)(-1.962,3.000){2}{\rule{0.250pt}{0.400pt}}
\put(188.17,182){\rule{0.400pt}{0.900pt}}
\multiput(189.17,182.00)(-2.000,2.132){2}{\rule{0.400pt}{0.450pt}}
\multiput(186.95,186.00)(-0.447,0.685){3}{\rule{0.108pt}{0.633pt}}
\multiput(187.17,186.00)(-3.000,2.685){2}{\rule{0.400pt}{0.317pt}}
\put(183.17,190){\rule{0.400pt}{0.700pt}}
\multiput(184.17,190.00)(-2.000,1.547){2}{\rule{0.400pt}{0.350pt}}
\put(181.17,193){\rule{0.400pt}{0.900pt}}
\multiput(182.17,193.00)(-2.000,2.132){2}{\rule{0.400pt}{0.450pt}}
\put(179.67,197){\rule{0.400pt}{0.964pt}}
\multiput(180.17,197.00)(-1.000,2.000){2}{\rule{0.400pt}{0.482pt}}
\put(178.17,201){\rule{0.400pt}{0.900pt}}
\multiput(179.17,201.00)(-2.000,2.132){2}{\rule{0.400pt}{0.450pt}}
\put(176.67,205){\rule{0.400pt}{0.723pt}}
\multiput(177.17,205.00)(-1.000,1.500){2}{\rule{0.400pt}{0.361pt}}
\put(175.67,208){\rule{0.400pt}{0.964pt}}
\multiput(176.17,208.00)(-1.000,2.000){2}{\rule{0.400pt}{0.482pt}}
\put(174.67,216){\rule{0.400pt}{0.723pt}}
\multiput(175.17,216.00)(-1.000,1.500){2}{\rule{0.400pt}{0.361pt}}
\put(176.0,212.0){\rule[-0.200pt]{0.400pt}{0.964pt}}
\put(174.67,231){\rule{0.400pt}{0.723pt}}
\multiput(174.17,231.00)(1.000,1.500){2}{\rule{0.400pt}{0.361pt}}
\put(175.0,219.0){\rule[-0.200pt]{0.400pt}{2.891pt}}
\put(175.67,238){\rule{0.400pt}{0.964pt}}
\multiput(175.17,238.00)(1.000,2.000){2}{\rule{0.400pt}{0.482pt}}
\put(176.67,242){\rule{0.400pt}{0.723pt}}
\multiput(176.17,242.00)(1.000,1.500){2}{\rule{0.400pt}{0.361pt}}
\put(178.17,245){\rule{0.400pt}{0.900pt}}
\multiput(177.17,245.00)(2.000,2.132){2}{\rule{0.400pt}{0.450pt}}
\put(179.67,249){\rule{0.400pt}{0.964pt}}
\multiput(179.17,249.00)(1.000,2.000){2}{\rule{0.400pt}{0.482pt}}
\put(181.17,253){\rule{0.400pt}{0.900pt}}
\multiput(180.17,253.00)(2.000,2.132){2}{\rule{0.400pt}{0.450pt}}
\put(183.17,257){\rule{0.400pt}{0.700pt}}
\multiput(182.17,257.00)(2.000,1.547){2}{\rule{0.400pt}{0.350pt}}
\multiput(185.61,260.00)(0.447,0.685){3}{\rule{0.108pt}{0.633pt}}
\multiput(184.17,260.00)(3.000,2.685){2}{\rule{0.400pt}{0.317pt}}
\put(188.17,264){\rule{0.400pt}{0.900pt}}
\multiput(187.17,264.00)(2.000,2.132){2}{\rule{0.400pt}{0.450pt}}
\multiput(190.00,268.61)(0.462,0.447){3}{\rule{0.500pt}{0.108pt}}
\multiput(190.00,267.17)(1.962,3.000){2}{\rule{0.250pt}{0.400pt}}
\multiput(193.61,271.00)(0.447,0.685){3}{\rule{0.108pt}{0.633pt}}
\multiput(192.17,271.00)(3.000,2.685){2}{\rule{0.400pt}{0.317pt}}
\multiput(196.00,275.60)(0.481,0.468){5}{\rule{0.500pt}{0.113pt}}
\multiput(196.00,274.17)(2.962,4.000){2}{\rule{0.250pt}{0.400pt}}
\multiput(200.00,279.60)(0.481,0.468){5}{\rule{0.500pt}{0.113pt}}
\multiput(200.00,278.17)(2.962,4.000){2}{\rule{0.250pt}{0.400pt}}
\multiput(204.00,283.61)(0.685,0.447){3}{\rule{0.633pt}{0.108pt}}
\multiput(204.00,282.17)(2.685,3.000){2}{\rule{0.317pt}{0.400pt}}
\multiput(208.00,286.60)(0.481,0.468){5}{\rule{0.500pt}{0.113pt}}
\multiput(208.00,285.17)(2.962,4.000){2}{\rule{0.250pt}{0.400pt}}
\multiput(212.00,290.60)(0.481,0.468){5}{\rule{0.500pt}{0.113pt}}
\multiput(212.00,289.17)(2.962,4.000){2}{\rule{0.250pt}{0.400pt}}
\multiput(216.00,294.61)(0.909,0.447){3}{\rule{0.767pt}{0.108pt}}
\multiput(216.00,293.17)(3.409,3.000){2}{\rule{0.383pt}{0.400pt}}
\multiput(221.00,297.60)(0.774,0.468){5}{\rule{0.700pt}{0.113pt}}
\multiput(221.00,296.17)(4.547,4.000){2}{\rule{0.350pt}{0.400pt}}
\multiput(227.00,301.60)(0.627,0.468){5}{\rule{0.600pt}{0.113pt}}
\multiput(227.00,300.17)(3.755,4.000){2}{\rule{0.300pt}{0.400pt}}
\multiput(232.00,305.60)(0.774,0.468){5}{\rule{0.700pt}{0.113pt}}
\multiput(232.00,304.17)(4.547,4.000){2}{\rule{0.350pt}{0.400pt}}
\multiput(238.00,309.61)(1.355,0.447){3}{\rule{1.033pt}{0.108pt}}
\multiput(238.00,308.17)(4.855,3.000){2}{\rule{0.517pt}{0.400pt}}
\multiput(245.00,312.60)(0.920,0.468){5}{\rule{0.800pt}{0.113pt}}
\multiput(245.00,311.17)(5.340,4.000){2}{\rule{0.400pt}{0.400pt}}
\multiput(252.00,316.60)(0.920,0.468){5}{\rule{0.800pt}{0.113pt}}
\multiput(252.00,315.17)(5.340,4.000){2}{\rule{0.400pt}{0.400pt}}
\multiput(259.00,320.61)(1.579,0.447){3}{\rule{1.167pt}{0.108pt}}
\multiput(259.00,319.17)(5.579,3.000){2}{\rule{0.583pt}{0.400pt}}
\multiput(267.00,323.60)(1.066,0.468){5}{\rule{0.900pt}{0.113pt}}
\multiput(267.00,322.17)(6.132,4.000){2}{\rule{0.450pt}{0.400pt}}
\multiput(275.00,327.60)(1.212,0.468){5}{\rule{1.000pt}{0.113pt}}
\multiput(275.00,326.17)(6.924,4.000){2}{\rule{0.500pt}{0.400pt}}
\multiput(284.00,331.61)(1.802,0.447){3}{\rule{1.300pt}{0.108pt}}
\multiput(284.00,330.17)(6.302,3.000){2}{\rule{0.650pt}{0.400pt}}
\multiput(293.00,334.60)(1.358,0.468){5}{\rule{1.100pt}{0.113pt}}
\multiput(293.00,333.17)(7.717,4.000){2}{\rule{0.550pt}{0.400pt}}
\multiput(303.00,338.60)(1.505,0.468){5}{\rule{1.200pt}{0.113pt}}
\multiput(303.00,337.17)(8.509,4.000){2}{\rule{0.600pt}{0.400pt}}
\multiput(314.00,342.60)(1.505,0.468){5}{\rule{1.200pt}{0.113pt}}
\multiput(314.00,341.17)(8.509,4.000){2}{\rule{0.600pt}{0.400pt}}
\multiput(325.00,346.61)(2.695,0.447){3}{\rule{1.833pt}{0.108pt}}
\multiput(325.00,345.17)(9.195,3.000){2}{\rule{0.917pt}{0.400pt}}
\multiput(338.00,349.60)(1.797,0.468){5}{\rule{1.400pt}{0.113pt}}
\multiput(338.00,348.17)(10.094,4.000){2}{\rule{0.700pt}{0.400pt}}
\multiput(351.00,353.60)(1.943,0.468){5}{\rule{1.500pt}{0.113pt}}
\multiput(351.00,352.17)(10.887,4.000){2}{\rule{0.750pt}{0.400pt}}
\multiput(365.00,357.61)(3.365,0.447){3}{\rule{2.233pt}{0.108pt}}
\multiput(365.00,356.17)(11.365,3.000){2}{\rule{1.117pt}{0.400pt}}
\multiput(381.00,360.60)(2.382,0.468){5}{\rule{1.800pt}{0.113pt}}
\multiput(381.00,359.17)(13.264,4.000){2}{\rule{0.900pt}{0.400pt}}
\multiput(398.00,364.60)(2.528,0.468){5}{\rule{1.900pt}{0.113pt}}
\multiput(398.00,363.17)(14.056,4.000){2}{\rule{0.950pt}{0.400pt}}
\multiput(416.00,368.60)(2.967,0.468){5}{\rule{2.200pt}{0.113pt}}
\multiput(416.00,367.17)(16.434,4.000){2}{\rule{1.100pt}{0.400pt}}
\multiput(437.00,372.61)(4.704,0.447){3}{\rule{3.033pt}{0.108pt}}
\multiput(437.00,371.17)(15.704,3.000){2}{\rule{1.517pt}{0.400pt}}
\multiput(459.00,375.60)(3.698,0.468){5}{\rule{2.700pt}{0.113pt}}
\multiput(459.00,374.17)(20.396,4.000){2}{\rule{1.350pt}{0.400pt}}
\multiput(485.00,379.60)(3.990,0.468){5}{\rule{2.900pt}{0.113pt}}
\multiput(485.00,378.17)(21.981,4.000){2}{\rule{1.450pt}{0.400pt}}
\multiput(513.00,383.61)(7.160,0.447){3}{\rule{4.500pt}{0.108pt}}
\multiput(513.00,382.17)(23.660,3.000){2}{\rule{2.250pt}{0.400pt}}
\multiput(546.00,386.60)(5.599,0.468){5}{\rule{4.000pt}{0.113pt}}
\multiput(546.00,385.17)(30.698,4.000){2}{\rule{2.000pt}{0.400pt}}
\multiput(585.00,390.60)(6.915,0.468){5}{\rule{4.900pt}{0.113pt}}
\multiput(585.00,389.17)(37.830,4.000){2}{\rule{2.450pt}{0.400pt}}
\multiput(633.00,394.60)(8.670,0.468){5}{\rule{6.100pt}{0.113pt}}
\multiput(633.00,393.17)(47.339,4.000){2}{\rule{3.050pt}{0.400pt}}
\multiput(693.00,398.61)(17.877,0.447){3}{\rule{10.900pt}{0.108pt}}
\multiput(693.00,397.17)(58.377,3.000){2}{\rule{5.450pt}{0.400pt}}
\multiput(774.00,401.60)(18.613,0.468){5}{\rule{12.900pt}{0.113pt}}
\multiput(774.00,400.17)(101.225,4.000){2}{\rule{6.450pt}{0.400pt}}
\multiput(902.00,405.60)(45.810,0.468){5}{\rule{31.500pt}{0.113pt}}
\multiput(902.00,404.17)(248.620,4.000){2}{\rule{15.750pt}{0.400pt}}
\put(176.0,234.0){\rule[-0.200pt]{0.400pt}{0.964pt}}
\put(65,410){\usebox{\plotpoint}}
\put(65.0,410.0){\rule[-0.200pt]{329.792pt}{0.400pt}}
\put(65,40){\usebox{\plotpoint}}
\put(65.0,40.0){\rule[-0.200pt]{329.792pt}{0.400pt}}
\end{picture}
\caption{\label{figura2} The {\em grim reaper} relative to $e_1$.}
\end{figure}
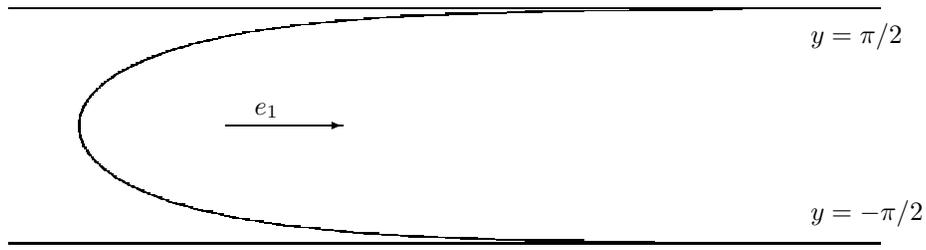

\begin{lemma}\label{lemtra}
If a closed, unbounded and embedded triod $\TTT$ in $\R^2$ without 
end points satisfies $k=\langle w\,\vert\,\nu\rangle$ at every point,
with $w\not=0$, then its curves are halflines parallel to $w$ or
translated copies of pieces of
the {\em grim reaper} relative to $w$ (see Figure~\ref{figura2}),
concurring at the 3--point with angles of $120$
degrees (it clearly follows that at most one curve is a halfline).
\end{lemma}
\begin{figure}[h]
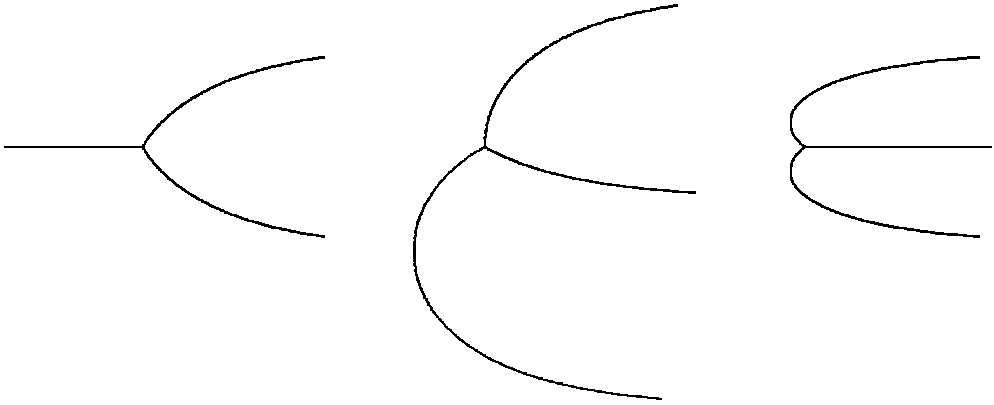

\setlength{\unitlength}{0.240900pt}
\ifx\plotpoint\undefined\newsavebox{\plotpoint}\fi

\caption{\label{figura3} Some examples of translating triods.}
\end{figure}
\begin{proof}
At every inner point of the curves composing the triod,
differentiating in arclength the equation above we get
$k_s=-k\langle w\,\vert\,\tau\rangle$, hence, if the curvature is zero
at some point, it has to be zero on all the curve, so this latter is  a
halfline parallel to $w$.\\
If $k$ is non zero, then differentiating again we obtain
$k_{ss}=-k_s\langle w\,\vert\,\tau\rangle-k^2\langle
w\,\vert\,\nu\rangle=k_s^2/k-k^3$ which is the equation of a 
{\em  grim reaper}. Integrating such equation we get the thesis.
\end{proof}

With the same arguments we can also prove the following lemmas.

\begin{lemma}\label{lemtra2}
If a closed, unbounded and embedded curve $\gamma$ in $\R^2$ without 
end points satisfies $k=\langle w\,\vert\,\nu\rangle$ at every point, 
then it is either a straight line parallel to the vector
$w$ or a translated copy of the {\em grim reaper} relative to $w$.
\end{lemma}

\begin{lemma}\label{lemtra3}
If a closed, unbounded and embedded curve $\gamma$ with only one end 
satisfies $k=\langle w\,\vert\,\nu\rangle$ at every point, 
then it is either a halfline parallel to the vector $w$
or a translated copy of a piece of the {\em grim reaper} relative to $w$.
\end{lemma}

Hence, we conclude as above.

\begin{prop}\label{proptra} If $\TTT_t$ is a flow of unbounded
  triods (curves) without end points (without or with only one end point)
  translating during the motion by curvature, then every $\TTT_t$ is one of
  the sets of Lemma~\ref{lemtra} (\ref{lemtra2} or~\ref{lemtra3}).
\end{prop}

\section{Type I Singularities}
\label{typeIsecn}

In this section we study the {\em Type I} singularities of a smooth
flow of triods $\TTT_t$, described by a map
$F:\TTT\times[0,T)\to\overline{\Omega}$ 
(see after Definition~\ref{geocomp}), in a bounded and strictly convex
$\Omega\subset\R^2$, on a maximal time interval $[0,T)$.\\
We remark that in this section, we will often consider 
the triods as subsets of $\R^2$, without mentioning $\Omega$.

By the {\em Type I} assumption, there exists a constant $C>0$ such that
\begin{equation}\label{stimadallalto}
\lim_{t\to T}\max_{\TTT_t} k^2=+\infty\qquad \text{ and }\qquad
k^2(p,t)\leq \frac{C}{T-t}
\end{equation}
for every $p\in\TTT$ and $t\in[0,T)$. 

\begin{lemma}\label{casi}
One of the following three (mutually non exclusive) possibilities holds:
\begin{enumerate}

\item There is an index $i\in\{1,2,3\}$ and a sequence of times
   $t_j\nearrow T$ such that 
\begin{equation*}
[k^i(O,t_j)]^2\geq \frac{C}{T-t_j}
\end{equation*}
for some positive constant $C$.

\item We have the estimate
$$
\max_{\TTT_t}k^2=\frac{o(1)}{T-t}
$$
as $t\to T$.

\item The maximum of $k^2$ is definitely not taken at the point $O$ 
and 
$$
\max_{\TTT_t}k^2\geq \frac{1}{2(T-t)}
$$
for every $t$ larger than some $t_0\in[0,T)$.
\end{enumerate}
\end{lemma}
\begin{proof}
We consider the non negative, locally Lipschitz functions
$f(t)=\max_{\TTT_t}k^2$ and $g(t)=(T-t)\max_{\TTT_t}k^2$.\\
If we are not in the first case, $\lim_{t\to T}(T-t)[k^i(O,t)]^2=0$
for every index $i\in\{1, 2, 3\}$. 
Setting $g_O(t)=\max_{i\in\{1,2,3\}}(T-t)[k^i(O,t)]^2$, if
$g(t)>g_O(t)$ then the maximum of 
$k^2$ at time $t$ is taken in the interior
of one of the three curves ($k$ is zero at the three end points).\\
If $g>g_O$ does not hold definitely, there exists a time $t_0$ such
that $g_O(t)<1/4$ for every $t>t_0$ and another time $t_1>t_0$ such that
$g(t_1)=g_O(t_1)<1/4$.\\
Then, following Huisken~\cite{huisk3}, at every time $t$ such that
$g(t)>g_O(t)$ by the parabolic maximum 
principle the function $f(t)=\max_{{\TTT}_t} k^2$ satisfies the
differential inequality 
\begin{equation}\label{maxkappa3}
\frac{d\,}{dt}f(t)\leq 2 \max_{\TTT_t}k^4 =2f^2(t)
\end{equation}
and the function $g(t)$,
$$
\frac{d\,}{dt}g(t)\leq 2 (T-t)\max_{\TTT_t}k^4
-\max_{\TTT_t}k^2 =(2g-1)f(t)\,.
$$
If $J$ is the set of times $t\in[t_1,T)$ such that $g(t)=g_O(t)<1/4$, 
then in $[t_1,T)\setminus J$, which is a union of open intervals, 
the function $g$ is decreasing, by a simple ODE's argument. 
It follows that $g(t)\leq g_O(\sup (J\cap [0,t]))$ 
and since we supposed that $\lim_{t\to T}\sup (J\cap[0,t])=T$, we
conclude $\lim_{t\to T}g(t)=0$. This means that we are in
the second case.\\
If instead $g>g_O$ definitely, then after some time the maximum of
$k^2$ is always taken at some inner points, hence 
inequality~\eqref{maxkappa3} holds and since $f$ goes to $+\infty$ as
$t\to T$, by integration we get the bound from below
(Huisken~\cite{huisk3})
$$
\max_{\TTT_t}k^2\geq \frac{1}{2(T-t)}
$$
for every $t$ larger than some $t_0$, that is, we are in the third
case.
\end{proof}

In the second situation above we can get something more.

\begin{lemma}\label{preterzo} If the curvature of the triods satisfies
$$
\max_{\TTT_t}k^2=\frac{o(1)}{T-t}\,.
$$
then there exists an index $i\in\{1,2,3\}$ and a sequence of times
$t_j\nearrow T$ such that
$$
k^2(p,t)\leq [k^i(O,t_j)]^2\, \text{ for every $t\leq t_j$,
  $p\in\TTT$ and } [k^i(O,t_j)]^2\nearrow+\infty\,.
$$
\end{lemma}
\begin{proof}
We start studying the non decreasing hull 
$\widetilde{f}(t)=\sup_{\xi\in[0,t]} f(\xi)$ of the function
$f(t)=\max_{{\TTT}_t} k^2=\varepsilon(t)/(T-t)$, with $\lim_{t\to
  T}\varepsilon(t)=0$. 
Clearly $\widetilde{f}(t)\geq f(t)$, hence also $\widetilde{f}$ goes to $+\infty$
as $t\to T$.\\
We notice that since $f$ is locally Lipschitz, the same hold for
$\widetilde{f}$, moreover the derivative of this latter, which 
exists at almost every time $t\in[0,T)$, is zero or coincides with
the derivative of $f$, if this happen $\widetilde{f}(t)=f(t)$ also
holds.\\
We saw in inequality~\eqref{maxkappa3} 
that when the maximum of $k^2$ is taken at some points
different from $O$, then $f^\prime\leq 2f^2$, 
so under the same hypothesis, it also holds $\widetilde{f}^\prime\leq
f^\prime\leq 2f^2\leq 2\widetilde{f}^2$.\\
Let us consider the set $J=\{t\in[0,T)\,\vert\,
\widetilde{f}(t)=\max_{i\in\{1, 2, 3\}}[k^i(O,t)]^2\}$, 
if a value $t$ does not belong to $J$ then either the maximum of
$k^2$ is taken at some points different from $O$ or
$f(t)<\widetilde{f}(t)$ and the derivative of $\widetilde{f}$ at $t$
is zero. In both cases the inequality $\widetilde{f}^\prime\leq
2\widetilde{f}^2$ holds at every time such that the derivatives of $f$
and $\widetilde{f}$ exist (almost everywhere).\\
If $T$ is not a limit point of $J$, then definitely $\widetilde{f}^\prime\leq
2\widetilde{f}^2$ almost everywhere, integrating then 
the distributional derivative of $1/\widetilde{f}$, 
we obtain, as in the previous lemma for $f$, the bound from below
$\widetilde{f}(t)\geq \frac{1}{2(T-t)}$.\\
This gives a contradiction since it implies that definitely
$\widetilde{f}>f$, hence $\widetilde{f}$ would be constant against the
fact that $f(t)$ goes to $+\infty$ as $t\to T$.\\
Thus, we can assume that there exists an index $i\in\{1,2,3\}$ and a 
sequence of times $t_j\nearrow T$ belonging to $J$, hence 
satisfying $\widetilde{f}(t_j)=[k^i(O,t_j)]^2$. It follows, by
construction, that 
$$
k^2(p,t)\leq [k^i(O,t_j)]^2\, \text{ for every $t\leq t_j$ and  $p\in\TTT$,}
$$
and $[k^i(O,t_j)]^2\nearrow+\infty$.
\end{proof}

We will deal with the first and the last case of Lemma~\ref{casi} 
by means of the rescaling procedure of Huisken~\cite{huisk3},
adapting the line of Stone in~\cite{stone1,stone2} to our situation.\\
By means of an argument of Ilmanen in~\cite[Lemma~8]{ilman3}
and~\cite[Section~3]{ilman4}, we will show in
Proposition~\ref{nosing0} that the second case also cannot happen.

We start with the analogue in our situation of Huisken's 
{\em monotonicity formula} for mean curvature flow 
(see~\cite{huisk3}).

Let $x_0\in\R^2$ and define the {\em backward heat kernel} relative to
$(x_0,T)$ as
$$
\rho_{x_0}(x,t)=\frac{e^{-\frac{\vert
    x-x_0\vert^2}{4(T-t)}}}{\sqrt{4\pi(T-t)}}\,.
$$

\begin{lemma}\label{monforuno}
For every $t\in[0,T)$ and $i\in\{1,2,3\}$ the following identity holds
\begin{align*}
\frac{d\,}{dt}\int_{{\gamma^i(\cdot , t)}} \rho_{x_0}(x,t)\,ds=
&\,-\int_{{\gamma^i(\cdot , t)}}
\left\vert\,\underline{k}+\frac{(x-x_0)^{\perp}}{2(T-t)}\right\vert^2
\rho_{x_0}(x,t)\,ds\\
&\,+\left\langle\,\frac{P^i-x_0}{2(T-t)}\,\biggl\vert\,\tau^i(1,t)\right\rangle
\rho_{x_0}(P^i,t)\\
&\,+\left\langle\,\frac{\gamma^i(0,t)-x_0}{2(T-t)}\,\biggl\vert\,\tau^i(0,t)\right\rangle
\rho_{x_0}(\gamma^i(0,t),t)\\
&\,+\lambda^i(0,t)\rho_{x_0}(\gamma^i(0,t),t)\,.
\end{align*}
\end{lemma}
\begin{proof}
The proof goes like in~\cite[Theorem~3.1]{huisk3} with the extra
boundary terms coming from the application of the {\em first variation
formula} (see~\cite{simon}).\\
By such formula, for every $C^1$ vector field $X$ we have
$$
\int_{\gamma^i(\cdot , t)}\div^\top\!X\,ds
=-\int_{\gamma^i(\cdot ,  t)}\langle
X\,\vert\,\underline{k}\,\rangle\,ds + \langle X(\gamma^i(1,t))\,\vert\tau^i(1,t)\rangle -
\langle X(\gamma^i(0,t))\,\vert\tau^i(0,t)\rangle
$$
where $\div^\top$ means {\em tangential divergence}
(see~\cite{simon}), and following Huisken~\cite{huisk3},
\begin{align*}
\frac{d\,}{dt}\int_{{\gamma^i(\cdot , t)}} \rho_{x_0}\,ds
=&\,\int_{{\gamma^i(\cdot , t)}} \rho_{x_0}\div^\top\!\underline{v}\,ds+
\int_{{\gamma^i(\cdot , t)}}
\langle\nabla\rho_{x_0}\,\vert\,\underline{v}\,\rangle+\frac{\partial\rho_{x_0}}{\partial
  t}\,ds\\
=&\,-\int_{{\gamma^i(\cdot , t)}} \rho_{x_0}k^2 - \div^\top(\rho_{x_0}\lambda\tau)\,ds+
\int_{{\gamma^i(\cdot , t)}}\langle\nabla^\perp\rho_{x_0}\,\vert\,\underline{v}\,\rangle
+\frac{\partial\rho_{x_0}}{\partial  t}\,ds\\
=&\,-\int_{{\gamma^i(\cdot , t)}} \rho_{x_0}\left\{k^2
  -\frac{1}{2(T-t)}+\frac{\langle
    x-x_0\,\vert\,\underline{k}\,\rangle}{2(T-t)}+\frac{\vert
    x-x_0\vert^2}{4(T-t)^2}\right\}\,ds\\
&\, - \lambda^i(0,t)\rho_{x_0}(\gamma^i(0,t),t)\\
=&\,-\int_{{\gamma^i(\cdot , t)}} \left\vert\,\underline{k} +
  \frac{x-x_0}{2(T-t)}\right\vert^2\rho_{x_0}\,ds
- \lambda^i(0,t)\rho_{x_0}(\gamma^i(0,t),t)\\
&\, +\int_{{\gamma^i(\cdot , t)}}\frac{\rho_{x_0}}{2(T-t)}+\frac{\langle
  x-x_0\,\vert\,\underline{k}\,\rangle}{2(T-t)}\rho_{x_0}\,ds\\
\end{align*}
since $\lambda^i(1,t)=0$.\\
Now we use again the first variation formula with the field
$X(x)=(x-x_0)\rho_{x_0}(x,t)/2(T-t)$ whose tangential divergence is given
explicitly by
\begin{align*}
\frac{\sum_{l=1}^2\nabla_j^\top\bigl((x^l-x_0^l)\rho_{x_0}(x,t)\bigr)}{2(T-t)}=
&\, \frac{\langle x-x_0\,\vert\,\nabla^\top\rho_{x_0}(x,t)\rangle + 
\sum_{j=1}^2\langle e_j\,\vert
e_j^\top\rangle\rho_{x_0}(x,t)}{2(T-t)}\\
=&\,\frac{\rho_{x_0}(x,t)}{2(T-t)}-\frac{\vert
  (x-x_0)^\top\vert^2}{4(T-t)^2}\rho_{x_0}(x,t)
\end{align*}
then,
\begin{align*}
\int_{\gamma^i(\cdot , t)}\div^\top\!X\,ds=
&\,-\int_{\gamma^i(\cdot , t)}\frac{\langle
  x-x_0\,\vert\,\underline{k}\,\rangle}{2(T-t)}\rho_{x_0}\,ds\\
&\, +\frac{\langle P^i-x_0\,\vert\tau^i(1,t)\rangle\rho_{x_0}(P^i,t) -
  \langle\gamma^i(0,t)-x_0\,\vert\tau^i(0,t)\rangle\rho_{x_0}(\gamma^i(0,t),t)}
{2(T-t)}\\
=&\, \int_{\gamma^i(\cdot , t)}\frac{\rho_{x_0}}{2(T-t)}-\frac{\vert
  (x-x_0)^\top\vert^2}{4(T-t)^2}\rho_{x_0}\,ds\,.
\end{align*}
Hence,
\begin{align*}
\frac{d\,}{dt}\int_{{\gamma^i(\cdot , t)}} \rho_{x_0}\,ds
=&\,-\int_{{\gamma^i(\cdot , t)}} \left\vert\,\underline{k} +
  \frac{x-x_0}{2(T-t)}\right\vert^2\rho_{x_0}\,ds
- \lambda^i(0,t)\rho_{x_0}(\gamma^i(0,t),t)\\
&\, +\frac{\langle P^i-x_0\,\vert\tau^i(1,t)\rangle\rho_{x_0}(P^i,t) -
  \langle\gamma^i(0,t)-x_0\,\vert\tau^i(0,t)\rangle\rho_{x_0}(\gamma^i(0,t),t)}
{2(T-t)}\\
&\, + \int_{\gamma^i(\cdot , t)}\frac{\vert
  (x-x_0)^\top\vert^2}{4(T-t)^2}\rho_{x_0}\,ds\\
=&\,-\int_{{\gamma^i(\cdot , t)}} \left\vert\,\underline{k} +
  \frac{(x-x_0)^\perp}{2(T-t)}\right\vert^2\rho_{x_0}\,ds
- \lambda^i(0,t)\rho_{x_0}(\gamma^i(0,t),t)\\
&\, +\frac{\langle P^i-x_0\,\vert\tau^i(1,t)\rangle\rho_{x_0}(P^i,t) -
  \langle\gamma^i(0,t)-x_0\,\vert\tau^i(0,t)\rangle\rho_{x_0}(\gamma^i(0,t),t)}
{2(T-t)}
\end{align*}
and reordering the terms we get the claimed identity.
\end{proof}

\begin{prop}[Monotonicity Formula]
For every $t\in[0,T)$ the following identity holds
\begin{align}
\frac{d\,}{dt}\int_{{\TTT_t}} \rho_{x_0}(x,t)\,ds=
&\,-\int_{{\TTT_t}} \left\vert\,\underline{k}+\frac{(x-x_0)^{\perp}}{2(T-t)}\right\vert^2
\rho_{x_0}(x,t)\,ds\label{eqmonfor}\\
&\,+\sum_{i=1}^3\left\langle\,\frac{P^i-x_0}{2(T-t)}\,\biggl\vert\,\tau^i(1,t)\right\rangle
\rho_{x_0}(P^i,t)\,.\nonumber
\end{align}
Integrating between $t_1$ and $t_2$ with $0\leq t_1\leq t_2<T$ we get
\begin{align}
\int_{t_1}^{t_2}\int_{{\TTT_t}}
\left\vert\,\underline{k}+\frac{(x-x_0)^{\perp}}{2(T-t)}\right\vert^2
\rho_{x_0}(x,t)\,ds\,dt = 
&\,\int_{{\TTT_{t_1}}} \rho_{x_0}(x,t_1)\,ds -
\int_{{\TTT_{t_2}}} \rho_{x_0}(x,t_2)\,ds\label{eqmonfor-int}\\
&\,+\sum_{i=1}^3\int_{t_1}^{t_2}
\left\langle\,\frac{P^i-x_0}{2(T-t)}\,\biggl\vert\,\tau^i(1,t)\right\rangle
\rho_{x_0}(P^i,t)\,dt\,.\nonumber
\end{align}
\end{prop}
\begin{proof}
We simply add the contributions for the three curves given by
Lemma~\ref{monforuno}, recalling that $\gamma^i(0,t)=O$ and $\sum_{i=1}^3
\tau^i(0,t)=\sum_{i=1}^3 \lambda^i(0,t)=0$.
\end{proof}

\begin{lemma}\label{stimadib} Setting $\vert P^i- x_0\vert=d^i$, 
for every index $i\in\{1,2,3\}$ the following estimate holds 
\begin{equation*}
\left\vert\int_{t}^{T}\left\langle\,\frac{P^i-x_0}{2(T-\xi)}\,\biggl\vert\,
\tau^i(1,\xi)\right\rangle\rho_{x_0}(P^i,\xi)\,d\xi\,\right\vert\leq
\frac{1}{\sqrt{2\pi}}\int\limits_{d^i/\sqrt{2(T-t)}}^{+\infty}e^{-y^2/2}\,dy\leq
1/2\,.
\end{equation*}
Then, for every $x_0\in\R^2$,
\begin{equation*}
\lim_{t\to  T}\sum_{i=1}^3\int_{t}^{T}
\left\langle\,\frac{P^i-x_0}{2(T-\xi)}\,\biggl\vert\,\tau^i(1,\xi)\right\rangle
\rho_{x_0}(P^i,\xi)\,d\xi=0\,.
\end{equation*}
\end{lemma}
\begin{proof}
If $d^i>0$, we estimate
\begin{align*}
\left\vert\int_{t}^{T}\left\langle\,\frac{P^i-x_0}{2(T-\xi)}\,\biggl\vert\,\tau^i(1,\xi)\right\rangle
\rho_{x_0}(P^i,\xi)\,d\xi\,\right\vert\leq
&\,\int_{t}^{T}\left\vert\left\langle\,\frac{P^i-x_0}{2(T-\xi)}\,\biggl\vert\,\tau^i(1,\xi)\right\rangle
\right\vert\rho_{x_0}(P^i,\xi)\,d\xi\\
\leq &\,\frac{1}{\sqrt{2\pi}}\int_{t}^{T}\frac{d^i}{[2(T-\xi)]^{3/2}}e^{-\frac{[d^i]^2}{4(T-\xi)}}\,d\xi\\
= &\,\frac{1}{\sqrt{2\pi}}\int\limits_{d^i/\sqrt{2(T-t)}}^{+\infty}e^{-y^2/2}\,dy\\
\leq &\,\frac{1}{\sqrt{2\pi}}\int_0^{+\infty}e^{-y^2/2}\,dy=1/2
\end{align*}
where we could change variable $y=d^i/\sqrt{2(T-\xi)}$ since
$d^i>0$.\\
Since the integral on $[0,+\infty)$ is finite, looking at the third line we
have also that 
$$
\lim_{t\to T}
\left\vert\int_{t}^{T}\left\langle\,\frac{P^i-x_0}{2(T-\xi)}\,\biggl\vert\,\tau^i(1,\xi)
\right\rangle\rho_{x_0}(P^i,\xi)\,d\xi\,\right\vert=0\,.
$$
In the special case $d^i=0$, that is, $x_0$ coincides with the end
point $P^i$,  the corresponding integral is zero for every
$t\in[0,T)$, so the thesis follows.
\end{proof}

\begin{prop}\label{valori0}
For every $x_0\in\R^2$ the limit $\lim_{t\to
  T}\int_{\TTT_t}\rho_{x_0}(x,t)\,ds$ there exists.
\end{prop}
\begin{proof}
Fixed $x_0\in\R^2$ we look at the function $b:[0,T)\to\R$
$$
b(t)=\int_t^T\sum_{i=1}^3
\left\langle\,\frac{P^i-x_0}{2(T-\xi)}\,\biggl\vert\,\tau^i(1,\xi)\right\rangle
\rho_{x_0}(P^i,\xi)\,d\xi\,.
$$
As Lemma~\ref{stimadib} says that $b$ is bounded and 
$\lim_{t\to T}b(t)=0$, the monotonicity formula~\eqref{eqmonfor}
implies that the limit of the statement there exists.
\end{proof}

Now, we  introduce the the rescaling procedure of
Huisken~\cite{huisk3}.\\
Fixed $x_0\in\R^2$, let $\widetilde{F}_{x_0}:\TTT\times
[-1/2\log{T},+\infty)\to\R^2$ be the map
$$
\widetilde{F}_{x_0}(p,\tt)=\frac{F(p,t)-x_0}{\sqrt{2(T-t)}}\qquad
\tt(t)=-\frac{1}{2}\log{(T-t)}
$$
then, the rescaled triods are given by 
$$
\widetilde{\TTT}_{x_0,\tt}=\frac{\TTT_t-x_0}{\sqrt{2(T-t)}}
$$
and they evolve according to the equation 
$$
\frac{\partial\,}{\partial
  \tt}\widetilde{F}_{x_0}(p,\tt)=\widetilde{\underline{v}}(p,\tt)+\widetilde{F}_{x_0}(p,\tt)
$$
where 
$$
\widetilde{\underline{v}}(p,\tt)=\frac{\underline{v}(p,t(\tt))}{\sqrt{2(T-t(\tt))}}=
\widetilde{\underline{k}}+\widetilde{\underline{\lambda}}=
\widetilde{k}\nu+\widetilde{\lambda}\tau\qquad \text{ and }\qquad
t(\tt)=T-e^{-2\tt}\,.
$$
Notice that we did not put the ``tilde'' over the unit tangent and
normal, since they remain the same in the rescaling.\\
We will often write $\widetilde{O}(\tt)=\widetilde{F}_{x_0}(O,\tt)$ 
for the 3--point of the rescaled triod $\widetilde{\TTT}_{x_0,\tt}$,
when there is no ambiguity on the point $x_0$.\\
The rescaled curvature evolves according to the following equation,
\begin{equation}\label{evolriscforf}
{\partial_\tt} \widetilde{k}=
\widetilde{k}_{\sigma\sigma}+\widetilde{k}_\sigma\widetilde{\lambda} +
\widetilde{k}^3 -\widetilde{k}
\end{equation}
which can be obtained as in Section~\ref{kestimates} by means of the
commutation law
\begin{equation}\label{commutforf}
{\partial_\tt}{\partial_\sigma}={\partial_\sigma}{\partial_\tt} + (\widetilde{k}^2
-\widetilde{\lambda}_\sigma-1){\partial_\sigma}\,,
\end{equation}
where we denoted with $\sigma$ the arclength parameter for
$\widetilde{\TTT}_{x_0,\tt}$.

By a straightforward computation
(\cite{huisk3},~\cite[Lemma~2.3]{stone2}) we have the following
rescaled version of the monotonicity formula.

\begin{prop}[Rescaled Monotonicity Formula]
Let $x_0\in\R^2$ and set  
$$
\widetilde{\rho}(x)=e^{-\frac{\vert x\vert^2}{2}}
$$
For every $\tt\in[-1/2\log{T},+\infty)$ the following identity holds
\begin{equation}
\frac{d\,}{d\tt}\int_{\widetilde{\TTT}_{x_0,\tt}}
\widetilde{\rho}(x)\,d\sigma=
-\int_{\widetilde{\TTT}_{x_0,\tt}}\vert
\,\widetilde{\underline{k}}+x^\perp\vert^2\widetilde{\rho}(x)\,d\sigma
+\sum_{i=1}^3\left\langle\,{\widetilde{P}^i_{x_0,\tt}}
\,\Bigl\vert\,{\tau}^i(1,t(\tt))\right\rangle
\widetilde{\rho}(\widetilde{P}^i_{x_0,\tt})\label{reseqmonfor}
\end{equation}
where $\widetilde{P}^i_{x_0,\tt}=\frac{P^i-x_0}{\sqrt{2(T-t(\tt))}}$.\\  
Integrating between $\tt_1$ and $\tt_2$ with 
$-1/2\log{T}\leq \tt_1\leq \tt_2<+\infty$ we get
\begin{align} 
\int_{\tt_1}^{\tt_2}\int_{\widetilde{\TTT}_{x_0,\tt}}\vert
\,\widetilde{\underline{k}}+x^\perp\vert^2\widetilde{\rho}(x)\,d\sigma\,d\tt=
&\, \int_{\widetilde{\TTT}_{x_0,\tt_1}}\widetilde{\rho}(x)\,d\sigma - 
\int_{\widetilde{\TTT}_{x_0,\tt_2}}\widetilde{\rho}(x)\,d\sigma\label{reseqmonfor-int}\\
&\,+\sum_{i=1}^3\int_{\tt_1}^{\tt_2}\left\langle\,{\widetilde{P}^i_{x_0,\tt}}\,
\Bigl\vert\,{\tau}^i(1,t(\tt))\right
\rangle\widetilde{\rho}(\widetilde{P}^i_{x_0,\tt})\,d\tt\,.\nonumber
\end{align}
\end{prop}

Then, we have the analog of Lemma~\ref{stimadib} whose proof follows
in the same way, substituting the rescaled quantities.

\begin{lemma}\label{rescstimadib} For every index $i\in\{1,2,3\}$ 
the following estimate holds 
\begin{equation*}
\left\vert\int_{\tt}^{+\infty}\left\langle\,{\widetilde{P}^i_{x_0,\xi}}\,
\Bigl\vert\,{\tau}^i(1,t(\xi))\right
\rangle\widetilde{\rho}(\widetilde{P}^i_{x_0,\xi})\,d\xi\right\vert\leq \sqrt{\pi/2}\,.
\end{equation*}
Then, for every $x_0\in\R^2$,
\begin{equation*}
\lim_{\tt\to  +\infty}\sum_{i=1}^3\int_{\tt}^{+\infty}\left\langle\,{\widetilde{P}^i_{x_0,\xi}}\,
\Bigl\vert\,{\tau}^i(1,t(\xi))\right
\rangle\widetilde{\rho}(\widetilde{P}^i_{x_0,\xi})\,d\xi=0\,.
\end{equation*}
\end{lemma}

We need the following lemmas in order to study the possible limits of
the rescaled triods.

\begin{lemma}\label{Obehav}
Under the {\em Type I} hypothesis~\eqref{stimadallalto} there exists
$\lim_{t\to T}F(O,t)=\widehat{O}\in\R^2$.\\
The 3--point $\widetilde{F}_{x_0}(O,\tt)$ of the rescaled triods
either it is uniformly bounded or it goes to infinity as
$\tt\to+\infty$, according to the
fact that $x_0=\widehat{O}$ or not.
\end{lemma}
\begin{proof}
Since at the 3--point we have
$\sum_{i=1}^3({\lambda}^i)^2=\sum_{i=1}^3({k}^i)^2\leq C/(T-t)$ we have also
$\vert\underline v(O,t)\vert^2\leq C/(T-t)$ for some constant $C>0$
independent of $t\in[0,T)$. Then we get
$$
\vert F(O,t_1)-F(O,t_2)\vert=\left\vert\int_{t_1}^{t_2}\frac{\partial F}{\partial
    t}\,dt\right\vert\leq \int_{t_1}^{t_2}\vert
\underline{v}(O,t) \vert\,dt\leq\int_{t_1}^{t_2}\sqrt{\frac{C}{T-t}}\,dt
\leq 2\sqrt{C(T-t_1)}
$$
for every $0\leq t_1\leq t_2<T$. 
Hence, the limit $\lim_{t\to T}F(O,t)=\widehat{O}$ there exists and
\begin{equation}\label{distanza-limo5}
\vert F(O,t)-\widehat{O}\,\vert\leq 2\sqrt{C(T-t)}
\end{equation}
for every $t$.\\
Considering then the 3--points $\widetilde{F}_{\widehat{O}}(O,\tt)$ of
the triods  $\widetilde{\TTT}_{\widehat{O},\tt}$, we have 
$$
\vert\widetilde{F}_{\widehat{O}}(O,\tt)\vert=\frac{\vert
F(O,t)-\widehat{O}\vert}{\sqrt{2(T-t)}}
\leq \frac{2\sqrt{C(T-t)}}{\sqrt{2(T-t)}}=\sqrt{2C}
$$
hence, if $x_0=\widehat{O}$ then the 
3--point always belongs to the ball $B_{\sqrt{2C}}$.\\
By the same inequality~\eqref{distanza-limo5} it follows that if 
$x_0\not=\widehat{O}$ then the 3--points of the rescaled triods go to
infinity as $\tt\to+\infty$.
\end{proof}

\begin{lemma}\label{rescestim}
The curvature $\widetilde{k}$ of the rescaled triods
$\widetilde{\TTT}_{x_0,\tt}$ is uniformly bounded
in space and time. Moreover, for every ball $B_R$ centered at the
origin of $\R^2$, we have the following
estimates with a constant $C_R$ independent of 
$x_0\in\R^2$ and $\tt\in[-1/2\log{T},+\infty)$
$$
\HH^1({\widetilde{\TTT}_{x_0,\tt}}\cap B_R)\leq C_R\,.
$$
\end{lemma}
\begin{proof}
The maximum of the curvature of every rescaled triod is bounded by a 
uniform constant, by the 
assumption~\eqref{stimadallalto} on the blow up rate of the
curvature.\\
Moreover, using the rescaled monotonicity formula~\eqref{reseqmonfor}
we get
\begin{align*} 
\int_{\widetilde{\TTT}_{x_0,\tt}}\widetilde{\rho}\,d\sigma=&\,
\int\limits_{\widetilde{\TTT}_{x_0,-1/2\log{T}}}\widetilde{\rho}\,d\sigma
-\int\limits_{-1/2\log{T}}^{\tt}\int\limits_{\widetilde{\TTT}_{x_0,\xi}}
\vert\,\widetilde{\underline{k}}+x^\perp\vert^2\widetilde{\rho}\,d\sigma\,d\xi\\
&\,+\sum_{i=1}^3\int\limits_{-1/2\log{T}}^{\tt}\left\langle\,
{\widetilde{P}^i_{x_0,\xi}}\,\Bigl\vert\,{\tau}^i(1,t(\xi))\right\rangle
\widetilde{\rho}(\widetilde{P}^i_{x_0,\xi})\,d\xi\leq C
\end{align*}
where the last estimate follows from Lemma~\ref{rescstimadib}.\\
Hence, since $\widetilde{\rho}\geq e^{-R^2/2}$ in every ball $B_R$ 
centered at the origin of $\R^2$, we have a uniform bound
$\HH^1({\widetilde{\TTT}_{x_0,\tt}}\cap B_R)\leq C_R$ for some
constants $C_R$ independent of $\tt$ and $x_0$.
\end{proof}

\begin{dfnz} We say that a sequence of triods converges in the
  $C^r\loc$ topology if, after reparametrizing the curves composing
  the triods in arclength, they converge in $C^r$ in every compact of
  $\R^2$.\\
  The definition of convergence in $W^{n,p}\loc$ is analogous.
\end{dfnz}

\begin{lemma}\label{elemma} The function (see Section~\ref{disuellesec})
\begin{equation*}
E(\TTT)= \inf_{\genfrac{}{}{0pt}{}
{p,q\in\TTT}{p\not=q}}\frac{\vert p-q\vert^2}{A_{p,q}}\,,
\end{equation*}
defined on the class of $C^2$ triods without self--intersections 
(bounded or unbounded and with or without end points), is upper
semicontinuous with respect to the $C^2\loc$ convergence.\\
Moreover, $E$ is {\em dilation and translation invariant}.\\
Hence, every $C^2\loc$ limit of rescaled triods satisfies $E>C>0$
where the uniform positive constant $C$ is given by
Theorem~\ref{dlteo}. Thus, such limit has no self--intersections.
\end{lemma}
\begin{proof}
The proof is straightforward.
\end{proof}

\begin{lemma}\label{convo1}
For every $x_0\in\R^2$ and
every sequence of rescaled times $\tt_j\to\infty$ there exists a
subsequence $\tt_{j_l}$ such that the triods
$\widetilde{\TTT}_{x_0,\tt_{j_l}}$ converge in the 
$C^1\loc$ topology to a limit set
$\TTT_\infty$ which, {\em if non empty}, is a curve or a triod with at
most one end point and without self--intersections.\\
Moreover, the Radon measures $\HH^1\res
\widetilde{\TTT}_{x_0,\tt_{j_l}}$ weakly$^{{\displaystyle{\star}}}$ converge in $\R^2$
  to the Radon measure $\HH^1\res\TTT_\infty$.
\end{lemma}
\begin{proof}
Reparametrizing the triods in arclength, we have curves with uniformly
bounded first and second derivatives, moreover, by Lemma~\ref{rescestim}, 
for every ball $B_R$ centered at the origin of $\R^2$ we have a uniform bound
$\HH^1({\widetilde{\TTT}_{x_0,\tt}}\cap B_R)\leq C_R$ for some
constants $C_R$ independent of $\tt$.\\
Then, by standard compactness arguments (see~\cite{huisk3,langer2,simon}), 
the sequence ${\widetilde{\TTT}_{x_0,\tt_{j}}}$ of reparametrized
triods has a subsequence ${\widetilde{\TTT}_{x_0,\tt_{j_l}}}$
weakly$^{{\displaystyle{\star}}}$ converging in  $W^{2,\infty}\loc$, 
then in the $C^1\loc$ topology to a (possibly empty) set
$\TTT_{\infty}$ which, if $x_0$ is distinct from all the $P^i$, 
has no end points since they go to infinity during the rescaled flow. 
If $x_0=P^i$, the set $\TTT_\infty$ has a single end point at the
origin of $\R^2$.\\
In both cases the 3--point could be present or not, if it is present then
the angles are of $120$ degrees by the convergence of the
curves in $C^1\loc$. The only ``strange'' situation is if $x_0=P^1$, for
instance, and $\lim_{l\to\infty}\widetilde{O}(\tt_{j_l})=0$, 
which is in contradiction with the fact that the three 
lengths are uniformly positively bounded
from below, indeed in this situation the curve between $P^1$ and the
3--point has to collapse otherwise embeddedness, which we are going
to show now, is lost.\\
The limit set, which we suppose non empty, 
has no self--intersections by Lemma~\ref{elemma}.\\
Finally, the embeddedness of the limit and the $C^1$ convergence in
every compact imply that the Radon measures $\HH^1\res
\widetilde{\TTT}_{x_0,\tt_{j_l}}$ weakly$^{{\displaystyle{\star}}}$
converge in $\R^2$  to the Radon measure $\HH^1\res\TTT_\infty$.
\end{proof}

\begin{lemma}\label{deltabuono0} If the 3--point of the rescaled triod 
$\widetilde{\TTT}_{x_0,\tt}$ does not belong to the ball
$B_{2R}(z)\subset\R^2$ with radius $2R$ and center $z\in\R^2$, then
there exist constants $\delta_R>0$ and $D_R$, independent of 
  the points $z, x_0\in\R^2$ and the time 
  $\tt\in[-1/2\log{T},+\infty)$, such that 
$\widetilde{k}_\sigma$, $\widetilde{k}_{\sigma\sigma}$ and 
$\widetilde{k}_\tt$ for the family 
$\bigl\{\widetilde{\TTT}_{x_0,r}\,\bigl\vert\,
  r\in[\tt,\tt+\delta_R]\bigl\}$ are uniformly bounded by $D_R$ 
in the smaller ball $B_R(z)$.
\end{lemma}
\begin{proof}
By the control on the curvature, the velocity of the 3--point is
bounded by a uniform constant, hence, for some $\delta_R>0$ in the
time interval $[\tt-\delta_R,\tt+\delta_R]$ the 3--point cannot enter in the
ball $B_{3R/2}(z)$.\\
Then, as the 3--points are far from the ball, by the uniform bound on the
curvature of all the rescaled triods, repeating the interior estimates
of Ecker and Huisken in~\cite{eckhui2}  (see
also~\cite[Section~2]{eck1}) for the rescaled flow and recalling the
evolution equation for the rescaled curvature 
${\partial_\tt}\widetilde{k}=\widetilde{k}_{\sigma\sigma}+\widetilde{\eta}
\widetilde{k}_\sigma+\widetilde{k}^3-\widetilde{k}$, we get the
thesis, possibly choosing a smaller $\delta_R>0$.
\end{proof}

\begin{rem}\label{deltabuono2} 
  The same conclusion clearly holds for a family of
  triods moving by curvature (not rescaled) if we have a uniform 
  bound on the curvature.
\end{rem} 
 
\begin{prop}\label{resclimit}
For every $x_0\in\R^2$ and
every sequence of rescaled times $\tt_j\to\infty$ there exists a
subsequence $\tt_{j_l}$ such that the triods
$\widetilde{\TTT}_{x_0,\tt_{j_l}}$ converge in the 
$C^2\loc$ topology to a limit set
$\TTT_\infty$ which, {\em if non empty}, is one of the following:
\begin{itemize}
\item a triod composed of three halflines originating from $0\in\R^2$, 
\item a halfline from $0\in\R^2$,
\item a straight line from $0\in\R^2$.
\end{itemize}
Moreover, the Radon measures $\HH^1\res
\widetilde{\TTT}_{x_0,\tt_{j_l}}$ weakly$^{{\displaystyle{\star}}}$ converge in $\R^2$
  to the Radon measure $\HH^1\res\TTT_\infty$.
\end{prop}
\begin{proof}
By Lemma~\ref{convo1}, we have to show that, supposing the $C^1\loc$--limit
$\TTT_\infty=\lim_{l\to\infty}{\widetilde{\TTT}_{x_0,\tt_{j_l}}}$ non 
empty, it is among the ones of the statement.

Putting $\tt_1=-1/2\log T$ and sending $\tt_2$ to $+\infty$ in 
the rescaled monotonicity formula~\ref{reseqmonfor-int},  by
Lemma~\ref{stimadib} we get
$$
\int\limits_{-1/2\log{T}}^{+\infty}\int\limits_{\widetilde{\TTT}_{x_0,\tt}}
\vert\,\widetilde{\underline{k}}+x^\perp\vert^2\widetilde{\rho}\,d\sigma\,d\tt<+\infty\,,
$$
hence, extracting from the sequence of times $\tt_{j_l}$ a 
subsequence (not relabelled) with $\tt_{j_{l+1}}>\tt_{j_l}+1/l$, 
we see that there exists an increasing sequence $r_{j_l}$ such that
$\tt_{j_l}\leq r_{j_l}\leq \tt_{j_l}+1/l$ and
on a subsequence of the $r_{j_l}$ (again not relabelled) we
have
$$
\lim_{l\to\infty}\int\limits_{\widetilde{\TTT}_{x_0,r_{j_l}}}
\vert\,\widetilde{\underline{k}}+x^\perp\vert^2\widetilde{\rho}\,d\sigma^1=0\,.
$$
Reapplying Lemma~\ref{convo1}, we can assume that also the triods
$\widetilde{\TTT}_{x_0,r_{j_l}}$ converges to some limit
  $\overline{\TTT}_\infty$ in $W^{2,\infty}\loc$ (possibly empty), 
  and since the integral above is lower semicontinuous in this convergence
  (see~\cite{simon}), the limit $\overline{\TTT}_\infty$ satisfies 
  $\widetilde{\underline{k}}+x^\perp=0$, distributionally. 
  The limit set is composed of curves in
  $W^{2,\infty}\loc$ but from the relation
  $\widetilde{\underline{k}}+x^\perp=0$ it follows that 
  $\widetilde{\underline{k}}$ is continuous, since the curves are
  $C^1\loc$.\\
  Such limit set is a triod or a curve and the end point is present or not 
  according to the choice of the point $x_0$.\\
  As the relation above implies $\widetilde{k}=-\langle
  x\,\vert\,\nu\rangle$ at every point $x\in\overline{\TTT}_\infty$, 
  the classification Lemmas~\ref{lemomt},~\ref{lemomt2},~\ref{lemomt3} 
  show that in any case the curvature of the limit set is zero 
  everywhere and that $\overline{\TTT}_\infty$ is among the sets of
  the statement. Indeed, when such limit is an halfline, we have 
  necessarily that the point $x_0$ coincides with one of the
  fixed end points of the triods, hence the limit halfline starts from
  the origin of $\R^2$.

We show now that ${\TTT}_\infty=\overline{\TTT}_\infty$ and that the
convergence is in the $C^2\loc$ topology, proving the
proposition.\\
We consider a point $y\in\TTT_{\infty}$ such that the distance
$\dist(y,\overline{\TTT}_\infty)>0$.
If we denote with $y_\tt=\widetilde{F}_{x_0}(p_\tt,\tt)$ 
the point of minimum distance from $y$ in 
the rescaled triod $\widetilde{\TTT}_{x_0,\tt}$ and we look at the
function
$f(\tt)=\dist(y,\widetilde{\TTT}_{x_0,\tt})=\vert y - y_\tt\vert$ we
have (with the usual remarks about differentiability),
$$
\frac{df(\tt)}{d\tt}=\frac{\Bigl\langle y-y_\tt\,\Bigr\vert\,
\frac{\partial\widetilde{F}_{x_0}(p_\tt,\tt)}{\partial\tt}\Bigl\rangle}{\vert y-y_\tt\vert}
=\frac{\Bigl\langle y-y_\tt\,\Bigr\vert\,
    \widetilde{\underline{v}}(p_\tt,\tt)+
\widetilde{F}_{x_0}(p_\tt,\tt)\Bigr\rangle}{\vert y-y_\tt\vert}
\leq C\vert\widetilde{k}\vert+\vert\widetilde{F}_{x_0}(p_\tt,\tt)\vert
$$
where we substituted the velocity with the curvature since, if
$y_\tt$ is in the interior, the vector $(y-y_\tt)$ is orthogonal to
$\widetilde{\TTT}_{x_0,\tt}$ by minimality, if $y_\tt$ is an end point
then $\widetilde{\underline{v}}(p_\tt,\tt)=0$, finally if $y_\tt$ is the 3--point, by
the usual relation
$\sum_{i=1}^3(\widetilde{\lambda}^i)^2=\sum_{i=1}^3(\widetilde{k}^i)^2$,
the velocity is controlled by a constant multiple of the curvature.\\
Since $\widetilde{k}$ is uniformly bounded and the triangular
inequality gives $\vert\widetilde{F}_{x_0}(p_\tt,\tt)\vert=\vert y_\tt\vert
\leq \vert y\vert+\vert y-y_\tt\vert\leq C+f(\tt)$, we conclude
$$
\frac{df(\tt)}{d\tt}\leq f(\tt)+C\,.
$$
Integrating this differential inequality on the interval
$[\tt_{j_l},r_{j_l}]$ we get 
$$
f(r_{j_l})\leq
e^{r_{j_l}-\tt_{j_l}}(f(\tt_{j_l})+C)-C\leq 
e^{1/l}f(\tt_{j_l}) + C(e^{1/l}-1)
$$
so, if $l\to\infty$, since we know that
$\lim_{l\to\infty}f(\tt_{j_l})=0$ we have also 
$\lim_{l\to\infty}f(r_{j_l})=0$, thus $y\in\overline{\TTT}_\infty$.\\
This clearly implies that $\overline{\TTT}_\infty$ cannot be empty,
then, inverting the roles of the two limit sets and repeating this
argument we conclude that they must coincide.

Now we show the $C^2\loc$ convergence.\\
If the limit set $\TTT_\infty$ is a straight line, by the $C^1$
convergence in every ball $B_R$ and the uniform bound on the curvature
$\widetilde{k}$, we can apply Lemma~\ref{deltabuono0} to get 
a uniform bound on the norm
$\Vert\widetilde{k}_\sigma\Vert_{L^\infty(B_R)}$ 
for all the triods $\widetilde{\TTT}_{x_0,\tt_{j_l}}$, 
hence the $C^2\loc$ convergence follows.\\
With the same argument we have a uniform bound on
$\widetilde{k}_\sigma$ in the subset $B_{nR}\setminus B_{2R}$, for every
$n\in\NN$ greater than 2, in the case that the 3--point of
$\widetilde{\TTT}_{x_0,\tt}$ definitely belongs to the ball $B_R$ 
(see Lemma~\ref{Obehav}), or the limit $\TTT_\infty$
is a halfline from the origin of $\R^2$.\\
Then, we work out a local version of the estimates leading to
Proposition~\ref{topolino5} in order to deal with these two situations.
By means of equations~\eqref{evolriscforf} and~\eqref{commutforf} we
have
\begin{align*}
{\partial_\tt} \widetilde{k}=
&\,\widetilde{k}_{\sigma\sigma}+\widetilde{k}_\sigma\widetilde{\lambda} +
\widetilde{k}^3 -\widetilde{k}\\
{\partial_\tt} \widetilde{k_\sigma}=
&\,\widetilde{k}_{\sigma\sigma\sigma}+\widetilde{k}_{\sigma\sigma}\widetilde{\lambda}
+ 4\widetilde{k}^2\widetilde{k}_\sigma - 2\widetilde{k}_\sigma\\
\partial_\tt\widetilde{k}_{\sigma\sigma}=
&\,\widetilde{k}_{\sigma\sigma\sigma\sigma}
+\widetilde{\lambda}\widetilde{k}_{\sigma\sigma\sigma}
+5\widetilde{k}^2\widetilde{k}_{\sigma\sigma}
+8\widetilde{k}\widetilde{k}^2_{\sigma}
-3\widetilde{k}_{\sigma\sigma}
\end{align*}
hence, for every smooth function $\varphi:\R^2\times[0,\delta]\to[0,1]$ with
compact support we compute, 
\begin{align*}
\frac{d\,}{d\tt}\int\limits_{\widetilde{\TTT}_{x_0,\tt}}
&\,\left(\widetilde{k}^2+\tt\widetilde{k}_\sigma^2+\tt^2\widetilde{k}_{\sigma\sigma}^2/2\right)
\varphi^2\,d\sigma\\
=&\,-\int\limits_{\widetilde{\TTT}_{x_0,\tt}}
\left(\widetilde{k}^2_\sigma+\tt\widetilde{k}_{\sigma\sigma}^2
+\tt^2\widetilde{k}_{\sigma\sigma\sigma}^2\right)
\varphi^2\,d\sigma\\
&\,-\varphi^2(\widetilde{O},\tt)\sum_{i=1}^32\widetilde{k}\widetilde{k}_\sigma+
2\tt\widetilde{k}_{\sigma}\widetilde{k}_{\sigma\sigma}+
\tt^2\widetilde{k}_{\sigma\sigma}\widetilde{k}_{\sigma\sigma\sigma}
\,\biggr\vert_{\text{{      at the 3--point}}}\\
 &\,+\int\limits_{\widetilde{\TTT}_{x_0,\tt}}
\left(2\widetilde{k}\widetilde{k}_\sigma+
2\tt\widetilde{k}_{\sigma}\widetilde{k}_{\sigma\sigma}+
\tt^2\widetilde{k}_{\sigma\sigma}\widetilde{k}_{\sigma\sigma\sigma}\right)
\widetilde{\lambda}\varphi^2\,d\sigma\\
&\,+\int\limits_{\widetilde{\TTT}_{x_0,\tt}}
\left(\widetilde{k}^2+\tt\widetilde{k}_{\sigma}^2
+\tt^2\widetilde{k}_{\sigma\sigma}^2/2\right)
\left(\varphi^2\widetilde{\lambda}_\sigma+
2\varphi\langle\nabla\varphi\,\vert\,\tau\rangle
\widetilde{\lambda}\right)\,d\sigma\\
&\,+\int\limits_{\widetilde{\TTT}_{x_0,\tt}}
\left(\widetilde{k}^2+\tt\widetilde{k}_{\sigma}^2
+\tt^2\widetilde{k}_{\sigma\sigma}^2/2\right)
\left(1-\widetilde{k}^2\right)\varphi^2\,d\sigma\\
&\,+2\int\limits_{\widetilde{\TTT}_{x_0,\tt}}
\left[\widetilde{k}^4-\widetilde{k}^2+\tt\left(4\widetilde{k}^2\widetilde{k}_{\sigma}^2
    -2\widetilde{k}_{\sigma}^2\right)
+\frac{\tt^2}{2}\left(8\widetilde{k}\widetilde{k}_\sigma^2\widetilde{k}_{\sigma\sigma}-
3\widetilde{k}_{\sigma\sigma}^2\right)\right]\varphi^2\,d\sigma\\
&\,+\int\limits_{\widetilde{\TTT}_{x_0,\tt}}
\left(\widetilde{k}^2+\tt\widetilde{k}_{\sigma}^2
+\tt^2\widetilde{k}_{\sigma\sigma}^2/2\right)
\left(\langle\nabla\varphi^2\,\vert\,\widetilde{\underline{k}}\rangle+2\varphi\varphi_\tt
  \right)\,d\sigma
\end{align*}
where we 
already forgot the end points contributions, by Lemma~\ref{evenly}.\\
After integrating by parts the terms containing $\widetilde{\lambda}$,
if
$\langle\nabla\varphi^2\,\vert\,\widetilde{\underline{k}}\rangle+2\varphi\varphi_\tt$
is non positive for $\tt\in[0,\delta]$ and taking into account that
$\widetilde{k}$ is bounded, we get
\begin{align*}
\frac{d\,}{d\tt}\int\limits_{\widetilde{\TTT}_{x_0,\tt}}
&\,\left(\widetilde{k}^2+\tt\widetilde{k}_\sigma^2+\tt^2\widetilde{k}_{\sigma\sigma}^2/2\right)
\varphi^2\,d\sigma\\
\leq&\,-\int\limits_{\widetilde{\TTT}_{x_0,\tt}}
\left(\widetilde{k}^2_\sigma+\tt\widetilde{k}_{\sigma\sigma}^2
+\tt^2\widetilde{k}_{\sigma\sigma\sigma}^2\right)
\varphi^2\,d\sigma\\
&\,-\varphi^2(\widetilde{O},\tt)\sum_{i=1}^32\widetilde{k}\widetilde{k}_\sigma+
2\tt\widetilde{k}_{\sigma}\widetilde{k}_{\sigma\sigma}+
\tt^2\widetilde{k}_{\sigma\sigma}\widetilde{k}_{\sigma\sigma\sigma}
\,\biggr\vert_{\text{{      at the 3--point}}}\\
&\,-\varphi^2(\widetilde{O},\tt)\sum_{i=1}^3\widetilde{k}^2\widetilde{\lambda} +
\tt\widetilde{k}_{\sigma}^2\widetilde{\lambda}+
\tt^2\widetilde{k}_{\sigma\sigma}^2\widetilde{\lambda}/2
\,\biggr\vert_{\text{{      at the 3--point}}}\\
&\,+C\int\limits_{\widetilde{\TTT}_{x_0,\tt}}
\left(\widetilde{k}^2+\tt\widetilde{k}_{\sigma}^2
+\tt^2\widetilde{k}_{\sigma\sigma}^2/2\right)\varphi^2\,d\sigma\\
&\,+2\int\limits_{\widetilde{\TTT}_{x_0,\tt}}
\left[\widetilde{k}^4-\widetilde{k}^2+\tt\left(4\widetilde{k}^2\widetilde{k}_{\sigma}^2
    -2\widetilde{k}_{\sigma}^2\right)
+\frac{\tt^2}{2}\left(8\widetilde{k}\widetilde{k}_\sigma^2\widetilde{k}_{\sigma\sigma}-
3\widetilde{k}_{\sigma\sigma}^2\right)\right]\varphi^2\,d\sigma\,.
\end{align*}
In doing this we choose a function $\varphi$ as follows. If the
3--point of the rescaled triods is bounded, according to
Lemma~\ref{Obehav}, and it is definitely contained in the ball of radius
$R$, supposing that the curvature is uniformly bounded by some
constant $C>0$, we set
$\varphi(x,\tt)=\sqrt{h(\psi(x)-C\tt\Vert\nabla\psi\Vert_{L^\infty})}$
choosing a radially monotone and symmetric 
smooth function $\psi:\R^2\to[0,1]$, with compact support, 
such that $\psi\vert_{B_{2R}}=1$ and
$\psi\vert_{\R^2\setminus B_{3R}}=0$ and a smooth 
increasing function $h:\R\to[0,1]$ such that $h(z)=0$
if $z\leq 0$ and $h(z)=1$ for $z\geq1$. Then, for $\delta>0$ small enough, 
depending only on $\Vert\nabla\psi\Vert_{L^\infty}$ and $C$, 
the function $\varphi$ satisfies the requirements above and possibly
choosing a smaller value $\delta>0$, there holds 
$\varphi\geq 1/2$ on $B_{2R}$ for every
$\tt\in[0,\delta]$.\\
Since by Lemma~\ref{Obehav} the 3--point of the rescaled triods is
definitely inside or outside the ball $B_{R}$ we
can consider $\varphi^2(\widetilde{O},\tt)$ constantly equal to one or zero.\\
Hence, if the 3--point is present, we are dealing
with the case where $\TTT_\infty$ is an unbounded triod, 
so, by the $C^1\loc$ convergence, the length of the curves in the ball
$B_{2R}$ are bounded from below by $R$ and we can
treat the 3--point term as before in proving
Proposition~\ref{topolino5}. Then, denoting with 
${\mathcal I}$ the boundary term, we get
\begin{align*}
\frac{d\,}{d\tt}\int\limits_{\widetilde{\TTT}_{x_0,\tt}}
&\,\left(\widetilde{k}^2+\tt\widetilde{k}_\sigma^2+\tt^2\widetilde{k}_{\sigma\sigma}^2/2\right)
\varphi^2\,d\sigma\\
\leq&\,-\int\limits_{\widetilde{\TTT}_{x_0,\tt}}
\left(\widetilde{k}^2_\sigma+\tt\widetilde{k}_{\sigma\sigma}^2
+\tt^2\widetilde{k}_{\sigma\sigma\sigma}^2\right)
\varphi^2\,d\sigma + {\mathcal I}\\
&\,+C\int\limits_{\widetilde{\TTT}_{x_0,\tt}}
\left(\widetilde{k}^2+\tt\widetilde{k}_{\sigma}^2
+\tt^2\widetilde{k}_{\sigma\sigma}^2/2\right)\varphi^2\,d\sigma\\
&\,+2\int\limits_{\widetilde{\TTT}_{x_0,\tt}}
\left[\widetilde{k}^4-\widetilde{k}^2+\tt\left(4\widetilde{k}^2\widetilde{k}_{\sigma}^2
    -2\widetilde{k}_{\sigma}^2\right)
+\frac{\tt^2}{2}\left(8\widetilde{k}\widetilde{k}_\sigma^2\widetilde{k}_{\sigma\sigma}-
3\widetilde{k}_{\sigma\sigma}^2\right)\right]\varphi^2\,d\sigma\\
\leq&\,-\int\limits_{\widetilde{\TTT}_{x_0,\tt}}
\left(\widetilde{k}^2_\sigma+\tt\widetilde{k}_{\sigma\sigma}^2
+\tt^2\widetilde{k}_{\sigma\sigma\sigma}^2\right)
\varphi^2\,d\sigma+{\mathcal I}+C\\
&\,+C\int\limits_{\widetilde{\TTT}_{x_0,\tt}}
\left(\widetilde{k}^2+\tt\widetilde{k}_{\sigma}^2
+\tt^2\widetilde{k}_{\sigma\sigma}^2/2\right)\varphi^2\,d\sigma\\
&\,+C\int\limits_{\widetilde{\TTT}_{x_0,\tt}}
\left(\tt\widetilde{k}_{\sigma}^2
+\tt^2\widetilde{k}_\sigma^2\widetilde{k}_{\sigma\sigma}
\right)\varphi^2\,d\sigma\\
\leq&\,-\int\limits_{\widetilde{\TTT}_{x_0,\tt}}
\left(\widetilde{k}^2_\sigma+\tt\widetilde{k}_{\sigma\sigma}^2
+\tt^2\widetilde{k}_{\sigma\sigma\sigma}^2\right)
\varphi^2\,d\sigma +{\mathcal I}+ C\\
&\,+C\int\limits_{\widetilde{\TTT}_{x_0,\tt}}
\left(\widetilde{k}^2+\tt\widetilde{k}_{\sigma}^2
+\tt^2\widetilde{k}_{\sigma\sigma}^2/2\right)\varphi^2\,d\sigma\\
&\,+C\int\limits_{\widetilde{\TTT}_{x_0,\tt}}
\left(\tt\widetilde{k}_{\sigma}^2
+\tt^2\widetilde{k}_{\sigma\sigma}^2 +\tt^2\widetilde{k}_{\sigma}^4
\right)\varphi^2\,d\sigma
\end{align*}
and choosing $\delta<1$ such that $\delta C<1/2$, it follows
\begin{align*}
\frac{d\,}{d\tt}\int\limits_{\widetilde{\TTT}_{x_0,\tt}}
\left(\widetilde{k}^2+\tt\widetilde{k}_\sigma^2+\tt^2\widetilde{k}_{\sigma\sigma}^2/2\right)
\varphi^2\,d\sigma\leq&\,-1/2\int\limits_{\widetilde{\TTT}_{x_0,\tt}}
\left(\widetilde{k}^2_\sigma+\tt\widetilde{k}_{\sigma\sigma}^2
+\tt^2\widetilde{k}_{\sigma\sigma\sigma}^2\right)
\varphi^2\,d\sigma + {\mathcal I}+C\\
&\,+C\int\limits_{\widetilde{\TTT}_{x_0,\tt}}
\left(\widetilde{k}^2+\tt\widetilde{k}_{\sigma}^2
+\tt^2\widetilde{k}_{\sigma\sigma}^2/2\right)\varphi^2\,d\sigma\\
&\,+C\int\limits_{\widetilde{\TTT}_{x_0,\tt}}
\tt^2\widetilde{k}_{\sigma}^4\varphi^2\,d\sigma\,.
\end{align*}
We break the last term into
\begin{align*}
\int\limits_{\widetilde{\TTT}_{x_0,\tt}}
\tt^2\widetilde{k}_{\sigma}^4\varphi^2\,d\sigma=
&\,\int\limits_{\widetilde{\TTT}_{x_0,\tt}\cap B_R}
\tt^2\widetilde{k}_{\sigma}^4\varphi^2\,d\sigma
+\int\limits_{\widetilde{\TTT}_{x_0,\tt}\cap (B_{3R}\setminus B_R)}
\tt^2\widetilde{k}_{\sigma}^4\varphi^2\,d\sigma\\
=&\,\int\limits_{\widetilde{\TTT}_{x_0,\tt}\cap B_R}
\tt^2\widetilde{k}_{\sigma}^4\,d\sigma
+\int\limits_{\widetilde{\TTT}_{x_0,\tt}\cap (B_{3R}\setminus B_R)}
\tt^2\widetilde{k}_{\sigma}^4\varphi^2\,d\sigma
\end{align*}
as $\varphi$ is zero outside the ball $B_{3R}$.\\
The second integral is bounded since, by the argument based on the
interior estimates of Ecker and Huisken, 
discussed previously, $\widetilde{k}_\sigma$ is uniformly bounded
for ${\widetilde{\TTT}_{x_0,\tt}\cap (B_{3R}\setminus B_R)}$.\\
The first integral is controlled by
interpolating between a possibly large multiple of 
$\int\limits_{\widetilde{\TTT}_{x_0,\tt}\cap B_R}
\tt^2\widetilde{k}^2\,d\sigma$, which is bounded, and a small fraction
of $\int\limits_{\widetilde{\TTT}_{x_0,\tt}\cap B_R}
\tt^2\widetilde{k}_{\sigma\sigma\sigma}^2\,d\sigma$ which is less than
$\int\limits_{\widetilde{\TTT}_{x_0,\tt}}
\tt^2\widetilde{k}_{\sigma\sigma\sigma}^2\varphi^2\,d\sigma$.\\
Hence, we conclude (in the case $\TTT_\infty$ is a halfline the
3--point contribution is not present)
\begin{align*}
\frac{d\,}{d\tt}\int\limits_{\widetilde{\TTT}_{x_0,\tt}}
\left(\widetilde{k}^2+\tt\widetilde{k}_\sigma^2+\tt^2\widetilde{k}_{\sigma\sigma}^2/2\right)
\varphi^2\,d\sigma
\leq&\,+C\int\limits_{\widetilde{\TTT}_{x_0,\tt}}
\left(\widetilde{k}^2+\tt\widetilde{k}_{\sigma}^2
+\tt^2\widetilde{k}_{\sigma\sigma}^2/2\right)\varphi^2\,d\sigma + C\\
&\,-\sum_{i=1}^32\widetilde{k}\widetilde{k}_\sigma+
2\tt\widetilde{k}_{\sigma}\widetilde{k}_{\sigma\sigma}+
\tt^2\widetilde{k}_{\sigma\sigma}\widetilde{k}_{\sigma\sigma\sigma}
\,\biggr\vert_{\text{{      at the 3--point}}}\\
&\,-\sum_{i=1}^32\widetilde{k}^2\widetilde{\lambda} +
\tt\widetilde{k}_{\sigma}^2\widetilde{\lambda}+
\tt^2\widetilde{k}_{\sigma\sigma}^2\widetilde{\lambda}/2
\,\biggr\vert_{\text{{      at the 3--point}}}
\end{align*}
and dealing with the 3--point terms as in Section~\ref{kestimates} we
finally get
$$
\int\limits_{\widetilde{\TTT}_{x_0,\tt}}
\left(\widetilde{k}^2+\tt\widetilde{k}_\sigma^2+\tt^2\widetilde{k}_{\sigma\sigma}^2/2\right)
\varphi^2\,d\sigma\leq C
$$
for every $\tt\in[0,\delta]$, with uniform constants $\delta$ and $C$.\\
Reparametrizing the flow in time, in such a way that
$\tt_{j_l}$ converges to $\delta$, we have
$$
\int\limits_{\widetilde{\TTT}_{x_0,\tt_{j_l}}\cap B_R}
\widetilde{k}_\sigma^2\,d\sigma\leq C/\delta
$$
for every $l\in\NN$.\\
Then, since $\Vert\widetilde{k}_\sigma\Vert_{L^2(B_R)}$ and
$\widetilde{k}$ are uniformly  
bounded, we can finally extract a subsequence of triods converging in
the $C^2\loc$ topology to $\TTT_\infty$.
\end{proof}

By means of this proposition we can exclude the first case of 
Lemma~\ref{casi}.

\begin{prop}\label{tripsing} {\em Type I} singularities such that for an
  index $i\in\{1, 2, 3\}$ there exists a sequence of times $t_j\to T$ satisfying
\begin{equation*}
[k^i(O,t_j)]^2\geq \frac{C}{T-t_j}
\end{equation*}
for some positive constant $C$, are not present.
\end{prop}
\begin{proof}
We have seen in Lemma~\ref{Obehav} that $\lim_{t\to T}F(O,t)=\widehat{O}$.\\
Considering the triods $\widetilde{\TTT}_{\widehat{O},\tt_j}$,  where
$\tt_j=\tt(t_j)$, by Proposition~\ref{resclimit} we can extract a 
subsequence (not relabelled) $C^2$ converging in the ball of radius
$\sqrt{2C}$ to a limit set with zero curvature, moreover this limit is
not empty because the 3--point  $\widetilde{F}_{\widehat{O}}(O,\tt_j)$
belongs to such ball for every $j\in\NN$, again by
Lemma~\ref{Obehav}.\\
We have now a contradiction because at these points
$$
[\widetilde{k}^i(O,\tt_j)]^2
=2(T-t_j)[{k}^i(O,t_j)]^2\geq{2(T-t_j)}\frac{C}{T-t_j}=2C>0\,,
$$
hence, as the convergence is in the $C^2$ topology, 
the curvature of the limit cannot be zero.
\end{proof}

Following Ilmanen, we deal now 
with the second case of Lemma~\ref{casi}.

\begin{prop}\label{nosing0} There are no {\em Type I} singularities 
such that 
$$
\max_{\TTT_t}k^2=\frac{\varepsilon(t)}{T-t}
$$
and $\varepsilon:[0,T)\to\R$ goes to zero as $t\to T$.
\end{prop}
\begin{proof}
By means of Lemma~\ref{preterzo} we know that 
there exists an index $i\in\{1,2,3\}$ and a sequence of times
$t_j\nearrow T$ such that
$$
k^2(p,t)\leq [k^i(O,t_j)]^2\, \text{ for every $t\leq t_j$,
  $p\in\TTT$ and } [k^i(O,t_j)]^2\nearrow+\infty\,.
$$
By Lemma~\ref{Obehav}, 
the limit $\lim_{t\to T}F(O,t)=\widehat{O}$ exists and
repeating the computation in its proof we get
\begin{equation*}
\vert F(O,t)-\widehat{O}\,\vert\leq 2\sqrt{C(T-t)\varepsilon(t)}
\end{equation*}
for every $t\in[0,T)$.\\
We can suppose that $\widehat{O}=0\in\R^2$ and we consider the
sequence of positive values $\alpha_j=\max_{\TTT_{t_j}}\vert k\vert
\nearrow+\infty$, then we rescale the
triods $\TTT_t$ as follows, let 
$$
\TTT_\zz^{j}=\alpha_j({\TTT_{t_j+\zz/\alpha_j^2}})\qquad\text{ {
    for} $\zz\in[-t_j\alpha_j^2,(T-t_j)\alpha_j^2)$.}
$$
We see that $-t_j\alpha_j^2\to-\infty$ and, possibly passing to a
subsequence, we can assume that $(T-t_j)\alpha_j^2\searrow0$. Then,
notice that, for every $j\in\NN$ large enough, 
the rescaled triods $\TTT_\zz^{j}$
still move by curvature with end points $\alpha_j P^i$ for
$\zz\in[-1,0]$. Moreover if $O^{j}(\zz)$ is the 3--point of 
$\TTT_\zz^{j}$, we have
\begin{align}\label{boundonO}
\vert O^{j}(\zz)\vert=&\, 
\alpha_j\vert F(O,  t_j+\zz/\alpha_j^2)\vert\\
\leq&\,\alpha_j{2\sqrt{C(T-t_j-\zz/\alpha_j^2)\varepsilon(t_j+\zz/\alpha_j^2)}}\nonumber\\
=&\,{2\sqrt{C[(T-t_j)\alpha_j^2-\zz]\varepsilon(t_j+\zz/\alpha_j^2)}}\nonumber
\end{align}
which goes to zero when $j\to\infty$ uniformly for $\zz\in[-1,0]$, 
since $[(T-t_j)\alpha_j^2-\zz]\to
-\zz$ and $\varepsilon(t_j+\zz/\alpha_j^2)\to0$.\\
We set $\varepsilon_j=\varepsilon(t_j)\geq0$ and 
we denote with
$$
\rho_j(x,\zz)=\frac{e^{-\frac{\vert
    x\vert^2}{4(\varepsilon_j-\zz)}}}{\sqrt{4\pi(\varepsilon_j- \zz)}}\,,
$$
the Huisken's backward heat kernel relative to
$(0,\varepsilon_j)\in\R^2\times\R$.\\
Computing as in~\cite[Lemma~8]{ilman3} and~\cite[Section~3]{ilman4} we
show that a subsequence of the flows $\TTT_\zz^{j}$ converges 
to a curvature flow $\TTT_\zz^\infty$ on $[-1,0]$ which is homothetic.
\begin{align*}
\int_{-1}^{0}\int_{\TTT_\zz^{j}}
\left\vert\,\underline{k}+\frac{x^{\perp}}{2(\varepsilon_j-\zz)}\right\vert^2
&\,\rho_j(x,\zz)\,d\sigma\,d\zz\\ 
=&\,\int_{-1}^{0}\int_{\TTT_{t_j+\zz/\alpha_j^2}}\alpha_j
\left\vert\,\frac{\underline{k}}{\alpha_j}+\frac{\alpha_j
    y^{\perp}}{2(\varepsilon_j-\zz)}\right\vert^2
\frac{e^{-\frac{\alpha_j^2\vert
    y\vert^2}{4(\varepsilon_j-\zz)}}}{\sqrt{4\pi(\varepsilon_j- \zz)}}
\,ds\,d\zz\\
=&\,\int_{t_j-1/\alpha_j^2}^{t_j} \alpha_j^2\int_{\TTT_{t}}
\alpha_j\left\vert\,\frac{\underline{k}}{\alpha_j}
+\frac{\alpha_j y^{\perp}}{2(T-t)\alpha_j^2}\right\vert^2
\frac{e^{-\frac{\alpha_j^2\vert
    y\vert^2}{4(T-t)\alpha_j^2}}}{\sqrt{4\pi(T-t)\alpha_j^2}}
\,ds\,dt\\
=&\,\int_{t_j-1/\alpha_j^2}^{t_j}\int_{\TTT_{t}}
\left\vert\,\underline{k}+\frac{y^{\perp}}{2(T-t)}\right\vert^2
\rho_0(y,t)\,ds\,dt
\end{align*}
where $\rho_0$ is the backward heat kernel relative to
$(0,T)\in\R^2\times\R$.\\
By the integrated monotonicity formula~\eqref{eqmonfor-int}, this last
term is equal to
$$
\int_{{\TTT_{t_j-1/\alpha_j^2}}} \rho_0(y,t_j-1/\alpha_j^2)\,ds-
\int_{{\TTT_{t_j}}} \rho_0(y,t_j)\,ds
+\sum_{i=1}^3\int_{t_j-1/\alpha_j^2}^{t_j}
\frac{\left\langle P^i\,\vert\,\tau^i(1,t)\right\rangle}{2(T-t)}
\rho_0(P^i,t)\,dt
$$
and this last expression goes to zero when $j\to\infty$, 
by Lemma~\ref{valori0} and~\ref{stimadib}.\\
Thus, for almost
every $\zz\in[-1,0]$ we have that
\begin{equation}\label{carlo100}
\lim_{j\to\infty} \int_{\TTT_\zz^{j}}
\left\vert\,\underline{k}+\frac{x^{\perp}}{2(\varepsilon_j-\zz)}\right\vert^2
\rho_j(x,\zz)\,d\sigma=0\,.
\end{equation}
Now, all the flows $\TTT_\zz^{j}$ have uniformly bounded curvature for 
$\zz\in[-1,0]$ since by construction,
$$
\max_{\TTT_\zz^{j}} k^2=\alpha_j^{-2}\max_{\TTT_{t_j+\zz\alpha_j^2}}
k^2\leq \alpha_j^{-2}\max_{\TTT_{t_j}} k^2 =\alpha_j^{-2}[k^i(O,t_j)]^2=1
$$
hence, repeating the local length estimate of Lemma~\ref{rescestim} and the 
convergence argument in the proof of
Proposition~\ref{resclimit} (interior estimates of Ecker and
Huisken plus the special treatment of the 3--point terms, see
Remark~\ref{deltabuono2}), we can extract a
subsequence of the flows converging in the $C^2\loc$ topology 
to a curvature flow of triods $\TTT_\zz^\infty$ 
in $[-,1,0]$ which, by limit~\eqref{carlo100} must 
satisfy
$$
\int_{\TTT_\zz^{\infty}}
\left\vert\,\underline{k}- x^{\perp}/2\zz\right\vert^2
\rho_j(x,\zz)\,d\sigma=0
$$
for almost every $\zz\in[-1,0]$. By Lemma~\ref{lemomt}, it follows
that {\em for every} ${\zz\in[-1,0]}$, every triod
$\TTT_\zz^\infty$ is composed of three halflines by the origin of
$\R^2$ and it has zero curvature.\\
Looking in particular at $\TTT_0^\infty$ which is the limit in the 
$C^2$ topology of the triods $\alpha_j\TTT_{t_j}$, we finally 
have a contradiction since these latter have at least one squared
curvature  equal to one at their 3--points (which converge to the
origin) by the initial choice of the sequence $t_j\nearrow T$.
\end{proof}

The rest of the section is concerned with the last case of
Lemma~\ref{casi},  so, from now on,
we suppose that the maximum of the curvature is
not taken at the 3--point and that $\lim_{t\to T}
(T-t)k^2(O,t)=0$, otherwise we are in the first case.
Since in this case we have no bound on the tangential velocity at
a maximum point of the curvature, we modify our function $F$ to keep
it under control.\\
Let $\theta^i:[0,1]\times[0,T)\to[0,1]$ be the smooth solutions of the
following ODE's,
$$
\frac{\partial \theta^i(x,t)}{\partial
  t}=\frac{\lambda^i(0,t)(1-x)-\lambda^i(\theta^i(x,t),t)}{\vert
  \gamma^i_x(\theta^i(x,t),t)\vert}
$$
with the initial conditions $\theta^i(x,0)=x$. These functions are
well defined on the interval $[0,T)$ since the curves 
$\gamma^i$ and all their derivatives are smooth and 
$\gamma^i_x\not=0$, moreover 
$\theta^i(0,t)=0$, $\theta^i(1,t)=1$ are barrier solutions, as
$\lambda^i(1,t)=0$.\\
Defining naturally, via $\theta^i$, a function 
$\theta:\TTT\times[0,T)\to\TTT$,  
the function $G(p,t)=F(\theta(p,t),t)$ satisfies
\begin{align*}
\frac{\partial G(p,t)}{\partial t} = 
\underline{k}+\underline{\lambda}+F_p
\frac{\partial  \theta}{\partial  t}
=&\,\underline{k}+\lambda\tau+\frac{F_p}
{\vert  F_p\vert}(\lambda(O,t)(1-x(p))-\lambda)\\
=&\,\underline{k}(G(p,t))+{\lambda}(O,t)(1-x(p))\tau(G(p,t))
\end{align*}
where $x(p)$ denotes the value in the interval $[0,1]$ such that for a
certain index $i\in\{1,2,3\}$, $F(p,t)=\gamma^i(x(p),t)$.\\
Notice that, by construction, $G(O,t)=F(O,t)$ and $G(P^i,t)=F(P^i,t)$.\\ 

{\em From now on in all the rest of this section 
we will refer all the quantities $\tau$, $\nu$, $k$, $\lambda$,
etc... to the new parametrization of the flow $G$, that is, for instance 
$k(p,t)$ is the curvature of the triod $\TTT_t$ at the point
$G(p,t)$ and so on.}\\

Denoting with $\underline{w}(p,t)$ the new velocity of the point
$G(p,t)\in\TTT_t$ and with $\underline{\eta}(p,t)=\eta(p,t)\tau(p,t)$
its tangential part, that is,
$\underline{w}=\underline{k}+\underline{\eta}$, we have
$$
\underline{\eta}(p,t)=\eta(p,t)\tau(p,t)={\lambda}(O,t)(1-x(p))\tau(p,t)\,,
$$
then clearly $\eta^i(O,t)=\lambda^i(O,t)$ and $\eta(P^i,t)=0$.\\
By assumption~\eqref{stimadallalto} and the relations between $k$ and
$\lambda$ at the 3--point, the new velocity $w(p,t)=
\underline{k}(p,t)+{\lambda}(O,t)(1-x(p))\tau(p,t)$ 
is uniformly bounded by $C/\sqrt{T-t}$.

The evolution equation for the curvature has to be modified in
$$
{\partial_t}k=k_{ss}+\eta
k_{s}+k^3
$$
and the commutation law as follows,
$$
\partial_t\partial_s
=\partial_s\partial_t 
+ (k^2 -\eta_s)\partial_s\,.
$$
Rescaling $G$, like we did before for $F$, 
we define the map $\widetilde{G}_{x_0}:\TTT\times
[-1/2\log{T},+\infty)\to\R^2$ as 
$$
\widetilde{G}_{x_0}(p,\tt)=\frac{G(p,t)-x_0}{\sqrt{2(T-t)}}\qquad
\tt(t)=-\frac{1}{2}\log{(T-t)}\,.
$$
Notice that the rescaled triods are not changed
$$
\widetilde{\TTT}_{x_0,\tt}=\frac{\TTT_t-x_0}{\sqrt{2(T-t)}}
$$
but now they evolve according to the equation 
$$
\frac{\partial\,}{\partial
  \tt}\widetilde{G}_{x_0}(p,\tt)=\widetilde{\underline{w}}(p,\tt)+\widetilde{G}_{x_0}(p,\tt)
$$
where $\widetilde{\underline{w}}(p,\tt)$ is the rescaled velocity 
given by
$$
\widetilde{\underline{w}}(p,\tt)=\frac{\underline{w}(p,t(\tt))}{\sqrt{2(T-t(\tt))}}=
\widetilde{\underline{k}}+\widetilde{\lambda}(O,t(\tt))(1-x(p))\tau
=\widetilde{k}\nu+\widetilde{\eta}\tau\,.
$$
The rescaled curvature evolves according to the following equation,
\begin{equation}\label{evolriscforg}
{\partial_\tt}\widetilde{k}=\widetilde{k}_{\sigma\sigma}+\widetilde{\eta}
\widetilde{k}_\sigma+\widetilde{k}^3-\widetilde{k}
\end{equation}
and the commutation law is
$$
{\partial_\tt}{\partial_\sigma}={\partial_\sigma}{\partial_\tt} + (\widetilde{k}^2
-\widetilde{\eta}_\sigma-1){\partial_\sigma}\,.
$$

We call $p_t\not=O$ a point in $\TTT$ where the
curvature achieves its maximum, so by hypothesis
$$
\frac{1}{2(T-t)}\leq k^2(p_t,t)\leq\frac{C}{2(T-t)}\,.
$$

We define the map $\,\,\, \widehat{ }: \TTT\to\R^2$ as follows
$$
\widehat{p}=\lim_{t\to T} G(p,t)
$$
for every $p\in\TTT$.\\
Such limit exists for every $p\in\TTT$ since, for every $0\leq
t_1\leq t_2<T$, we have
$$
\vert G(p,t_1)-G(p,t_2)\vert=\left\vert\int_{t_1}^{t_2}\frac{\partial G}{\partial
    t}\,dt\right\vert\leq \int_{t_1}^{t_2}\vert
\underline{w}(p,t)\vert\,dt\leq\int_{t_1}^{t_2}\sqrt{\frac{C}{T-t}}\,dt
\leq 2\sqrt{C(T-t_1)}
$$
and then
\begin{equation}\label{distanza-lim}
\vert G(p,t)-\widehat{p}\,\vert\leq 2\sqrt{C(T-t)}
\end{equation}
for every $t$. Notice that the definition of $\widehat{O}$ coincides
with the one in Lemma~\ref{Obehav}, since $G(O,t)=F(O,t)$.

\begin{lemma}\label{puntopalla}\ \\
\begin{enumerate}
\item The map $\,\,\, \widehat{ }: \TTT\to\R^2$ is continuous.
\item For every $p\in\TTT$ all the triods
  $\widetilde{\TTT}_{\widehat{p},\tt}$ intersect the closed ball
  $B_{\sqrt{2C}}$ centered at the origin of $\R^2$, indeed,
  $\widetilde{G}_{\widehat{p}}(p,\tt)$ belongs to such ball for every
  $\tt\in[-1/2\log{T},+\infty)$.
\item The point $\widehat O$ is different from $P^1$, $P^2$, $P^3$
  (which coincide respectively with $\widehat{P}^1$, $\widehat{P}^2$,
  $\widehat{P}^3$).
\end{enumerate}
\end{lemma}
\begin{proof}\ \\
\begin{enumerate} 
\item Since all the maps $G(\cdot , t)$ are continuous and
inequality~\eqref{distanza-lim} says that they converge uniformly as
$t\to T$, the map
$\,\,\, \widehat{ }: \TTT\to\R^2$ is also continuous.
\item By inequality~\eqref{distanza-lim}, we have
$$
\vert\widetilde{G}_{\widehat{p}}(p,\tt)\vert=\frac{\vert
G(p,t(\tt))-\widehat{p}\vert}{\sqrt{2(T-t(\tt))}}
\leq \frac{2\sqrt{C(T-t(\tt))}}{\sqrt{2(T-t(\tt))}}=\sqrt{2C}
$$
and the second statement is proved.
\item The point $\widehat{O}$ cannot coincides with any of the end
  points $P^i$, let us say $P^1$, otherwise, rescaling the triods
  around such end point, by Proposition~\ref{resclimit} we can find
  a subsequence of the rescaled triods converging in the $C^2\loc$
  topology to a halfline from the origin of $\R^2$, but
  Lemma~\ref{Obehav} says that since $\widehat{O}=P^1$ the 3--point
  cannot disappear in the limit, which is a contradiction.
\end{enumerate}
\end{proof}

Following Stone~\cite{stone1}, we define the 
function $\Theta(p,t)$ as
$$
\Theta(p,t_j)=\int_{\TTT_{t_j}}\rho_{\widehat{p}}\,ds=
\frac{1}{\sqrt{2\pi}}\int_{\widetilde{\TTT}_{\widehat{p},\tt_j}}\widetilde{\rho}\,d\sigma
$$
for every $p\in\TTT$ and the {\em limiting heat
  density} as the limit for $t\to T$ of
$$
\widehat{\Theta}(p)=\lim_{t\to T}\Theta(p,t)\,,
$$
if it exists.\\
Notice that since $p\mapsto\widehat{p}\,$ is continuous, all the maps 
$p\mapsto\Theta(p,t)$ are also continuous, for every $t\in[0,T)$

\begin{prop}\label{valori}
The limit $\widehat{\Theta}(p)$ exists and it is finite for   every
$p\in\TTT$, moreover, $\widehat{\Theta}(p)$ can take only the values
  $1/2$, $1$, $3/2$.
\end{prop}
\begin{proof}
{As in Lemma~\ref{valori0}, we define the function $b:\TTT\times[0,T)\to\R$ as follows,
\begin{equation}\label{defb}
b(p,t)=\int_t^T\sum_{i=1}^3
\left\langle\,\frac{P^i-\widehat{p}}{2(T-\xi)}\,\biggl\vert\,\tau^i(1,\xi)\right\rangle
\rho_{\widehat{p}}(P^i,\xi)\,d\xi
\end{equation}
then Lemma~\ref{stimadib} says that $b$ is bounded and for every
$p\in\TTT$ we have $\lim_{t\to T}b(p,t)=0$. The monotonicity
formula~\eqref{eqmonfor} can be rewritten as
$$
\frac{d\,}{dt}\bigl(\Theta(p,t)+b(p,t)\bigr)=-\int_{{\TTT_t}}
\left\vert\,\underline{k}+\frac{(x-\widehat{p}\,)^{\perp}}{2(T-t)}\right\vert^2
\rho_{\widehat{p}}\,ds\leq 0\,,
$$
hence, being non increasing and bounded from below, the
functions $\bigl(\Theta(\cdot , t)+b(\cdot , t)\bigr)$ pointwise converge on all
$\TTT$ when $t\to T$. Since we have seen that $b(\cdot , t)$ also
pointwise converge to zero everywhere, the limit $\widehat{\Theta}(p)$
exists for every $p$.}

We consider now, by means of Proposition~\ref{resclimit}, 
a sequence of times $t_j\nearrow T$ such that the
rescaled triods ${\widetilde{\TTT}_{\widehat{p},\tt_j}}$ 
converge in the $C^2\loc$ topology to some zero curvature
set $\TTT_\infty$.\\
Taking into account that $\widetilde{k}$ is uniformly bounded, we
compute
\begin{align*}
\frac{d\,}{d\tt}\int_{\widetilde{\TTT}_{\widehat{p},\tt}}e^{-\vert x
\vert}\,d\sigma
=&\,\int_{\widetilde{\TTT}_{\widehat{p},\tt}}e^{-\vert x
\vert}\left(1-\widetilde{k}^2 -\frac{\langle
x\,\vert\, \widetilde{\underline{k}}+x\rangle}{\vert x\vert}\right)\,d\sigma\\
\leq &\,\int_{\widetilde{\TTT}_{\widehat{p},\tt}}e^{-\vert x
\vert}\left(C-\vert x\vert\right)\,d\sigma\\
\leq &\,C\int\limits_{\widetilde{\TTT}_{\widehat{p},\tt}\cap B_{C+1}}e^{-\vert x
\vert}\,d\sigma
-\int\limits_{\widetilde{\TTT}_{\widehat{p},\tt}\setminus B_{C+1}}e^{-\vert x
\vert}\,d\sigma\\
\leq &\,(C+1)\int\limits_{\widetilde{\TTT}_{\widehat{p},\tt}\cap B_{C+1}}e^{-\vert x
\vert}\,d\sigma
-\int\limits_{\widetilde{\TTT}_{\widehat{p},\tt}}e^{-\vert x
\vert}\,d\sigma\,.
\end{align*}
This means that if
$$
\frac{d\,}{d\tt}\int_{\widetilde{\TTT}_{\widehat{p},\tt}}e^{-\vert x
\vert}\,d\sigma\geq 0
$$
then
$$
\int\limits_{\widetilde{\TTT}_{\widehat{p},\tt}}e^{-\vert x
\vert}\,d\sigma
\leq(C+1)\int\limits_{\widetilde{\TTT}_{\widehat{p},\tt}\cap B_{C+1}}e^{-\vert x
\vert}\,d\sigma
$$
which clearly implies that for every $\tt\in[-1/2\log(T-t),+\infty)$ 
\begin{align*}
\int\limits_{\widetilde{\TTT}_{\widehat{p},\tt}}e^{-\vert x
\vert}\,d\sigma
\leq &\,(C+1)\int\limits_{\widetilde{\TTT}_{\widehat{p},\tt}\cap B_{C+1}}e^{-\vert x
\vert}\,d\sigma\\
\leq &\,(C+1)\int\limits_{\widetilde{\TTT}_{\widehat{p},\tt}\cap
  B_{C+1}}
e^{\frac{C^2-\vert x\vert^2}{2}}\,d\sigma\\
= &\,(C+1)e^{C^2/2}\int\limits_{\widetilde{\TTT}_{\widehat{p},\tt}\cap B_{C+1}}
e^{-\frac{\vert x\vert^2}{2}}\,d\sigma\\
= &\,\frac{(C+1)e^{C^2/2}}{\sqrt{2\pi}}
\Theta(p,t(\tt))\leq D
\end{align*}
where the constant $D$ does not depend on $\tt$, since we have seen that
$\Theta(p,t)$ converges as $t\to T$. 

Now, we subdivide $\widetilde{\TTT}_{\widehat{p},\tt_j}$ into annular  pieces
$\widetilde{\TTT}_{\widehat{p},\tt_j}^{n}$ for $n\in\NN$, as follows:
$$
{\widetilde{\TTT}}^{\,0}_{\widehat{p},\tt_j}=\widetilde{\TTT}_{\widehat{p},\tt_j}\cap
B_1\qquad
{\widetilde{\TTT}}^{n}_{\widehat{p},\tt_j}=\{x\in\widetilde{\TTT}_{\widehat{p},\tt_j}\,\vert\,
2^{n-1}\leq \vert x\vert<2^n\} \text{ \, { \, for every} $n\geq 1$.}
$$
By the computations above $\int_{\widetilde{\TTT}_{\widehat{p},\tt_j}}
e^{-\vert x\vert}\,d\sigma$ is bounded independently of $j\in\NN$,
then trivially
$$
\HH^1({\widetilde{\TTT}_{\widehat{p},\tt}^{n}})\leq e^{(2^n)}
\int_{\widetilde{\TTT}_{\widehat{p},\tt}^{n}}
e^{-\vert x \vert}\,d\sigma\leq Ae^{(2^n)}\,,
$$
hence, for every $n\geq 1$ we have
$$
\int_{\widetilde{\TTT}_{\widehat{p},\tt}^{n}}e^{-\frac{\vert x
\vert^2}{2}}\,d\sigma\leq e^{-\frac{1}{2}(2^{n-1})^2} A e^{(2^n)} =
Ae^{(2^n-2^{2n-3})}\,.
$$
Now, for every $\varepsilon>0$ we can find a number $n_0\in\NN$ such that
$\sum_{n=n_0}^\infty Ae^{(2^n-2^{2n-3})}\leq\varepsilon$, that is, 
if $R\geq 2^{n_0-1}$,
$$
\int\limits_{\widetilde{\TTT}_{\widehat{p},\tt}\setminus B_R}e^{-\frac{\vert x
\vert^2}{2}}\,d\sigma\leq\varepsilon
$$
for every $\tt\in[-1/2\log(T-t),+\infty)$.\\
Since we have
$$
\Theta(p,t_j)=\int_{\TTT_{t_j}}\rho_{\widehat{p}}\,ds=
\frac{1}{\sqrt{2\pi}}\int_{\widetilde{\TTT}_{\widehat{p},\tt_j}}\widetilde{\rho}\,d\sigma
$$
we get
$$
\left\vert \Theta(p,t_j) -
  \frac{1}{\sqrt{2\pi}}\int\limits_{\widetilde{\TTT}_{\widehat{p},\tt_j}
\cap B_{R}}\widetilde{\rho}\,d\sigma\right\vert
\leq\frac{1}{\sqrt{2\pi}}\int\limits_{\widetilde{\TTT}_{\widehat{p},\tt_j}\setminus
B_R}\widetilde{\rho}\,d\sigma\leq\frac{\varepsilon}{\sqrt{2\pi}}\,.
$$
Since $\Theta(p,t_j)$ converges to $\widehat{\Theta}(p)$ 
and the measures $\HH^1\res{\widetilde{\TTT}_{\widehat{p},\tt_j}}$ go to
$\HH^1\res{\TTT}_\infty$, by Proposition~\ref{resclimit}, we obtain
$$
\left\vert \widehat{\Theta}(p) - \frac{1}{\sqrt{2\pi}}
\int\limits_{{\TTT}_\infty\cap B_{R}}\widetilde{\rho}\,d\sigma\right\vert
\leq\frac{\varepsilon}{\sqrt{2\pi}}
$$
and sending $R$ to $+\infty$,
$$
\left\vert \widehat{\Theta}(p) - \frac{1}{\sqrt{2\pi}}
\int_{{\TTT}_\infty}\widetilde{\rho}\,d\sigma\right\vert
\leq\frac{\varepsilon}{\sqrt{2\pi}}\,,
$$
hence, by the arbitrariness of $\varepsilon$,
$$
\widehat{\Theta}(p) = \frac{1}{\sqrt{2\pi}}
\int_{{\TTT}_\infty}\widetilde{\rho}\,d\sigma\,.
$$ 
As ${\TTT}_\infty$ is one of the sets of 
Proposition~\ref{resclimit}, the claim follows from the fact that
$\int_0^{+\infty}e^{-y^2/2}\,dy=\sqrt{\pi/2}$.
\end{proof}

We consider the following three parts of $\TTT$,
\begin{align*}
&W^- =\left\{p\in\TTT\,\biggl\vert\,\widehat{\Theta}(p)=1/2\right\}\\
&W =\left\{p\in\TTT\,\biggl\vert\,\widehat{\Theta}(p)=1\right\}\\
&W^+ =\left\{p\in\TTT\,\biggl\vert\,\widehat{\Theta}(p)=3/2\right\}\,.
\end{align*}
Since by Lemma~\ref{puntopalla} 
the point $\widehat{O}$  is distinct from $P^1$, $P^2$ and $P^3$ 
(which respectively coincide with  $\widehat{P}^1$, $\widehat{P}^2$
and $\widehat{P}^3$),  
if $\, \widehat{p}$ is different from all the $P^i$ and
$\widehat{O}$, taking into account Proposition~\ref{resclimit}, every
$C^2\loc$ converging subsequence of the rescaled triods
$\widetilde{\TTT}_{\widehat{p},\tt}$ must go to a straight line, 
as the 3--point and the end points go to infinity (recall
Lemma~\ref{Obehav}). Hence $\widehat{\Theta}(p)=1$.\\
If $\, \widehat{p}=P^1$, for instance, then
$\widehat{\Theta}(p)=\widehat{\Theta}(P^1)$ which is equal to $1/2$ since the
3--point of the rescaled triods goes to infinity and the end 
point $P^1$ is fixed, hence the limit is a halfline from the origin of
$\R^2$.\\
Finally, with the same argument (now the three end points 
go to infinity) if $\, \widehat{p}=\widehat{O}$ then
$\widehat{\Theta}(p)=3/2$.\\
We get then the following equalities
\begin{align*}
&W^- =\biggl\{p\in\TTT\,\biggl\vert\,\widehat{p}\in\{P^i\}\Bigr\}\\
&W =\biggl\{p\in\TTT\,\biggl\vert\,\widehat{p}\notin\{P^i, \widehat{O}\}\Bigr\}\\
&W^+ =\left\{p\in\TTT\,\biggl\vert\,\widehat{p}=\widehat{O}\right\}
\end{align*}
moreover, notice that they are non empty with $W^-$, $W^+$ closed
and $W$ open in $\TTT$, as the map $p\mapsto\widehat{p}$ is
continuous.

Fixed $x_0\in\R^2$,  if $\limnf_{t\to T} \frac{\vert\widehat{p}_t-
  x_0\vert^2}{T-t}<+\infty$ (we recall that $p_t\in\TTT$ is a 
maximum point of $k^2$ at time $t\in[0,T)$) we can find a 
sequence of points $p_j=p_{t_j}\in\TTT$ such that $\vert
  \widehat{p}_j-x_0\vert^2\leq 2D(T-t_j)$, for some constant $D$.\\
Rescaling around $x_0$ like in the proof of
Proposition~\ref{tripsing},  we see that the limit triod is not empty
because the point $\widetilde{G}_{x_0}(p_j,\tt_j)$ belongs
to the closed ball of radius $\sqrt{2C}+\sqrt{D}$, indeed
\begin{align*}
\vert\widetilde{G}_{x_0}(p_j,\tt_j)\vert=\frac{\vert
G(p_j,t_j)-x_0\vert}{\sqrt{2(T-t_j)}}
\leq&\, \frac{\vert G(p_j,t_j)-\widehat{p}_j\vert}{\sqrt{2(T-t_j)}}+
\frac{\vert \widehat{p}_j-x_0\vert}{\sqrt{2(T-t_j)}}\\
\leq&\, \frac{2\sqrt{C(T-t_j)}}{\sqrt{2(T-t_j)}}+
\frac{\sqrt{2D(T-t_j)}}{\sqrt{2(T-t_j)}}\\
=&\,\sqrt{2C}+\sqrt{D}\,,
\end{align*}
by Lemma~\ref{puntopalla}. Moreover, 
$$
\widetilde{k}^2(p_j,\tt_j)^2=2(T-t_j){k}^2(p_j,t_j)^2
=2(T-t_j)\max_{\TTT_{t_j}}k^2
\geq{2(T-t_j)}\frac{1}{2(T-t_j)}=1\,,
$$
by our initial assumption, hence, we can exclude this situation as we
did in Proposition~\ref{tripsing}.\\
Thus, we can assume that for every $x_0\in\R^2$
\begin{equation}\label{assump999}
\lim_{t\to T} \frac{\vert \widehat{p}_t- x_0\vert^2}{T-t}=+\infty\,,
\end{equation}
in particular for $x_0\in\{P^i, \widehat{O}\}$. This implies that 
the points $p_t$ belong definitely to the set
$W$, that is, $\widehat{\Theta}(p_t)=1$.

We consider then a sequence of times $t_j\nearrow T$, such that the
sequence of points $p_j=p_{t_j}\in\TTT$ converges 
to some $p\in\TTT$. By construction,  
$\max_{\TTT_{t_j}}k^2=k^2(p_j,t_j)$, moreover we can suppose that 
all the points $p_j$ belong to the same curve of $\TTT$ and 
to the set $W$, that is, $\widehat{\Theta}(p_j)$ is constantly
equal to $1$. By the uniform continuity of  the map
$G$, we have $G(p_j,t_j)\to\widehat{p}\,$ and 
$\,\widehat{p}_j\to\widehat{p}\,$.\\
Assumption~\eqref{assump999}  says also that for every radius $R>0$ 
there exists an index $j_0\in\NN$ such
that all the sets $\widetilde{\TTT}_{\widehat{p}_j,\tt_j}\cap B_{3R}$ do not
contain neither the 3--point nor the three end points, for
every $j\geq j_0$. Indeed,
\begin{align*}
\vert\widetilde{G}_{\widehat{p}_j}(O,\tt_j)\vert=\frac{\vert
G(O,t_j)-\widehat{p}_j\vert}{\sqrt{2(T-t_j)}}
\geq&\, \frac{\vert \widehat{p}_j-\widehat{O}\vert}{\sqrt{2(T-t_j)}}
-\frac{\vert G(O,t_j)-\widehat{O}\vert}{\sqrt{2(T-t_j)}}\\
\geq&\, \frac{\vert \widehat{p}_j-\widehat{O}\vert}{\sqrt{2(T-t_j)}}
- \frac{2\sqrt{C(T-t_j)}}{\sqrt{2(T-t_j)}}\\
\geq&\, \frac{\vert \widehat{p}_j-\widehat{O}\vert}{\sqrt{2(T-t_j)}}
- \sqrt{2C}
\end{align*}
which goes to $+\infty$ as $j\to\infty$. The argument is the same for
the three end points $P^i$.

\begin{lemma}\label{intornobuono} For every $\varepsilon>0$ there
  exists an open neighborhood $U_\varepsilon$ of $W^-$ in $\TTT$ and a time
  $t_\varepsilon$ such that $\Theta(p,t)\leq 1+\varepsilon$ for every
  $p\in U_\varepsilon$ and $t\in[t_\varepsilon,T)$.
\end{lemma}
\begin{proof}
Suppose, by contradiction, that $t_j\nearrow T$, $p_j=p_{t_j}\to p\in
W^-$ and $\Theta(p_j,t_j)>1+\varepsilon$ for every $j\in\NN$. 
Then, we can also assume that $\,\widehat{p}=P^1$ 
hence, $d_j^2, d_j^3>C$ for a positive constant
$C$ and $d_j^1\to0$, where $d_j^i=\vert P^i-\widehat{p}_j\vert$.\\
Since on $W^-$ the functions $b(\cdot , t)$ are locally constant,
hence continuous, and by Lemma~\ref{stimadib} they 
converge uniformly to zero in $W^-$ which 
is compact, by Dini's Theorem, the functions
$\Theta(\cdot,t)$, restricted to $W^-$, converge uniformly to
$1/2$. Hence, fixing $\varepsilon>0$, we can find a time
$t_\varepsilon$ and an open neighborhood $U_\varepsilon$ of $W^-$ 
such that $\Theta(p,t_\varepsilon)\leq 1/2+\varepsilon/2$ for every
  $p\in U_\varepsilon$ and $\Theta(p,t)\leq 1/2+\varepsilon$ for every
  $p\in W^-$ and $t\geq t_\varepsilon$. 
Moreover, we can choose $t_\varepsilon$ such that
$$
\frac{1}{\sqrt{2\pi}}\int\limits_{C/\sqrt{2(T-t_\varepsilon)}}^{+\infty}e^{-y^2/2}\,dy\leq
\varepsilon/4\,.
$$
If now $p_j\in U_\varepsilon\setminus W^-$ and $t_\varepsilon\leq t_j<T$ we
estimate $\Theta(p_j,t_j)$ with equation~\eqref{eqmonfor-int} as follows,
\begin{align*}
\Theta(p_j,t_j)= &\,\,\Theta(p_j,t_\varepsilon) + \bigl(\Theta(p_j,t_j)- \Theta(p_j,t_\varepsilon)\bigr)\\
\leq &\,1/2+\varepsilon/2+\sum_{i=1}^3\left\vert\int_{t_\varepsilon}^{t_j}
\left\langle\,\frac{P^i-\widehat{p}_j}{2(T-\xi)}\,\biggl\vert\,\tau^i(1,\xi)\right\rangle
\rho_{\widehat{p}_j}(P^i,\xi)\,d\xi\,\right\vert\\
\leq &\,1/2+\varepsilon/2+\sum_{i=1}^3\frac{1}{\sqrt{2\pi}}
\int\limits_{d_j^i/\sqrt{2(T-t_\varepsilon)}}^{+\infty}e^{-y^2/2}\,dy\\
\leq &\,1/2+\varepsilon/2+\frac{1}{\sqrt{2\pi}}
\int\limits_{0}^{+\infty}e^{-y^2/2}\,dy
+\sum_{i=2}^3\frac{1}{\sqrt{2\pi}}
\int\limits_{d_j^i/\sqrt{2(T-t_\varepsilon)}}^{+\infty}e^{-y^2/2}\,dy\\
\leq &\,1+\varepsilon/2
+ \frac{2}{\sqrt{2\pi}}\int\limits_{C/\sqrt{2(T-t_\varepsilon)}}^{+\infty}e^{-y^2/2}\,dy\\
\leq &\,1+\varepsilon\,.
\end{align*}
Since $p_j\in W^-$ and $t_j\geq t_\varepsilon$ imply 
$\Theta(p_j,t_j)\leq 1/2+\varepsilon$, we have a contradiction.
\end{proof}

We are now ready to exclude the last case of Lemma~\ref{casi}.

\begin{prop}\label{nosing1} There are no {\em Type I} singularities
  such that the third case of Lemma~\ref{casi} holds.
\end{prop}
\begin{proof}
After the previous discussion, we have three situations to consider.
\medskip

{\em The case $\widehat{\Theta}(p)=1$.}\\
The points $p_j$ and $p$ stay in the open set $W$, as
$\widehat{\Theta}(p)=1$, 
hence we can find a closed interval 
$I$ in $W$ containing $p$, so definitely the sequence $p_j$.\\
Recalling the definition of the function
$b(\cdot , t)$ in formula~\eqref{defb}, it is easy to see that, since
for every point $q\in W$ we have $\,\widehat{q}\not\in\{P^i\}$, it
follows that $b(\cdot, t)$ is continuous on $W$ for every
$t\in[0,T)$.\\
As the functions $\Theta(\cdot , t)+b(\cdot , t)$ converge
monotonically on the compact $I$ to the constant function $1$ when
$t\to T$, by the monotonicity formula~\eqref{eqmonfor}, by applying 
Dini's Theorem such convergence is uniform. This implies that
$$
\lim_{j\to\infty}\Theta(p_j,t_j)+b(p_j,t_j)=1
$$
and coming back to the integrated monotonicity
formula~\eqref{eqmonfor-int} we have
$$
\lim_{j\to\infty}\int_{t_j}^T\int_{{\TTT_t}}
\left\vert\,\underline{k}+\frac{(x-\widehat{p}_j)^{\perp}}{2(T-t)}\right\vert^2
\rho_{\widehat{p}_j}\,ds\,dt= 0\,.
$$
Passing to the rescaled version~\eqref{reseqmonfor-int} we get
$$
\lim_{j\to\infty}\int_{\tt_j}^{+\infty}\int_{\widetilde{\TTT}_{\widehat{p}_j,\tt}}
\vert\,\widetilde{\underline{k}}+x^{\perp}\vert^2\widetilde{\rho}\,d\sigma\,d\tt=0\,.
$$
Since we know that there exists an
  index $j_0\in\NN$ such that the 3--point is outside the ball
  $B_{\sqrt{8C}}$, by Lemma~\ref{deltabuono0} we have that 
$\widetilde{k}_\sigma$ and $\widetilde{k}_\tt$ for the family 
$\bigl\{\widetilde{\TTT}_{\widehat{p}_j,r}\,\bigl\vert\, j>j_0$ and 
$r\in[\tt_j,\tt_j+\delta_{\sqrt{2C}}]\bigl\}$ are uniformly bounded by
a constant $D_R$ in the ball $B_{\sqrt{2C}}$ (notice that this is true
because the constant $D_R$ is independent of $x_0$ in
Lemma~\ref{deltabuono0}). Hence, considering a value $\delta>0$
smaller than the constant $\delta_{\sqrt{2C}}$ with 
$\delta D_{\sqrt{2C}}<1/2$, possibly passing to a subsequence, 
there exist times $r_j\in[\tt_j,\tt_j+\delta]$ such that
\begin{equation}\label{kiloju}
\lim_{j\to\infty}\int_{\widetilde{\TTT}_{\widehat{p}_j,r_j}}
\vert\,\widetilde{\underline{k}}+x^{\perp}\vert^2\widetilde{\rho}\,d\sigma=0\,.
\end{equation}
By the choice of $\delta<\delta_{\sqrt{2C}}$ 
and Lemma~\ref{rescestim}, 
we can assume that the two sequences of rescaled triods
$\widetilde{\TTT}_{\widehat{p}_j,\tt_j}$ and
$\widetilde{\TTT}_{\widehat{p}_j,r_j}$  both 
converge in the $C^2$ topology
in the ball $B_{\sqrt{2C}}$. Moreover, by limit~\eqref{kiloju} and 
Lemma~\ref{lemomt}, the
sequence $\widetilde{\TTT}_{\widehat{p}_j,r_j}$ 
converges to a straight line passing through the origin in 
the ball $B_{\sqrt{2C}}$, hence with zero curvature.\\
The points
$\widetilde{G}_{\widehat{p}_j}(p_j,r_j)\in{\widetilde{\TTT}_{\widehat{p}_j,r_j}}$
all belong to the closed ball $B_{\sqrt{2C}}$, according to
Lemma~\ref{puntopalla}, then, as $\widetilde{k}^2(p_j,t_j)\geq 1$, we
have
$$
\widetilde{k}^2(p_j,r_j)\geq
\widetilde{k}^2(p_j,\tt_j)-(r_j-\tt_j)\vert\partial_\tt \widetilde{k}^2(p_j,\xi)\vert
\geq1-\delta D_{\sqrt{2C}}\geq 1/2
$$
for some $\xi\in(\tt_j,r_j)$.\\
Since the convergence of the rescaled triods is in the $C^2$ topology 
in $B_{\sqrt{2C}}$, such estimate from below is clearly in contradiction 
with the fact that the limit triod has zero curvature,
proving the thesis in this case.
\medskip

{\em The case $\widehat{\Theta}(p)=1/2$.}\\
We suppose that $\, \widehat{p}=P^1$ hence, by continuity,
$\, \widehat{p}_j\to P^1$.\\
By the rescaled version of the monotonicity formula we have
\begin{align*}
\frac{1}{\sqrt{2\pi}}\int_{\tt_j}^{+\infty}\int_{\widetilde{\TTT}_{\widehat{p}_j,\tt}}\vert
\,\widetilde{\underline{k}}+x^\perp\vert^2\widetilde{\rho}\,d\sigma\,d\tt
=&\,\,\Theta(p_j,t_j)-\widehat{\Theta}(p_j)\\
&\,+\int_{t_j}^{T}\sum_{i=1}^3
\left\langle\,\frac{P^i-\widehat{p}_j}{2(T-t)}\,\biggl\vert\,{\tau}^i(1,t)\right\rangle
\,\rho_{\widehat{p}_j}(P^i,t)\,dt\\
\leq&\,\,\Theta(p_j,t_j) - 1 +
\frac{1}{\sqrt{2\pi}}\sum_{i=1}^3\int\limits^{+\infty}_{{d_j^i}/\sqrt{2(T-t_j)}}e^{-y^2/2}\,dy
\end{align*}
where we applied Lemma~\ref{stimadib} and we substituted 
$\widehat{\Theta}(p_j)=1$.\\
Taking $\varepsilon>0$, since $p_j\to p\in W^-$ and 
$t_j\to T$, by Lemma~\ref{intornobuono} we have $\Theta(p_j,t_j)\leq
1+\varepsilon$ if $j\in\NN$ is large enough, hence
$$
\int_{\tt_j}^{+\infty}\int_{\widetilde{\TTT}_{\widehat{p}_j,\tt}}\vert
\,\widetilde{\underline{k}}+x^\perp\vert^2\widetilde{\rho}\,d\sigma\,d\tt\leq
\varepsilon +
\frac{1}{\sqrt{2\pi}}\sum_{i=1}^3\int\limits^{+\infty}_{{d_j^i}
/\sqrt{2(T-t_j)}}e^{-y^2/2}\,dy\,.
$$
Now, using the hypothesis that
$\lim_{j\to\infty}\frac{d_j^i}{\sqrt{2(T-t_j)}}=+\infty$, we obtain
$$
\lim_{j\to\infty}\int_{\tt_j}^{+\infty}\int_{\widetilde{\TTT}_{\widehat{p}_j,\tt}}\vert
\,\widetilde{\underline{k}}+x^\perp\vert^2\widetilde{\rho}\,d\sigma\,d\tt=0\,,
$$
then we conclude like in the previous case.
\medskip

{\em The case $\widehat{\Theta}(p)=3/2$.}\\
We consider the evolution of only two curves of the triods $\TTT_t$, one
of them containing all the points $p_j$, let us say the ones given by
the curves $\gamma^1$ and $\gamma^2$. Thus, we restrict the map $G$
to the set $\VVV$ composed of these two curves of $\TTT$ and we call $\VVV_t$
the evolving sets in $\R^2$. The monotonicity formula in this
situation reads
\begin{align}\label{veqmonfor}
\frac{d\,}{dt}\int_{{\VVV_t}} \rho_{x_0}(x,t)\,ds=
&\,-\int_{{\VVV_t}} \left\vert\,\underline{k}+\frac{(x-x_0)^{\perp}}{2(T-t)}\right\vert^2
\rho_{x_0}(x,t)\,ds\\
&\,+\sum_{i=1}^2\left\langle\,\frac{P^i-x_0}{2(T-t)}\,\biggl\vert\,\tau^i(1,t)\right\rangle
\rho_{x_0}(P^i,t)\nonumber\\
&\,-\left\langle\,\frac{O-x_0}{2(T-t)}\,\biggl\vert\,\tau^3(0,t)\right\rangle
\rho_{x_0}(O,t)\nonumber\\
&\,-\lambda^3(0,t)\rho_{x_0}(O,t)\nonumber\,.
\end{align}
Indeed, it can be obtained by adding the contributions of the two curves
$\gamma^1$ and $\gamma^2$ 
given by Lemma~\ref{monforuno}, recalling that
$\gamma^1(0,t)=\gamma^2(0,t)=G(O,t)$ and $\sum_{i=1}^3
\tau^i(0,t)=\sum_{i=1}^3 \lambda^i(0,t)=0$.

We define the functions $\Theta_\VVV(\cdot , t)$ analogous to
$\Theta(\cdot , t)$ but relative to the sets $\VVV_t$.\\
For every $q\in\VVV$, if $t_j$ is any sequence of
times going to $T$, the rescaled triods
$\widetilde{\TTT}_{\widehat{q},\tt_j}$ locally converge in $C^2\loc$ (up to a
subsequence) to a subset of the limit of the triods 
$\widetilde{\TTT}_{\widehat{q},\tt_j}$ which can  be composed of a
straight line or one or two halflines from the origin of
$\R^2$. Hence, it follows that $\widehat{\Theta}_\VVV(q)=\lim_{t\to
  T}\Theta_\VVV(q,t_j)=\widehat{\Theta}(q)$ 
if $q\not\in W^+$ and $\widehat{\Theta}_\VVV(q)=1$ if $\,
\widehat{q}=\widehat{O}$.

Rescaling everything as before, equation~\eqref{veqmonfor} becomes
\begin{align*}
\frac{1}{\sqrt{2\pi}}\int_{\tt}^{\tt^*}\int_{\widetilde{\VVV}_{\widehat{p},\xi}}\vert
\,\widetilde{\underline{k}}+x^\perp\vert^2\widetilde{\rho}\,d\sigma\,d\xi
=&\,\,\Theta_\VVV(p,t)-\Theta_\VVV(p,t^*)\\
&\,+\int_{t}^{t^*}\sum_{i=1}^2
\left\langle\,\frac{P^i-\widehat{p}}{2(T-\xi)}\,\biggl\vert\,{\tau}^i(1,\xi)\right\rangle
\,\rho_{\widehat{p}}(P^i,\xi)\,d\xi\\
&\,-\int_{t}^{t^*}
\left\langle\,\frac{O-\widehat{p}}{2(T-\xi)}\,\biggl\vert\,{\tau}^3(0,\xi)\right\rangle
\,\rho_{\widehat{p}}(O,\xi)\,d\xi\\
&\,-\int_{t}^{t^*}
\lambda^i(0,\xi)\,\rho_{\widehat{p}}(O,\xi)\,d\xi\,.
\end{align*}
Since when $t^*\to T$, $\Theta_\VVV(p,t^*)\to\widehat{\Theta}_\VVV(p)$ and the
integral
$$
\frac{1}{\sqrt{2\pi}}\int_{\tt}^{\tt^*}\int_{\widetilde{\VVV}_{\widehat{p},\xi}}\vert
\,\widetilde{\underline{k}}+x^\perp\vert^2\widetilde{\rho}\,d\sigma\,d\xi
$$
is non decreasing and bounded above by the analogous integral for
the complete triod, the function
\begin{align*}
b_\VVV(p,t)=&\,\int_{t}^{T}\sum_{i=1}^2
\left\langle\,\frac{P^i-\widehat{p}}{2(T-\xi)}\,\biggl\vert\,{\tau}^i(1,\xi)\right\rangle
\,\rho_{\widehat{p}}(P^i,\xi)\,d\xi\\
&\,-\int_{t}^{T}\left[
\left\langle\,\frac{O-\widehat{p}}{2(T-\xi)}\,\biggl\vert\,{\tau}^3(0,\xi)\right\rangle
+\lambda^i(0,\xi)\right]\,\rho_{\widehat{p}}(O,\xi)\,d\xi
\end{align*}
is well defined, moreover $\lim_{t\to T}b_\VVV(p,t)=0$.

Restricted to a closed interval around $O\in\VVV$ 
contained in the set $\VVV\cap(W\cup W^+)$, all the functions $b_\VVV(\cdot, t)$
are continuous, at least for $t$ close enough to $T$, otherwise we
could find a sequence of times $t_l\nearrow T$ and points $q_l\to q\in
\VVV\cap(W\cup W^+)$ such that 
$$
\frac{\vert G(O,t_l)-\widehat{q}_l\vert}{2(T-t_l)}\to 0
$$
which would imply
$$
\frac{\vert \widehat{O}-\widehat{q}_l\vert}{2(T-t_l)}\leq
\frac{\vert \widehat{O}-G(O,t_l)\vert}{2(T-t_l)}+
\frac{\vert G(O,t_l)-\widehat{q}_l\vert}{2(T-t_l)}\to\sqrt{2C}
$$
as $l\to\infty$, clearly in contradiction with
assumption~\eqref{assump999}.\\
Reasoning as before, the continuous functions
$\Theta_\VVV(\cdot , t)+b_\VVV(\cdot , t)\vert_{\VVV\cap(W\cup W^+)}$ 
converge monotonically to the function $\widehat{\Theta}_\VVV$, 
which is constantly equal to $1$ on a compact
interval around $p\in\VVV$. Hence, by Dini's Theorem, such 
convergence is uniform and
\begin{align*}
\frac{1}{\sqrt{2\pi}}\int_{\tt_j}^{\infty}\int_{\widetilde{\VVV}_{\widehat{p}_j,\tt}}\vert
\,\widetilde{\underline{k}}+x^\perp\vert^2\widetilde{\rho}\,d\sigma\,d\tt
=&\,\,\Theta_\VVV(p_j,t_j) + b_\VVV(p_j,t_j) -\widehat{\Theta}_\VVV(p_j)\\
=&\,\,\Theta_\VVV(p_j,t_j) + b_\VVV(p_j,t_j) -1
\end{align*}
implies
$$
\lim_{j\to\infty}\int_{\tt_j}^{+\infty}\int_{\widetilde{\VVV}_{\widehat{p}_j,\tt}}\vert
\,\widetilde{\underline{k}}+x^\perp\vert^2\widetilde{\rho}\,d\sigma\,d\tt=0\,.
$$
Then we conclude again with the same argument of the previous cases,
looking at the curvature of the limit set in the ball $B_{\sqrt{2C}}$.
\end{proof}

Collecting together Propositions~\ref{tripsing}, \ref{nosing0},
\ref{nosing1} we conclude this section with the following theorem. 

\begin{teo} {\em Type I} singularities cannot develop during the
  smooth flow $\TTT_t$ of an embedded triod 
in a bounded and strictly convex domain $\Omega\subset\R^2$ if the
lengths of the three curves are uniformly positively bounded from
below.
\end{teo}

\begin{rem}
All this section, if the results on the function $E$ discussed in
Section~\ref{disuellesec} (in particular Theorem~\ref{dlteo}) 
hold for some network of curves with many 3--points, can be
extended. The reason is that if we assume a uniform bound from below
on the lengths of all the curves of such a network, in every blow up
argument all but at most one of the 3--points go to infinity and
vanish in the limit (the same for the family of end points).\\
The situation changes dramatically without the control on the
lengths, for instance in Proposition~\ref{resclimit} one should
consider as possible limits also all the unbounded networks
composed of straight segments and curves of Abresch and Langer (see
Lemma~\ref{lemomt}) forming $120$ degrees at the 3--points.
\end{rem}

\section{Type II Singularities}
\label{typeIIsecn}

We suppose to be in the {\em Type II} singularity case, that is, 
$$
\limup_{t\to T} (T-t)\max_{{\TTT}_t} k^2 = +\infty\,.
$$
We employ Hamilton's trick to get an eternal solution to the mean
curvature flow. 

Let us choose a sequence of times $t_n\in[0,T-1/n]$ and points $p_n\in\TTT$
such that
\begin{equation}\label{rateII}
k^2(p_n,t_n)(T-1/n-t_n)=\max_{\genfrac{}{}{0pt}{}
{t\in[0,T-1/n]}{p\in\TTT}}k^2(p,t)(T-1/n-t)\,.
\end{equation}
As the curvature blows up as $t\to T$, we can choose a 
(not relabelled) non decreasing subsequence $t_n$ converging to $T$, 
moreover, the curvature squared $k^2(p_n,t_n)$ goes monotonically to
$+\infty$ and we have
$$
k^2(p_n,t_n)t_n\to+\infty\,,\qquad k^2(p_n,t_n)(T-1/n-t_n)\to+\infty\,.
$$
The first limit is obvious, about the second, by
assumption~\eqref{rateII}, for every $M>0$ there
exists $\overline{t}<T$ and $\overline{p}\in\TTT$ such that
$k^2(\overline{p},\overline{t})(T-\overline{t})>2M$, then choosing $n\in\NN$
large enough such that $1/n<T-\overline{t}$ and
$k^2(\overline{p},\overline{t})(T-\overline{t}-1/n)>M$ we get 
$$
k^2(p_n,t_n)(T-1/n-t_n)\geq
k^2(\overline{p},\overline{t})(T-\overline{t}-1/n)>M\,.
$$
Since $k^2(p_n,t_n)(T-1/n-t_n)$ is an increasing sequence and $M$ was
arbitrary,
$$
\lim_{n\to\infty} k^2(p_n,t_n)(T-1/n-t_n)=+\infty\,.
$$

We rescale now the triods as follows: let
$F_n:\TTT\times[-k^2(p_n,t_n)t_n, k^2(p_n,t_n)(T-1/n-t_n)]\to\R^2$ be
the evolution given by
$$
F_n(p,\tt)=k(p_n,t_n)[F(p,\tt/k^{2}(p_n,t_n)+t_n)-F(p_n,t_n)]
$$
and $\TTT^n_\tt=F_n(\TTT,\tt)$.
\begin{itemize}
\item $F_n(p_n,0)=0$, $k_n(p_n,0)=1$;
\item for every $\varepsilon>0$ and $\omega>0$ there exists
  $\overline{n}$ such that $\max_{\TTT^n_\tt} k^2\leq 1+\varepsilon$ for
  every $n\geq \overline{n}$ and $\tt\in[-k^2(p_n,t_n)t_n,\omega]$.
\end{itemize}
Indeed (the first point is immediate), by the choice of the pair 
$(p_n,t_n)$ we get
\begin{align*}
k^2_n(p,t)&\,= [k(p_n,t_n)]^{-2}k^2(p,\tt/k^{2}(p_n,t_n)+t_n)\\
&\,\leq [k(p_n,t_n)]^{-2}k^2(p_n,t_n)
\frac{T-1/n-t_n}{T-1/n-t_n-\tt/k^2(p_n,t_n)}\\
&\,= \frac{k^2(p_n,t_n)(T-1/n-t_n)}{k^2(p_n,t_n)(T-1/n-t_n)-\tt}\,,
\end{align*}
if $\tt\in[-k^2(p_n,t_n)t_n, k^2(p_n,t_n)(T-1/n-t_n)]$.\\
The claim follows as $k^2(p_n,t_n)(T-1/n-t_n)\to+\infty$.

\begin{prop}\label{typeIIblow}  The family of flows $F_n$ converges in
  the $C^2\loc$ topology to an evolution by 
  curvature of unbounded triods $\TTT^\infty_\tt$ 
  (or curves, or curves with a single end point) without
  self--intersections, in the time interval
  $(-\infty,+\infty)$. Such a flow is
  called {\em eternal}.\\
Moreover, the modulus of the curvature is uniformly bounded in space
and time and it takes its absolute maximum, which is 1, at time
$\tt=0$ at the origin of $\R^2$.
\end{prop}
\begin{proof}
By the previous discussion, on every bounded interval of time, the
evolutions $F_n$ have uniformly bounded curvature, hence, reasoning
like in the previous section (interior estimates of
Ecker and Huisken and treatment of the 3--point terms as in the proof of
Proposition~\ref{resclimit}), we have a uniform local bound also on
$k_\sigma$ for the family of rescaled flows.\\
Passing to the limit in the $C^2\loc$ topology we then get the flow
$\TTT^\infty_\tt$ which clearly satisfies the claimed properties on the
curvature, again by the computation above.\\
Finally, the limit flow is embedded by Lemma~\ref{elemma}.
\end{proof}

\begin{rem}
It can be shown that such a limit flow is actually smooth (see
Definition~\ref{smttt}), Indeed, we know by standard estimates that 
it is smooth far from the 3--point and, being the curves
$C^{2,\alpha}$ by the uniform local estimate of $k_\sigma$, 
we can apply (locally around the 3--point) a 
version of the small time existence and
uniqueness Theorem~\ref{smoothexist}, substituting the boundary
conditions regarding the fixed end points with 
the smooth motion of three chosen points of the triods
$\TTT^\infty_\tt$.
\end{rem}

In the case of the evolution of a closed curve in the plane, it is
possible to show that the limit flow arising from the analogous of
this proposition is a translating unbounded curve. This conclusion is
reached in two steps: first, one shows that at every time in
$(-\infty,+\infty)$ the limit is a convex curve, then, by means of the
Harnack estimate proved by Hamilton in~\cite{hamilton4}, it follows
that it is a translating flow by curvature.\\
In our situation convexity means that the curvature is never
zero. Unfortunately, the presence of the 3--point makes troublesome the
extension of these results to the evolution of triods and we can only state
the following two conjectures.

\begin{conge}\label{congconvex} If the curvature is zero at some point of
  the triod (or curve) $\TTT^\infty_\tt$, then it is zero along all
  the curve containing such a point (everywhere).
\end{conge}

\begin{conge}\label{congtrans} The triods (or curves) $\TTT^\infty_\tt$ move by
  translation.
\end{conge}

If this last conjecture is true, we can exclude {\em Type II} singularities,
indeed, it would follows that $\TTT_\tt^\infty$ is one of the sets of 
Lemmas~\ref{lemtra},~\ref{lemtra2} or~\ref{lemtra3}, hence, considering
a pair of points $p,q$ on two curves of the triods
with opposite convexity such that the segment $[p,q]$ is orthogonal to
the velocity vector $w\in\R^2$ and sending them both to infinity, we
can see that $E(\TTT^\infty_\tt)=0$, which is in contradiction with
Lemma~\ref{elemma} (where the function $E$ is defined). Indeed, 
the distance $[p,q]$ is bounded by a constant (the three curves
are pieces of the same {\em grim reaper}) and the area $A_{p,q}$
diverges.\\
The same argument works if the limit $\TTT^\infty_\tt$ is a single
curve, noticing that in this case such a curve cannot have an end point,
because otherwise the curvature would be zero there 
(such point comes from an end point of $\TTT_t$), 
hence the equation $k=\langle w\,\vert\,\nu\rangle$ would imply that 
the curve is a static halfline parallel to $w$, so with null curvature, contradicting
Proposition~\ref{typeIIblow} which says that the curvature cannot be
identically zero.

\begin{rem} The same argument, but with a different geometric
  quantity, is used by Huisken in~\cite{huisk2} to
  exclude {\em Type II} singularities during the motion of a single
  curve.
\end{rem}

\begin{prop} If Conjecture~\ref{congtrans} is true then 
{\em Type II} singularities cannot develop during the
  smooth flow $\TTT_t$ of an embedded triod 
in a bounded and strictly convex domain $\Omega\subset\R^2$.
\end{prop}

Hamilton's proof that a convex blow up of a {\em Type II} singularity 
(in the standard smooth case) is translating, 
is heavily based on the maximum principle which, as we said, is
difficult to be applied in the case of triods. 
So, it could happen that only Conjecture~\ref{congconvex}
can be proved and one could possibly exclude {\em Type II}
singularities without actually show that the limit flow
$\TTT_\tt^\infty$ is translating. For instance, if the curvature is
always or never zero on each curve of a smooth triod, 
then its three curves have asymptotic tangents (indeed, they are
all convex but some of them in the opposite way, by the fact that 
the sum of the curvatures at the 3--point is zero), hence in
order to apply the previous argument based on the function $E$, it would
be enough to show that two of these limit tangents, belonging to curves
with opposite convexity, coincide.

\begin{conge}\label{nnscong} The flow $\TTT_t$ of an embedded triod 
in a bounded and strictly convex domain $\Omega\subset\R^2$, 
such that the lengths of the three curves are 
uniformly positively bounded from
below, does not develop singularities at all.
\end{conge}

\section{Open Problems}
\label{conclusione}

The main problem left open in the paper is Conjecture~\ref{congtrans},
whose validity would imply Conjectures~\ref{nnscong}
and~\ref{congettura}, indeed, if the flow is defined for every
positive time, we can show the convergence of the triods to the
minimal connection between the three end points $P^i$.\\
As the total length of the triods decreases, we have the estimate 
$$
\int_0^{+\infty}\int_{\TTT_t}k^2\,ds\,dt\leq
\sum_{i=1}^3L^i(0)<+\infty\,,
$$
then, supposing that for a sequence of times $t_j\nearrow+\infty$ we have 
$\int_{\TTT_{t_j}}k^2\,ds\geq \delta$, for some
$\delta>0$, and the lengths
are uniformly positively bounded from below, by
inequality~\eqref{evolint999} we get that $\int_{\TTT_t} k^2\,ds\geq
\delta/2$ for every $t$ belonging to a uniform time interval around
every $t_j$. This is clearly in contradiction with the previous
estimate.\\
Hence, as $\lim_{t\to+\infty}\int_{\TTT_{t}}k^2\,ds=0$, by
Proposition~\ref{topolino5}, all the derivatives of
the curvature are uniformly bounded. By the usual convergence
arguments, this implies that (possibly after reparametrization) the
triods converge in the $C^\infty$ topology to a limit triod with zero
curvature, which is clearly the minimal connection between the three
end points (Steiner configuration).

Other questions in the paper that we would like to set are concerned
with the extension of the results to all the networks, possibly with
loops, with only 3--points (no 4--points or higher order points), in
particular, proving an analogous of Theorem~\ref{dlteo}. This would
also make superfluous the requirement that the
ambient set $\Omega\subset\R^2$ is strictly convex.

Finally, we conclude listing some, naturally arising, research
directions.
\begin{enumerate}
\item The problem of the representation/uniqueness of Brakke flows (smooth/with
  equality) of triods (or networks), discussed in Remark~\ref{oprob1}.

\item The problem of the existence/uniqueness 
  of a (smooth) Brakke flow with equality for an initial triod not
  satisfying the $120$ degrees condition at the 3--point (see the
  introduction and Remark~\ref{oprob1}).

\item The study of the singularities such that the curvature blows
  up but the lengths are not positively bounded from below. 
  Such analysis requires new estimates and the
  classification of homothetic and translating networks moving by
  curvature.

\item The ``definitory'' problem of the motion of the multi--points
  and the analysis of the collapsing situations with change of
  topology (see the introduction and the papers of
  De~Giorgi~\cite{degio7} and Caraballo~\cite{caraballo1}).

\end{enumerate}

\bibliographystyle{amsplain}
\bibliography{tripunto}

\providecommand{\bysame}{\leavevmode\hbox to3em{\hrulefill}\thinspace}
\providecommand{\MR}{\relax\ifhmode\unskip\space\fi MR }
\providecommand{\MRhref}[2]{%
  \href{http://www.ams.org/mathscinet-getitem?mr=#1}{#2}
}
\providecommand{\href}[2]{#2}
\begin{thebibliography}{10}

\bibitem{ablang1}
U.~Abresch and J.~Langer, \emph{The normalized curve shortening flow and
  homothetic solutions}, J. Diff. Geom. \textbf{23} (1986), no.~2, 175--196.

\bibitem{adams}
R.~Adams, \emph{Sobolev {S}paces}, Academic Press, New York, 1975.

\bibitem{altawa}
F.~J. Almgren, J.~E. Taylor, and L.~Wang, \emph{Curvature driven flows: a
  variational approach}, SIAM J. Cont. and Opt. \textbf{31} (1993), 387--438.

\bibitem{altgra}
S.~Altschuler and M.~Grayson, \emph{Shortening space curves and flow through
  singularities}, J. Diff. Geom. \textbf{35} (1992), 283--298.

\bibitem{altsch}
S.~J. Altschuler, \emph{Singularities of the curve shrinking flow for space
  curves}, J. Diff. Geom. \textbf{34} (1991), no.~2, 491--514.

\bibitem{angen2}
S.~Angenent, \emph{Parabolic equations for curves on surfaces: curves with
  {$p$}--integrable curvature}, Ann. of Math. (2) \textbf{132} (1990),
  451--483.

\bibitem{angen3}
\bysame, \emph{Parabolic equations for curves on surfaces: intersections, blow
  up and generalized solutions}, Ann. of Math. (2) \textbf{133} (1991),
  171--215.

\bibitem{aubin0}
T.~Aubin, \emph{Some {N}onlinear {P}roblems in {R}iemannian {G}eometry},
  Springer, 1998.

\bibitem{hermul}
J.~M. Ball, D.~Kinderlehrer, P.~Podio-Guidugli, and M.~Slemrod (eds.),
  \emph{Fundamental {C}ontributions to the {C}ontinuum {T}heory of {E}volving
  {P}hase {I}nterfaces in {S}olids}, Springer, Berlin, 1999.

\bibitem{brakke}
K.~A. Brakke, \emph{The {M}otion of a {S}urface by its {M}ean {C}urvature},
  Princeton University Press, Princeton, N.J., 1978.

\bibitem{bronsard}
L.~Bronsard and F.~Reitich, \emph{On three-phase boundary motion and the
  singular limit of a vector--valued {G}inzburg--{L}andau equation}, Arch. Rat.
  Mech. Anal. \textbf{124} (1993), no.~4, 355--379.

\bibitem{caraballo1}
D.~G. Caraballo, \emph{A variational scheme for the evolution of polycrystals
  by curvature}, Ph.D. thesis, Princeton University, 1996.

\bibitem{chzh}
K.-S. Chou and X.-P. Zhu, \emph{Shortening complete plane curves}, J. Diff.
  Geom. \textbf{50} (1998), no.~3, 471--504.

\bibitem{degio4}
E.~{De~Giorgi}, \emph{Barriers, boundaries and motion of manifolds}, Scuola
  Normale Superiore di Pisa, 1995.

\bibitem{degio7}
\bysame, \emph{Motions of partitions}, Variational methods for discontinuous
  structures (Como, 1994), Progr. Nonlinear Differential Equations Appl.,
  vol.~25, Birkh\"auser, Basel, 1996, pp.~1--5.

\bibitem{eck1}
K.~Ecker, \emph{Lectures on regularity for mean curvature flow}, preprint,
  2001.

\bibitem{eckhui2}
K.~Ecker and G.~Huisken, \emph{Interior estimates for hypersurfaces moving by
  mean curvature}, Invent. Math. \textbf{105} (1991), no.~3, 547--569.

\bibitem{es}
L.~C. Evans and J.~Spruck, \emph{Motion of level sets by mean curvature {I}},
  J. Diff. Geom. \textbf{33} (1991), 635--681.

\bibitem{gray1}
M.~A. Grayson, \emph{The heat equation shrinks embedded plane curves to round
  points}, J. Diff. Geom. \textbf{26} (1987), 285--314.

\bibitem{gurtin2}
M.~E. Gurtin, \emph{Thermomechanics of {E}volving {P}hase {B}oundaries in the
  {P}lane}, Oxford Science Publication, New York, 1993.

\bibitem{hamilton2}
R.~S. Hamilton, \emph{Four--manifolds with positive {R}icci curvature}, J.
  Diff. Geom. \textbf{24} (1986), no.~2, 153--179.

\bibitem{hamilton4}
\bysame, \emph{Harnack estimate for the mean curvature flow}, J. Diff. Geom.
  \textbf{41} (1995), no.~1, 215--226.

\bibitem{hamilton3}
\bysame, \emph{Isoperimetric estimates for the curve shrinking flow in the
  plane}, Modern methods in complex analysis (Princeton, NJ, 1992), Princeton
  University Press, Princeton, NJ, 1995, pp.~201--222.

\bibitem{huisk1}
G.~Huisken, \emph{Flow by mean curvature of convex surfaces into spheres}, J.
  Diff. Geom. \textbf{20} (1984), 237--266.

\bibitem{huisk3}
\bysame, \emph{Asymptotic behavior for singularities of the mean curvature
  flow}, J. Diff. Geom. \textbf{31} (1990), 285--299.

\bibitem{huisk2}
\bysame, \emph{A distance comparison principle for evolving curves}, Asian J.
  Math. \textbf{2} (1998), 127--133.

\bibitem{ilman1}
T.~Ilmanen, \emph{Elliptic regularization and partial regularity for motion by
  mean curvature}, Mem. Amer. Math. Soc., vol. 108(520), AMS, 1994.

\bibitem{ilman3}
\bysame, \emph{Singularities of mean curvature flow of surfaces},
  http://www.math.ethz.ch/{$\sim$}ilmanen/papers/sing.ps, 1995.

\bibitem{ilman4}
\bysame, \emph{Lectures on mean curvature flow and related equations},
  http://www.math.ethz.ch/{$\sim$}ilmanen/papers/notes.ps, 1998.

\bibitem{kinderliu}
D.~Kinderlehrer and C.~Liu, \emph{Evolution of grain boundaries}, Math. Models
  Methods Appl. Sci. \textbf{11} (2001), no.~4, 713--729.

\bibitem{langer2}
J.~Langer, \emph{A compactness theorem for surfaces with {$L_p$}--bounded
  second fundamental form}, Math. Ann. \textbf{270} (1985), 223--234.

\bibitem{lunardi1}
A.~Lunardi, \emph{Analytic semigroups and optimal regularity in parabolic
  problems}, Birkh\"auser Verlag, Basel, 1995.

\bibitem{lusiw}
A.~Lunardi, E.~Sinestrari, and W.~von Wahl, \emph{A semigroup approach to the
  time dependent parabolic initial--boundary value problem}, Differential
  Integral Equations \textbf{5} (1992), no.~6, 1275--1306.

\bibitem{simon}
L.~Simon, \emph{Lectures on {G}eometric {M}easure {T}heory}, Proc. Center Math.
  Anal., vol.~3, Australian National University, Canberra, 1983.

\bibitem{solonnikov1}
V.~A. Solonnikov, \emph{Boundary value problems of mathematical physics.
  {VIII}}, American Mathematical Society, Providence, R.I., 1975.

\bibitem{soner1}
H.~M. Soner, \emph{Motion of a set by the curvature of its boundary}, J. Diff.
  Eqs. \textbf{101} (1993), no.~2, 313--372.

\bibitem{stone1}
A.~Stone, \emph{A density function and the structure of singularities of the
  mean curvature flow}, Calc. Var. Partial Differential Equations \textbf{2}
  (1994), 443--480.

\bibitem{stone2}
\bysame, \emph{A boundary regularity theorem for mean curvature flow}, J. Diff.
  Geom. \textbf{44} (1996), 371--434.

\end{thebibliography}

\end{document}